\documentclass[aps,pra,letterpaper,notitlepage,longbibliography,reprint,onecolumn,superscriptaddress]{revtex4-1}
\usepackage[table]{xcolor}
\usepackage{graphicx,amsmath,amsfonts,amssymb,amsthm,bm,bbm}
\usepackage{tikz,pgfplots}
\pgfplotsset{compat=newest}
\pgfplotsset{plot coordinates/math parser=false}
\usetikzlibrary{arrows,decorations.pathmorphing,backgrounds,positioning,fit,calc,patterns,shapes,decorations.markings,shadings,snakes}
\usetikzlibrary{external}
\usepackage{ifthen}
\usepackage{multirow}
\usepackage{diagbox}
\usepackage{stmaryrd}
\usepackage{tabularx}
\usepackage{booktabs}
\newcommand{\greyhline}{\arrayrulecolor[gray]{.8}\hline\arrayrulecolor{black}}
\newcommand{\greycline}[1]{\arrayrulecolor[gray]{.8}\cline{#1}\arrayrulecolor{black}}

\newcommand{\cat}[1]{\ensuremath{\mathcal{#1}}}

\usepackage{pdftexcmds}

\makeatletter
\ifcase\pdf@shellescape
\or
\tikzset{external/system call={pdflatex \tikzexternalcheckshellescape -halt-on-error -interaction=batchmode -jobname "\image" "\texsource"; ps2pdf13 -dEmbedAllFonts=true -r100000  "\image".pdf "\image-13".pdf && cp "\image-13".pdf "\image".pdf && rm "\image-13".pdf && rm "\image".log && rm "\image".dpth && rm "\image"Notes.bib}}
\or
\fi
\makeatother

\newtheoremstyle{exampleA}% name of the style to be used
{\topsep}% measure of space to leave above the theorem. E.g.: 3pt
{10pt}% measure of space to leave below the theorem. E.g.: 3pt
{}% name of font to use in the body of the theorem
{}% measure of space to indent
{\bfseries}% name of head font
{\\\\}% punctuation between head and body
{ 10mm }% space after theorem head; " " = normal interword space
{\thmname{#1}\thmnumber{ #2}\thmnote{ (#3)}}

\theoremstyle{exampleA}

\usepackage[pdfpagelabels,pdftex,bookmarks,breaklinks]{hyperref}
\definecolor{darkblue}{RGB}{0,0,127} % choose colors
\definecolor{darkgreen}{RGB}{0,180,0}
\definecolor{darkred}{RGB}{180,0,0}
\hypersetup{
	colorlinks, 
	linkcolor=darkblue, 
	citecolor=darkgreen, 
	filecolor=red, 
	urlcolor=blue,
	pdftitle={Computing Defects Associated to Bounded Domain Wall Structures: The Z/pZ Case.}, 
	pdfauthor={Jacob C Bridgeman and Daniel Barter}
}

\definecolor{tctwistcolor}{RGB}{255,0,0}
\definecolor{tcmcolor}{RGB}{0,0,180}
\definecolor{tcecolor}{RGB}{0,180,0}
\definecolor{tcTTppcolor}{RGB}{184,134,11}
\definecolor{tcTTpmcolor}{RGB}{252,15,192}

\theoremstyle{plain}
\newtheorem{theorem}{Theorem}

\theoremstyle{definition}
\newtheorem{definition}[theorem]{Definition}

\newtheorem{example}[theorem]{Example}

\newtheorem{algorithm}[theorem]{Algorithm}

\newcommand{\ZZ}[1]{\mathbb{Z}/#1\mathbb{Z}}

\let\oldonlinecite\onlinecite
\renewcommand{\onlinecite}[1]{Ref.~[\oldonlinecite{#1}]}

\newcommand{\onlinecites}[1]{Refs.~[\oldonlinecite{#1}]}

\DeclareMathOperator{\ob}{ob}

\newcommand{\D}[1]{\mathcal{D}{} }

\newcommand{\restrict}[1]{\raise-.2ex\hbox{\ensuremath|}_{#1}}

\definecolor{tensorblue}{rgb}{0.8,0.8,1}
\definecolor{tensorred}{rgb}{1,0.5,0.5}
\definecolor{tensorgreen}{rgb}{0.6,1,0.6}
\definecolor{tensorpurp}{rgb}{1,0.5,1}

\tikzset{ten/.style={fill=tensorblue}}
\tikzset{tenred/.style={fill=tensorred}}
\tikzset{tengreen/.style={fill=tensorgreen}}
\tikzset{tenpurp/.style={fill=tensorpurp}}

\makeatletter
\newcommand{\vast}{\bBigg@{4}}
\newcommand{\Vast}{\bBigg@{9}}
\makeatother

\def\Put(#1,#2)#3{\leavevmode\makebox(0,0){\put(#1,#2){#3}}}

\makeatletter
\def\pgf@plot@curveto@handler@finish{%
  \ifpgf@plot@started%
    \pgfpathcurvebetweentimecontinue{0}{0.995}{\pgf@plot@curveto@first}{\pgf@plot@curveto@first@support}{\pgf@plot@curveto@second}{\pgf@plot@curveto@second}%
  \fi%
}
\makeatother

\newlength\figureheight 
\newlength\figurewidth 
  
\newcommand{\includeTikz}[2]
{
	\tikzifexternalizing
	{
		\includeTikzrm{#1}{#2}
	}
	{
		\IfFileExists{figures/#1.pdf}{
			\includegraphics{figures/#1}
		}
		{
			\includeTikzrm{#1}{#2}
		}
	}
}

\newcommand{\includeTikzrm}[2]{
	\tikzset{external/remake next}
	\tikzsetnextfilename{#1}
	#2
}
  
  \definecolor{nicegreena}{RGB}{1,115,16}
  \definecolor{nicegreenb}{RGB}{1,240,16}
  
  \definecolor{nicegreen}{RGB}{60,183,82}
  \definecolor{dandark}{HTML}{666666}
    \colorlet{ccred}{red!20}
    \colorlet{ccgreen}{green!50}
    \colorlet{ccblue}{blue!20}

  \tikzset{hexr/.style= {shape=regular polygon,regular polygon sides=6,minimum size=1cm, draw,inner sep=0,anchor=center,fill=red!50}}
  \tikzset{hexg/.style= {shape=regular polygon,regular polygon sides=6,minimum size=1cm, draw,inner sep=0,anchor=center,fill=green!50}}
  \tikzset{hexb/.style= {shape=regular polygon,regular polygon sides=6,minimum size=1cm, draw,inner sep=0,anchor=center,fill=blue!50}}
  
  \definecolor{tensorblue}{rgb}{0.8,0.8,1}
  \definecolor{tensorred}{rgb}{1,0.5,0.5}
  \definecolor{tensorpurp}{rgb}{1,0.5,1}
  
  \tikzset{nonesty/.style={fill=none,draw=none}}
  \tikzset{ten/.style={fill=tensorblue}}
  \tikzset{tenred/.style={fill=tensorred}}
  \tikzset{tengreen/.style={fill=green!50!black!50}}
  \tikzset{tenpurp/.style={fill=tensorpurp}}
  \tikzset{tengrey/.style={fill=black!20}}
  \tikzset{tenorange/.style={fill=orange!30}}
  \tikzset{u/.style={fill=blue!20,draw=black}}
  \tikzset{w/.style={fill=green!50!black!50,draw=black}}

\usetikzlibrary{calc}
\tikzstyle{inline text}=[text height=1.1ex, text depth=0.1ex, yshift=-.1ex]

\newcommand{\inflationalhs}[3]{
	\ifthenelse{\equal{#1}{}}{}
	{
		\draw[blue](-.75,1)--(0,.5);
		\node[above,inline text,blue] at (-.75,1) {#1};
	}
	\ifthenelse{\equal{#3}{}}{}
	{
		\draw[blue](.75,1)--(0,-.5);
		\node[above,inline text,blue] at (.75,1) {#3};
	}
	\draw[red,ultra thick] (0,-1)--(0,1);
	\node[below,inline text,red] at (0,-1) {#2};
}

\newcommand{\inflationarhs}[5]{
	\ifthenelse{\equal{#1}{}}{}
	{
		\draw[blue](-1.25,1)--(-.5,.5);
		\node[above,inline text,blue] at (-1.25,1) {#1};
	}
	\ifthenelse{\equal{#3}{}}{}
	{
		\draw[blue](-.5,-.75)--(.5,-.25);
		\node[above,inline text,blue] at (0,-.5) {#3};
	}
	\ifthenelse{\equal{#5}{}}{}
	{
		\draw[blue](1.25,1)--(.5,0);
		\node[above,inline text,blue] at (1.25,1) {#5};
	}
	\draw[red,ultra thick] (-.5,-1)--(-.5,1);
	\draw[orange,ultra thick] (.5,-1)--(.5,1);
	\node[below,inline text,red] at (-.5,-1) {#2};
	\node[below,inline text,orange] at (.5,-1) {#4};
}

\newcommand{\inflationarhss}[5]{
	\ifthenelse{\equal{#1}{}}{}
	{
		\draw[blue](-1.25,1)--(-.5,.5);
		\node[above,inline text,blue] at (-1.25,1) {#1};
	}
	\ifthenelse{\equal{#3}{}}{}
	{
		\draw[blue](-.5,-.75)--(.5,-.25);
		\node[above,inline text,blue] at (0,-.5) {#3};
	}
	\ifthenelse{\equal{#5}{}}{}
	{
		\draw[blue](1.25,1)--(.5,0);
		\node[above,inline text,blue] at (1.25,1) {#5};
	}
	\draw[red,ultra thick] (-.5,-1)--(-.5,1);
	\draw[red,ultra thick] (.5,-1)--(.5,1);
	\node[below,inline text,red] at (-.5,-1) {#2};
	\node[below,inline text,red] at (.5,-1) {#4};
}

\newcommand{\inflationblhs}[3]{
	\ifthenelse{\equal{#1}{}}{}
	{
		\draw[blue](-.75,-1)--(0,-.5);
		\node[below,inline text,blue] at (-.75,-1) {#1};
	}
	\ifthenelse{\equal{#3}{}}{}
	{
		\draw[blue](.75,-1)--(0,.5);
		\node[below,inline text,blue] at (.75,-1) {#3};
	}
	\draw[red,ultra thick] (0,-1)--(0,1);
	\node[above,inline text,red] at (0,1) {#2};
}

\newcommand{\inflationbrhs}[5]{
	\ifthenelse{\equal{#1}{}}{}
	{
		\draw[blue](-1.25,-1)--(-.5,-.5);
		\node[below,inline text,blue] at (-1.25,-1) {#1};
	}
	\ifthenelse{\equal{#3}{}}{}
	{
		\draw[blue](-.5,-.25)--(.5,.25);
		\node[above,inline text,blue] at (0,0) {#3};
	}
	\ifthenelse{\equal{#5}{}}{}
	{
		\draw[blue](1.25,-1)--(.5,.5);
		\node[below,inline text,blue] at (1.25,-1) {#5};
	}
	\draw[red,ultra thick] (-.5,-1)--(-.5,1);
	\draw[orange,ultra thick] (.5,-1)--(.5,1);
	\node[above,inline text,red] at (-.5,1) {#2};
	\node[above,inline text,orange] at (.5,1) {#4};
}

\newcommand{\inflationbrhss}[5]{
	\ifthenelse{\equal{#1}{}}{}
	{
		\draw[blue](-1.25,-1)--(-.5,-.5);
		\node[below,inline text,blue] at (-1.25,-1) {#1};
	}
	\ifthenelse{\equal{#3}{}}{}
	{
		\draw[blue](-.5,-.25)--(.5,.25);
		\node[above,inline text,blue] at (0,0) {#3};
	}
	\ifthenelse{\equal{#5}{}}{}
	{
		\draw[blue](1.25,-1)--(.5,.5);
		\node[below,inline text,blue] at (1.25,-1) {#5};
	}
	\draw[red,ultra thick] (-.5,-1)--(-.5,1);
	\draw[red,ultra thick] (.5,-1)--(.5,1);
	\node[above,inline text,red] at (-.5,1) {#2};
	\node[above,inline text,red] at (.5,1) {#4};
}

\newcommand{\annparamss}[4]
{
	\def\ta{#1};
	\def\tb{#2};
	\def\tap{#3};
	\def\tbp{#4};
}
\newcommand{\annss}[2]{
	\def\ra{.5};
	\def\rb{1.5};
	\ifthenelse{\equal{#1}{}}{}
		{
			\node[left,inline text,blue] at ($(0, 0) + (180:.75 and .75)$() {\footnotesize#1};
			\draw[blue,xscale=2.5] (0,-.75)to[out=90+45,in=270-45] (0,.75);
		}
	\ifthenelse{\equal{#2}{}}{}
	{
		\node[right,inline text,blue] at ($(0, 0) + (0:.75 and .75)$() {\footnotesize#2};
		\draw[blue,xscale=1.5] (0,-1.25)to[out=90-45,in=270+45] (0,1.25);
	}
	\draw[red,ultra thick] (0,\ra)--(0,\rb);
	\draw[red,ultra thick] (0,-\ra)--(0,-\rb);
	\node[below,inline text,red] at ($(0, 0) + (270:1.5 and 1.5)$() {\footnotesize\ta};
	\node[above,inline text,red] at ($(0, 0) + (90:1.5 and 1.5)$() {\footnotesize\tb};
	\node[above,inline text,red] at ($(0, 0) + (270:.5 and .5)$() {\footnotesize\tap};
	\node[below,inline text,red] at ($(0, 0) + (90:.5 and .5)$() {\footnotesize\tbp};
	\draw (0,0) circle (\ra);
	\draw (0,0) circle (\rb);
}
\newcommand{\annparamst}[4]
{
	\def\ta{#1};
	\def\tb{#2};
	\def\tap{#3};
	\def\tbp{#4};
}
\newcommand{\annst}[2]{
	\def\ra{.5};
	\def\rb{1.5};
	\ifthenelse{\equal{#1}{}}{}
	{
		\node[left,inline text,blue] at ($(0, 0) + (180:.75 and .75)$() {\footnotesize#1};
		\draw[blue,xscale=2.5] (0,-.75)to[out=90+45,in=270-45] (0,.75);
	}
	\ifthenelse{\equal{#2}{}}{}
	{
		\node[right,inline text,blue] at ($(0, 0) + (0:.75 and .75)$() {\footnotesize#2};
		\draw[blue,xscale=1.5] (0,-1.25)to[out=90-45,in=270+45] (0,1.25);
	}
	\draw[orange,ultra thick] (0,\ra)--(0,\rb);
	\draw[red,ultra thick] (0,-\ra)--(0,-\rb);
	\node[below,inline text,red] at ($(0, 0) + (270:1.5 and 1.5)$() {\footnotesize\ta};
	\node[above,inline text,orange] at ($(0, 0) + (90:1.5 and 1.5)$() {\footnotesize\tb};
	\node[above,inline text,red] at ($(0, 0) + (270:.5 and .5)$() {\footnotesize\tap};
	\node[below,inline text,orange] at ($(0, 0) + (90:.5 and .5)$() {\footnotesize\tbp};
	\draw (0,0) circle (\ra);
	\draw (0,0) circle (\rb);
}

\newcommand{\annstpq}[4]{
	\ifthenelse{\equal{#1}{}}{}
	{
		\node[left,inline text,blue] at ($(180:.6 and .75)$() {\footnotesize#1};
		\coordinate (A) at ({.5*cos(-140)},{1*sin(140)});
		\draw[blue] ({.5*cos(-140)},{1*sin(-140)}) to[out=90+45,in=270-45] (A);
	}
	\ifthenelse{\equal{#2}{}}{}
	{
		\node[above,inline text,blue] at ($(90:.5 and .9)$() {\footnotesize#2};
		\coordinate (B) at ({.5*cos(40)},{1.7*sin(40)});
		\draw[blue] ({.5*cos(140)},{1.2*sin(140)}) to[out=45,in=270-45] (B);
	}
	\ifthenelse{\equal{#3}{}}{}
	{
		\node[right,inline text,blue] at ($(0,.3)+(0:.75 and .75)$() {\footnotesize#3};
		\coordinate (C) at ({.5*cos(40)},{2*sin(40)});
		\draw[blue] ({.5*cos(40)},{1*sin(-40)}) to[out=90-45,in=270+45] (C);
	}
	\ifthenelse{\equal{#4}{}}{}
	{
		\node[above,inline text,blue] at ($(270:.5 and 1.1)$() {\footnotesize#4};
		\coordinate (A) at ({.5*cos(-40)},{1.5*sin(-40)});
		\draw[blue] ({.5*cos(-140)},{2*sin(-140)}) to[out=45,in=270-45] (A);
	}
	\draw[red,ultra thick] ({.5*cos(-140)},{.5*sin(-140)+.1})--({.5*cos(-140)},{-1.5*sin(acos(1/3*cos(-140)))});
	\draw[orange,ultra thick] ({.5*cos(-40)},{.5*sin(-40)+.1})--({.5*cos(-40)},{-1.5*sin(acos(1/3*cos(-40)))});
	\draw[darkgreen,ultra thick] ({.5*cos(140)},{.5*sin(140)-.1})--({.5*cos(140)},{1.5*sin(acos(1/3*cos(140)))});
	\draw[darkred,ultra thick] ({.5*cos(40)},{.5*sin(40)-.1})--({.5*cos(40)},{1.5*sin(acos(1/3*cos(40)))});
	\draw[fill=white] (0,0) ellipse (.5 and .5);
	\draw (0,0) circle (1.5 and 1.5);
	\node[below,inline text,red] at ({.5*cos(-140)},{-1.5*sin(acos(1/3*cos(-140)))}) {\footnotesize\ta};
	\node[below,inline text,orange] at ({.5*cos(-40)},{-1.5*sin(acos(1/3*cos(-40)))}) {\footnotesize\tb};
	\node[above,inline text,darkgreen] at ({.5*cos(140)},{1.5*sin(acos(1/3*cos(140)))}) {\footnotesize\tc};
	\node[above,inline text,darkred] at ({.5*cos(40)},{1.5*sin(acos(1/3*cos(40)))}) {\footnotesize\td};
}
\newcommand{\annsspq}[4]{
	\ifthenelse{\equal{#1}{}}{}
	{
		\node[left,inline text,blue] at ($(180:.6 and .75)$() {\footnotesize#1};
		\coordinate (A) at ({.5*cos(-140)},{1*sin(140)});
		\draw[blue] ({.5*cos(-140)},{1*sin(-140)}) to[out=90+45,in=270-45] (A);
	}
	\ifthenelse{\equal{#2}{}}{}
	{
		\node[above,inline text,blue] at ($(90:.5 and .9)$() {\footnotesize#2};
		\coordinate (B) at ({.5*cos(40)},{1.7*sin(40)});
		\draw[blue] ({.5*cos(140)},{1.2*sin(140)}) to[out=45,in=270-45] (B);
	}
	\ifthenelse{\equal{#3}{}}{}
	{
		\node[right,inline text,blue] at ($(0,.3)+(0:.75 and .75)$() {\footnotesize#3};
		\coordinate (C) at ({.5*cos(40)},{2*sin(40)});
		\draw[blue] ({.5*cos(40)},{1*sin(-40)}) to[out=90-45,in=270+45] (C);
	}
	\ifthenelse{\equal{#4}{}}{}
	{
		\node[above,inline text,blue] at ($(270:.5 and 1.1)$() {\footnotesize#4};
		\coordinate (A) at ({.5*cos(-40)},{1.5*sin(-40)});
		\draw[blue] ({.5*cos(-140)},{2*sin(-140)}) to[out=45,in=270-45] (A);
	}
	\draw[red,ultra thick] ({.5*cos(-140)},{.5*sin(-140)+.1})--({.5*cos(-140)},{-1.5*sin(acos(1/3*cos(-140)))});
	\draw[red,ultra thick] ({.5*cos(-40)},{.5*sin(-40)+.1})--({.5*cos(-40)},{-1.5*sin(acos(1/3*cos(-40)))});
	\draw[darkgreen,ultra thick] ({.5*cos(140)},{.5*sin(140)-.1})--({.5*cos(140)},{1.5*sin(acos(1/3*cos(140)))});
	\draw[darkred,ultra thick] ({.5*cos(40)},{.5*sin(40)-.1})--({.5*cos(40)},{1.5*sin(acos(1/3*cos(40)))});
	\draw[fill=white] (0,0) ellipse (.5 and .5);
	\draw (0,0) circle (1.5 and 1.5);
	\node[below,inline text,red] at ({.5*cos(-140)},{-1.5*sin(acos(1/3*cos(-140)))}) {\footnotesize\ta};
	\node[below,inline text,red] at ({.5*cos(-40)},{-1.5*sin(acos(1/3*cos(-40)))}) {\footnotesize\tb};
	\node[above,inline text,darkgreen] at ({.5*cos(140)},{1.5*sin(acos(1/3*cos(140)))}) {\footnotesize\tc};
	\node[above,inline text,darkred] at ({.5*cos(40)},{1.5*sin(acos(1/3*cos(40)))}) {\footnotesize\td};
}
\newcommand{\annstpp}[4]{
	\ifthenelse{\equal{#1}{}}{}
	{
		\node[left,inline text,blue] at ($(180:.6 and .75)$() {\footnotesize#1};
		\coordinate (A) at ({.5*cos(-140)},{1*sin(140)});
		\draw[blue] ({.5*cos(-140)},{1*sin(-140)}) to[out=90+45,in=270-45] (A);
	}
	\ifthenelse{\equal{#2}{}}{}
	{
		\node[above,inline text,blue] at ($(90:.5 and .9)$() {\footnotesize#2};
		\coordinate (B) at ({.5*cos(40)},{1.7*sin(40)});
		\draw[blue] ({.5*cos(140)},{1.2*sin(140)}) to[out=45,in=270-45] (B);
	}
	\ifthenelse{\equal{#3}{}}{}
	{
		\node[right,inline text,blue] at ($(0,.3)+(0:.75 and .75)$() {\footnotesize#3};
		\coordinate (C) at ({.5*cos(40)},{2*sin(40)});
		\draw[blue] ({.5*cos(40)},{1*sin(-40)}) to[out=90-45,in=270+45] (C);
	}
	\ifthenelse{\equal{#4}{}}{}
	{
		\node[above,inline text,blue] at ($(270:.5 and 1.1)$() {\footnotesize#4};
		\coordinate (A) at ({.5*cos(-40)},{1.5*sin(-40)});
		\draw[blue] ({.5*cos(-140)},{2*sin(-140)}) to[out=45,in=270-45] (A);
	}
	\draw[red,ultra thick] ({.5*cos(-140)},{.5*sin(-140)+.1})--({.5*cos(-140)},{-1.5*sin(acos(1/3*cos(-140)))});
	\draw[orange,ultra thick] ({.5*cos(-40)},{.5*sin(-40)+.1})--({.5*cos(-40)},{-1.5*sin(acos(1/3*cos(-40)))});
	\draw[darkgreen,ultra thick] ({.5*cos(140)},{.5*sin(140)-.1})--({.5*cos(140)},{1.5*sin(acos(1/3*cos(140)))});
	\draw[darkgreen,ultra thick] ({.5*cos(40)},{.5*sin(40)-.1})--({.5*cos(40)},{1.5*sin(acos(1/3*cos(40)))});
	\draw[fill=white] (0,0) ellipse (.5 and .5);
	\draw (0,0) circle (1.5 and 1.5);
	\node[below,inline text,red] at ({.5*cos(-140)},{-1.5*sin(acos(1/3*cos(-140)))}) {\footnotesize\ta};
	\node[below,inline text,orange] at ({.5*cos(-40)},{-1.5*sin(acos(1/3*cos(-40)))}) {\footnotesize\tb};
	\node[above,inline text,darkgreen] at ({.5*cos(140)},{1.5*sin(acos(1/3*cos(140)))}) {\footnotesize\tc};
	\node[above,inline text,darkgreen] at ({.5*cos(40)},{1.5*sin(acos(1/3*cos(40)))}) {\footnotesize\td};
}
\newcommand{\annsspp}[4]{
	\ifthenelse{\equal{#1}{}}{}
	{
		\node[left,inline text,blue] at ($(180:.6 and .75)$() {\footnotesize#1};
		\coordinate (A) at ({.5*cos(-140)},{1*sin(140)});
		\draw[blue] ({.5*cos(-140)},{1*sin(-140)}) to[out=90+45,in=270-45] (A);
	}
	\ifthenelse{\equal{#2}{}}{}
	{
		\node[above,inline text,blue] at ($(90:.5 and .9)$() {\footnotesize#2};
		\coordinate (B) at ({.5*cos(40)},{1.7*sin(40)});
		\draw[blue] ({.5*cos(140)},{1.2*sin(140)}) to[out=45,in=270-45] (B);
	}
	\ifthenelse{\equal{#3}{}}{}
	{
		\node[right,inline text,blue] at ($(0,.3)+(0:.75 and .75)$() {\footnotesize#3};
		\coordinate (C) at ({.5*cos(40)},{2*sin(40)});
		\draw[blue] ({.5*cos(40)},{1*sin(-40)}) to[out=90-45,in=270+45] (C);
	}
	\ifthenelse{\equal{#4}{}}{}
	{
		\node[above,inline text,blue] at ($(270:.5 and 1.1)$() {\footnotesize#4};
		\coordinate (A) at ({.5*cos(-40)},{1.5*sin(-40)});
		\draw[blue] ({.5*cos(-140)},{2*sin(-140)}) to[out=45,in=270-45] (A);
	}
	\draw[red,ultra thick] ({.5*cos(-140)},{.5*sin(-140)+.1})--({.5*cos(-140)},{-1.5*sin(acos(1/3*cos(-140)))});
	\draw[red,ultra thick] ({.5*cos(-40)},{.5*sin(-40)+.1})--({.5*cos(-40)},{-1.5*sin(acos(1/3*cos(-40)))});
	\draw[darkgreen,ultra thick] ({.5*cos(140)},{.5*sin(140)-.1})--({.5*cos(140)},{1.5*sin(acos(1/3*cos(140)))});
	\draw[darkgreen,ultra thick] ({.5*cos(40)},{.5*sin(40)-.1})--({.5*cos(40)},{1.5*sin(acos(1/3*cos(40)))});
	\draw[fill=white] (0,0) ellipse (.5 and .5);
	\draw (0,0) circle (1.5 and 1.5);
	\node[below,inline text,red] at ({.5*cos(-140)},{-1.5*sin(acos(1/3*cos(-140)))}) {\footnotesize\ta};
	\node[below,inline text,red] at ({.5*cos(-40)},{-1.5*sin(acos(1/3*cos(-40)))}) {\footnotesize\tb};
	\node[above,inline text,darkgreen] at ({.5*cos(140)},{1.5*sin(acos(1/3*cos(140)))}) {\footnotesize\tc};
	\node[above,inline text,darkgreen] at ({.5*cos(40)},{1.5*sin(acos(1/3*cos(40)))}) {\footnotesize\td};
}
\newcommand{\annssss}[4]{
	\ifthenelse{\equal{#1}{}}{}
	{
		\node[left,inline text,blue] at ($(180:.6 and .75)$() {\footnotesize#1};
		\coordinate (A) at ({.5*cos(-140)},{1*sin(140)});
		\draw[blue] ({.5*cos(-140)},{1*sin(-140)}) to[out=90+45,in=270-45] (A);
	}
	\ifthenelse{\equal{#2}{}}{}
	{
		\node[above,inline text,blue] at ($(90:.5 and .9)$() {\footnotesize#2};
		\coordinate (B) at ({.5*cos(40)},{1.7*sin(40)});
		\draw[blue] ({.5*cos(140)},{1.2*sin(140)}) to[out=45,in=270-45] (B);
	}
	\ifthenelse{\equal{#3}{}}{}
	{
		\node[right,inline text,blue] at ($(0,.3)+(0:.75 and .75)$() {\footnotesize#3};
		\coordinate (C) at ({.5*cos(40)},{2*sin(40)});
		\draw[blue] ({.5*cos(40)},{1*sin(-40)}) to[out=90-45,in=270+45] (C);
	}
	\ifthenelse{\equal{#4}{}}{}
	{
		\node[above,inline text,blue] at ($(270:.5 and 1.1)$() {\footnotesize#4};
		\coordinate (A) at ({.5*cos(-40)},{1.5*sin(-40)});
		\draw[blue] ({.5*cos(-140)},{2*sin(-140)}) to[out=45,in=270-45] (A);
	}
	\draw[red,ultra thick] ({.5*cos(-140)},{.5*sin(-140)+.1})--({.5*cos(-140)},{-1.5*sin(acos(1/3*cos(-140)))});
	\draw[red,ultra thick] ({.5*cos(-40)},{.5*sin(-40)+.1})--({.5*cos(-40)},{-1.5*sin(acos(1/3*cos(-40)))});
	\draw[red,ultra thick] ({.5*cos(140)},{.5*sin(140)-.1})--({.5*cos(140)},{1.5*sin(acos(1/3*cos(140)))});
	\draw[red,ultra thick] ({.5*cos(40)},{.5*sin(40)-.1})--({.5*cos(40)},{1.5*sin(acos(1/3*cos(40)))});
	\draw[fill=white] (0,0) ellipse (.5 and .5);
	\draw (0,0) circle (1.5 and 1.5);
	\node[below,inline text,red] at ({.5*cos(-140)},{-1.5*sin(acos(1/3*cos(-140)))}) {\footnotesize\ta};
	\node[below,inline text,red] at ({.5*cos(-40)},{-1.5*sin(acos(1/3*cos(-40)))}) {\footnotesize\tb};
	\node[above,inline text,red] at ({.5*cos(140)},{1.5*sin(acos(1/3*cos(140)))}) {\footnotesize\tc};
	\node[above,inline text,red] at ({.5*cos(40)},{1.5*sin(acos(1/3*cos(40)))}) {\footnotesize\td};
}

\newcommand{\generalpantsstpq}[8]
{
	\pgfmathsetmacro{\ra}{.25};
	\pgfmathsetmacro{\sa}{.75};
	\pgfmathsetmacro{\rb}{2};
	\pgfmathsetmacro{\rc}{1.5};
	\coordinate (A) at ({-\sa},{-\rc*sin(acos(-\sa/\rb))});
	\coordinate (B) at ({\sa},{-\rc*sin(acos(-\sa/\rb))});
	\coordinate (C) at ({-\sa},{\rc*sin(acos(\sa/\rb))});
	\coordinate (D) at ({\sa},{\rc*sin(acos(\sa/\rb))});
		\ifthenelse{\equal{#1}{}}{}
		{
			\coordinate (A1) at ($(A)!.9!(-\sa,-\ra)$);
			\coordinate (B1) at ({-\sa-1.5*\ra},{0});
			\coordinate (C1) at ($(C)!.9!(-\sa,\ra)$);
			\node[left,inline text,blue] at (A1) {\footnotesize#1};
			\draw[blue] (A1) to[out=90+45,in=270] (B1) to[out=90,in=180+45] (C1);
		}
		\ifthenelse{\equal{#2}{}}{}
		{
			\coordinate (A1) at ($(A)!.6!(-\sa,-\ra)$);
			\coordinate (B1) at ({-\sa+2*\ra},{0});
			\coordinate (C1) at ($(C)!.75!(-\sa,\ra)$);
			\node[right,inline text,blue] at (A1) {\footnotesize#2};
			\draw[blue] (A1) to[out=90-45,in=270] (B1) to[out=90,in=0-45] (C1);
		}
		\ifthenelse{\equal{#3}{}}{}
		{
			\coordinate (A1) at ($(B)!.9!(\sa,-\ra)$);
			\coordinate (B1) at ({\sa-2*\ra},{0});
			\coordinate (C1) at ($(D)!.9!(\sa,\ra)$);
			\node[left,inline text,blue] at (A1) {\footnotesize#3};
			\draw[blue] (A1) to[out=90+45,in=270] (B1) to[out=90,in=180+45] (C1);
		}
		\ifthenelse{\equal{#4}{}}{}
		{
			\coordinate (A1) at ($(B)!.75!(\sa,-\ra)$);
			\coordinate (B1) at ({\sa+1.5*\ra},{0});
			\coordinate (C1) at ($(D)!.75!(\sa,\ra)$);
			\node[right,inline text,blue] at (A1) {\footnotesize#4};
			\draw[blue] (A1) to[out=90-45,in=270] (B1) to[out=90,in=0-45] (C1);
		}
		\ifthenelse{\equal{#5}{}}{}
		{
			\coordinate (A1) at ($(A)!.3!(-\sa,-\ra)$);
			\coordinate (B1) at ({-\sa-3*\ra},{0});
			\coordinate (C1) at ($(C)!.5!(-\sa,\ra)$);
			\node[left,inline text,blue] at (A1) {\footnotesize#5};
			\draw[blue] (A1) to[out=90+45,in=270] (B1) to[out=90,in=180+45] (C1);
		}
		\ifthenelse{\equal{#6}{}}{}
		{
			\coordinate (A1) at ($(C)!.25!(-\sa,-\ra)$);
			\coordinate (C1) at ($(D)!.25!(\sa,\ra)$);
			\coordinate (B1) at ($(A1)!.5!(C1)$);
			\node[below,inline text,blue] at (B1) {\footnotesize#6};
			\draw[blue] (A1) to[out=45,in=180+45] (C1);
		}
		\ifthenelse{\equal{#7}{}}{}
		{
			\coordinate (A1) at ($(B)!.3!(\sa,-\ra)$);
			\coordinate (B1) at ({\sa+3*\ra},{0});
			\coordinate (C1) at ($(D)!.1!(\sa,\ra)$);
			\node[right,inline text,blue] at (A1) {\footnotesize#7};
			\draw[blue] (A1) to[out=90-45,in=270] (B1) to[out=90,in=-45] (C1);
		}
		\ifthenelse{\equal{#8}{}}{}
		{
			\coordinate (A1) at ($(A)!.15!(-\sa,-\ra)$);
			\coordinate (C1) at ($(B)!.15!(\sa,\ra)$);
			\coordinate (B1) at ($(A1)!.5!(C1)$);
			\node[above,inline text,blue] at (B1) {\footnotesize#8};
			\draw[blue] (A1) to[out=45,in=180+45] (C1);
		}
\begin{scope}[even odd rule]
	\clip (0,0) ellipse [x radius=\rb,y radius=\rc] (-\sa,0) circle (\ra) (\sa,0) circle (\ra);
	\draw [red,ultra thick] ({-\sa},{0})--({-\sa},{-2*\rb});
	\draw [orange,ultra thick] ({\sa},{0})--({\sa},{-2*\rb});
	\draw [darkgreen,ultra thick] ({-\sa},{0})--({-\sa},{2*\rb});
	\draw [darkred,ultra thick] ({\sa},{0})--({\sa},{2*\rb});
\end{scope}
	\draw (-\sa,0) circle (\ra);
	\draw (\sa,0) circle (\ra);
	\draw (0,0) ellipse [x radius=\rb,y radius=\rc];
	\node[red,inline text,below] at (A) {\footnotesize\ta};
	\node[orange,inline text,below] at (B) {\footnotesize\tb};
	\node[darkgreen,inline text,above] at (C) {\footnotesize\tc};
	\node[darkred,inline text,above] at (D) {\footnotesize\td};
}

\newcommand{\tube}[3]{
	\ifthenelse{\equal{#2}{}}{}
	{
		\draw[domain=-90:0,smooth,variable=\x,blue] plot ({cos(\x)},{.1*sin(\x)+(\x)/360-.25});
		\draw[domain=0:180,smooth,variable=\x,blue,dashed] plot ({cos(\x)},{.1*sin(\x)+(\x)/360-.25});
		\draw[domain=180:270,smooth,variable=\x,blue] plot ({cos(\x)},{.1*sin(\x)+(\x)/360-.25});
	}
	\node[below,inline text,red] at ($(0, -1) + (270:1 and .1)$() {\footnotesize#1};
	\node[above,inline text,red] at ($(0, 1) + (-270:1 and .1)$() {\footnotesize#3};
	\draw (0,1) ellipse (1 and .1);
	\draw (-1,1)--(-1,-1) (1,1)--(1,-1);
	\draw[dashed] ($(0, -1) + (180:1 and .1)$() arc (180:0:1 and .1);
	\draw ($(0, -1) + (180:1 and .1)$() arc (180:360:1 and .1);
	\ifthenelse{\equal{#2}{}}{}
	{
		\node[inline text,blue] at (.5,-.3) {\footnotesize#2};
	}
	\draw [red,thick] ($(0, -1) + (270:1 and .1)$()--($(0, 1) + (-90:1 and .1)$();
}

\newcommand{\ttube}[6]{
	\draw [red,thick,dashed] ($(0, -1) + (70:1 and .1)$()--($(0, 1) + (70:1 and .1)$();
	\ifthenelse{\equal{#4}{}}{}
	{
		\draw[domain=-110:0,smooth,variable=\x,blue] plot ({cos(\x)},{.1*sin(\x)+(\x)/360-.25});
		\draw[domain=0:70,smooth,variable=\x,blue,dashed] plot ({cos(\x)},{.1*sin(\x)+(\x)/360-.25});
	}
	\ifthenelse{\equal{#3}{}}{}
	{
		\draw[domain=70:180,smooth,variable=\x,blue,dashed] plot ({cos(\x)},{.1*sin(\x)+(\x)/360});
		\draw[domain=180:250,smooth,variable=\x,blue] plot ({cos(\x)},{.1*sin(\x)+(\x)/360});
	}
	\node[below,inline text,red] at ($(0, -1) + (250:1 and .1)$() {\footnotesize#1};
	\node[below,inline text,red] at ($(0, -1) + (-70:1 and .1)$() {\footnotesize#2};
	\node[above,inline text,red] at ($(0, 1) + (-250:1 and .1)$() {\footnotesize#5};
	\node[above,inline text,red] at ($(0, 1) + (70:1 and .1)$() {\footnotesize#6};
	\draw (0,1) ellipse (1 and .1);
	\draw (-1,1)--(-1,-1) (1,1)--(1,-1);
	\draw[dashed] ($(0, -1) + (180:1 and .1)$() arc (180:0:1 and .1);
	\draw ($(0, -1) + (180:1 and .1)$() arc (180:360:1 and .1);
	\ifthenelse{\equal{#4}{}}{}
	{
		\node[right,inline text,blue] at ({cos(0)},{.1*sin(0)+(0)/360-.25}) {\footnotesize#4};
	}
	\ifthenelse{\equal{#3}{}}{}
	{
		\node[left,inline text,blue] at ({cos(180)},{.1*sin(180)+(180)/360}) {\footnotesize#3};
	}
	\draw [red,thick] ($(0, -1) + (250:1 and .1)$()--($(0, 1) + (250:1 and .1)$();
}

\newcommand{\ttttube}[4]{
	\draw [red,thick,dashed] ($(0, -1) + (100:1 and .1)$()--($(0, 1) + (100:1 and .1)$();
	\draw [red,thick,dashed] ($(0, -1) + (40:1 and .1)$()--($(0, 1) + (40:1 and .1)$();
	\draw (0,1) ellipse (1 and .1);
	\draw (-1,1)--(-1,-1) (1,1)--(1,-1);
	\draw[dashed] ($(0, -1) + (180:1 and .1)$() arc (180:0:1 and .1);
	\draw ($(0, -1) + (180:1 and .1)$() arc (180:360:1 and .1);
	\ifthenelse{\equal{#1}{}}{}
	{
		\draw[domain=230-360:280-360,smooth,variable=\x,blue] plot ({cos(\x)},{.1*sin(\x)+(\x)/360-.25});
	}
	\ifthenelse{\equal{#2}{}}{}
	{
		\draw[domain=280-360:0,smooth,variable=\x,blue] plot ({cos(\x)},{.1*sin(\x)+(\x)/360-.15});
		\draw[domain=0:40,smooth,variable=\x,blue,dashed] plot ({cos(\x)},{.1*sin(\x)+(\x)/360-.15});
	}
	\ifthenelse{\equal{#3}{}}{}
	{
		\draw[domain=40:100,smooth,variable=\x,blue,dashed] plot ({cos(\x)},{.1*sin(\x)+(\x)/360-.05});
	}
	\ifthenelse{\equal{#4}{}}{}
	{
		\draw[domain=100:180,smooth,variable=\x,blue,dashed] plot ({cos(\x)},{.1*sin(\x)+(\x)/360+.05});
		\draw[domain=180:230,smooth,variable=\x,blue] plot ({cos(\x)},{.1*sin(\x)+(\x)/360+.05});
	}
	\draw [red,thick] ($(0, -1) + (230:1 and .1)$()--($(0, 1) + (230:1 and .1)$();
	\draw [red,thick] ($(0, -1) + (280:1 and .1)$()--($(0, 1) + (280:1 and .1)$();
	\node[below,inline text,red] at ($(0, -1) + (230:1 and .1)$() {\footnotesize\ta};
	\node[below,inline text,red] at ($(0, -1) + (280:1 and .1)$() {\footnotesize\tb};
	\node[below,inline text,red] at ($(0, -1) + (-100:1 and .1)$() {\footnotesize\tc};
	\node[below,inline text,red] at ($(0, -1) + (-40:1 and .1)$() {\footnotesize\td};
	\node[above,inline text,red] at ($(0, 1) + (-230:1 and .1)$() {\footnotesize\tap};
	\node[above,inline text,red] at ($(0, 1) + (-280:1 and .1)$() {\footnotesize\tbp};
	\node[above,inline text,red] at ($(0, 1) + (100:1 and .1)$() {\footnotesize\tcp};
	\node[above,inline text,red] at ($(0, 1) + (40:1 and .1)$() {\footnotesize\tdp};
	\ifthenelse{\equal{#1}{}}{}
	{
		\node[above,inline text,blue] at ({cos(-90)},{.1*sin(-90)+(-90)/360-.25}) {\footnotesize#1};
	}
	\ifthenelse{\equal{#2}{}}{}
	{
		\node[right,inline text,blue] at ({cos(0)},{.1*sin(0)+(0)/360-.15}) {\footnotesize#2};
	}
	\ifthenelse{\equal{#3}{}}{}
	{
		\node[above,inline text,blue] at ({cos(60)},{.1*sin(60)+(60)/360-.05}) {\footnotesize#3};
	}
	\ifthenelse{\equal{#4}{}}{}
	{
		\node[left,inline text,blue] at ({cos(180)},{.1*sin(180)+(180)/360+.05}) {\footnotesize#4};
	}
}

\usetikzlibrary{intersections}

\tikzset{
	dot/.style={circle,inner sep=1pt,fill},
}

\newcommand{\pants}[8]{
	\draw (0,2) ellipse (1.5 and .1);
	\draw[name path=left edge] (-2.5,-2)to[out=90,in=-90](-1.5,2);
	\draw[dashed] ($(-1.5, -2) + (180:1 and .1)$() arc (180:0:1 and .1);
	\draw ($(-1.5, -2) + (180:1 and .1)$() arc (180:360:1 and .1);
	\draw[name path=right edge] (2.5,-2)to[out=90,in=-90](1.5,2);
	\draw[dashed] ($(1.5, -2) + (180:1 and .1)$() arc (180:0:1 and .1);
	\draw ($(1.5, -2) + (180:1 and .1)$() arc (180:360:1 and .1);
	\draw[name path=centre] (-.5,-2)to[out=90,in=180](0,0)to[out=0,in=90](.5,-2);
	\draw[thick,red,dashed,name path=C] ($(-1.5, -2) + (70:1 and .1)$()to[out=90,in=-90]($(0, 2) + (100:1.5 and .1)$();
	\draw[thick,red,dashed,name path=D] ($(1.5, -2) + (70:1 and .1)$()to[out=90,in=-90]($(0, 2) + (40:1.5 and .1)$();
	\path[thick,red,name path=A] ($(-1.5, -2) + (250:1 and .1)$()to[out=90,in=-90]($(0, 2) + (230:1.5 and .1)$();
	\path[thick,red,name path=B] ($(1.5, -2) + (250:1 and .1)$()to[out=90,in=-90]($(0, 2) + (280:1.5 and .1)$();
	\ifthenelse{\equal{#1}{}}{}
	{
		\begin{scope}
			\path[name path=X] (-.5,-.25)--(-2.5,0);
			\path [name intersections={of=C and X,by={c1}}];
			\path [name intersections={of=X and left edge,by={c2}}];
			\draw[blue,dashed] (c1)to[out=140,in=70] (c2);
			\path[name path=Y] (c2) -- (0,-2);
			\path [name intersections={of=Y and A,by={d1}}];
			\draw[blue] (c2) to[out=-110,in=220]  (d1);
			\node[left,inline text,blue] at (c2) {\footnotesize#1};
		\end{scope}
	}
	\ifthenelse{\equal{#2}{}}{}
	{
		\begin{scope}
			\path[name path=X] (-2.5,-2)--(0,-.5);
			\path [name intersections={of=A and X,by={c1}}];
			\path [name intersections={of=X and centre,by={c2}}];
			\draw[blue] (c1) to[out=45,in=260]  (c2);
			\path[name path=Y] (c2) -- (-2.5,2);
			\path [name intersections={of=Y and C,by={d1}}];
			\draw[blue,dashed] (c2) to[out=95,in=-40]  (d1);
			\node[right,inline text,blue] at (c2) {\footnotesize#2};
		\end{scope}
	}
	\ifthenelse{\equal{#3}{}}{}
	{
		\begin{scope}
			\path[name path=X] (0,-.5)--(2.5,-.5);
			\path[name intersections={of=D and X,by={c1}}];
			\path[name intersections={of=X and centre,by={c2}}];
			\draw[blue,dashed] (c1)to[out=140,in=100] (c2);
			\path[name path=Y] (c2) -- (2,-2);
			\path [name intersections={of=Y and B,by={d1}}];
			\draw[blue] (c2) to[out=-80,in=220]  (d1);
			\node[left,inline text,blue] at (c2) {\footnotesize#3};
		\end{scope}
	}
	\ifthenelse{\equal{#4}{}}{}
	{
		\begin{scope}
			\path[name path=X] (0,-2.5)--(2.5,-1);
			\path [name intersections={of=B and X,by={c1}}];
			\path [name intersections={of=X and right edge,by={c2}}];
			\draw[blue] (c1) to[out=45,in=280]  (c2);
			\path[name path=Y] (c2) -- (0,.5);
			\path [name intersections={of=Y and D,by={d1}}];
			\draw[blue,dashed] (c2) to[out=110,in=-40]  (d1);
			\node[right,inline text,blue] at (c2) {\footnotesize#4};
		\end{scope}
	}
	\ifthenelse{\equal{#5}{}}{}
	{
		\begin{scope}
			\path[name path=X] (-2.5,2)--(0,.5);
			\path [name intersections={of=left edge and X,by={c1}}];
			\path [name intersections={of=X and A,by={c2}}];
			\draw[blue] (c1) to[out=-110,in=250]  (c2);
			\path[name path=Y] (c1) -- (0,1.5);
			\path [name intersections={of=Y and C,by={d1}}];
			\draw[blue,dashed] (c1) to[out=80,in=110]  (d1);
			\node[left,inline text,blue] at (c1) {\footnotesize#5};
		\end{scope}
	}
	\ifthenelse{\equal{#6}{}}{}
	{
		\begin{scope}
			\path[name path=X] (-2.5,.5)--(2.5,.5);
			\path [name intersections={of=A and X,by={c1}}];
			\path [name intersections={of=X and B,by={c2}}];
			\draw[blue] (c1) to[out=70,in=250]  (c2);
			\node[left,inline text,blue] at ($ (c2) + (0,.1) $) {\footnotesize#6};
		\end{scope}
	}
	\ifthenelse{\equal{#7}{}}{}
	{
		\begin{scope}
			\path[name path=X] (0,.5)--(2.5,1);
			\path [name intersections={of=B and X,by={c1}}];
			\path [name intersections={of=X and right edge,by={c2}}];
			\draw[blue] (c1) to[out=70,in=280]  (c2);
			\path[name path=Y] (c2) -- (0,2.2);
			\path [name intersections={of=Y and D,by={d1}}];
			\draw[blue,dashed] (c2) to[out=110,in=280]  (d1);
			\node[right,inline text,blue] at (c2) {\footnotesize#7};
		\end{scope}
	}
	\ifthenelse{\equal{#8}{}}{}
	{
		\begin{scope}
			\path[name path=X] (-2.5,1.3)--(2.5,1.3);
			\path [name intersections={of=D and X,by={c1}}];
			\path [name intersections={of=X and C,by={c2}}];
			\draw[blue,dashed] (c1) to[out=120,in=300]  (c2);
			\node[left,inline text,blue] at (c1) {\footnotesize#8};
		\end{scope}
	}
	\draw[thick,red] ($(-1.5, -2) + (250:1 and .1)$()to[out=90,in=-90]($(0, 2) + (230:1.5 and .1)$();
	\draw[thick,red] ($(1.5, -2) + (250:1 and .1)$()to[out=90,in=-90]($(0, 2) + (280:1.5 and .1)$();
	\node[below,inline text,red] at ($(-1.5, -2) + (250:1 and .1)$() {\footnotesize\ta};
	\node[below,inline text,red] at ($(1.5, -2) + (250:1 and .1)$() {\footnotesize\tb};
	\node[below,inline text,red] at ($(-1.5, -2) + (-70:1 and .1)$() {\footnotesize\tc};
	\node[below,inline text,red] at ($(1.5, -2) + (-70:1 and .1)$() {\footnotesize\td};
	\node[above,inline text,red] at ($(0, 2) + (-230:1.5 and .1)$() {\footnotesize\tap};
	\node[above,inline text,red] at ($(0, 2) + (-280:1.5 and .1)$() {\footnotesize\tbp};
	\node[above,inline text,red] at ($(0, 2) + (100:1.5 and .1)$() {\footnotesize\tcp};
	\node[above,inline text,red] at ($(0, 2) + (40:1.5 and .1)$() {\footnotesize\tdp};
}

\def\ladder[#1][#2][#3][#4][#5][#6]{
	\path(#1);
	\pgfgetlastxy{\XCoord}{\YCoord};
	\begin{scope}[xshift=\XCoord,yshift=\YCoord]
		\draw (-.6,-.5)--(-.6,.5);
		\draw (.6,-.5)--(.6,.5);
		\ifthenelse{\equal{#4}{}}{}
			{
				\draw (-.6,-.25)--(.6,.25);
			};
		\node[below,inline text] at (-.6,-.5) {\footnotesize#2};
		\node[below,inline text] at (.6,-.5) {\footnotesize#3};
		\node[below,inline text] at (0,0) {\footnotesize#4};
		\node[above,inline text] at (-.6,.5) {\footnotesize#5};
		\node[above,inline text] at (.6,.5) {\footnotesize#6};
	\end{scope}
}

\def\Lladder[#1][#2][#3][#4][#5][#6][#7]{
	\path(#1);
	\pgfgetlastxy{\XCoord}{\YCoord};
	\begin{scope}[xshift=\XCoord,yshift=\YCoord]
		\draw (-1.8,-.5)--(-1.8,.5);
		\draw (-.6,-.5)--(-.6,.5);
		\draw (.6,-.5)--(.6,.5);
		\ifthenelse{\equal{#4}{}}{}
		{
			\draw (-.6,-.25)--(.6,.25);
		};
		\node[below,inline text] at (-.6,-.5) {\footnotesize#2};
		\node[below,inline text] at (.6,-.5) {\footnotesize#3};
		\node[below,inline text] at (0,0) {\footnotesize#4};
		\node[above,inline text] at (-.6,.5) {\footnotesize#5};
		\node[above,inline text] at (.6,.5) {\footnotesize#6};
		\node[below,inline text] at (-1.8,-.5) {\footnotesize#7};
	\end{scope}
}

\def\Lmorphism[#1][#2][#3][#4][#5][#6][#7]{
	\path(#1);
	\pgfgetlastxy{\XCoord}{\YCoord};
	\begin{scope}[xshift=\XCoord,yshift=\YCoord]
		\draw (-1.8,-.5)to[out=90,in=190](-.6,.25);
		\draw (-.6,-.5)--(-.6,.5);
		\draw (.6,-.5)--(.6,.5);
		\ifthenelse{\equal{#4}{}}{}
		{
			\draw (-.6,-.25)--(.6,.25);
		};
		\node[below,inline text] at (-.6,-.5) {\footnotesize#2};
		\node[below,inline text] at (.6,-.5) {\footnotesize#3};
		\node[below,inline text] at (0,0) {\footnotesize#4};
		\node[above,inline text] at (-.6,.5) {\footnotesize#5};
		\node[above,inline text] at (.6,.5) {\footnotesize#6};
		\node[below,inline text] at (-1.8,-.5) {\footnotesize#7};
	\end{scope}
}

\def\Rladder[#1][#2][#3][#4][#5][#6][#7]{
	\path(#1);
	\pgfgetlastxy{\XCoord}{\YCoord};
	\begin{scope}[xshift=\XCoord,yshift=\YCoord]
		\draw (1.8,-.5)--(1.8,.5);
		\draw (-.6,-.5)--(-.6,.5);
		\draw (.6,-.5)--(.6,.5);
		\ifthenelse{\equal{#4}{}}{}
		{
			\draw (-.6,-.25)--(.6,.25);
		};
		\node[below,inline text] at (-.6,-.5) {\footnotesize#2};
		\node[below,inline text] at (.6,-.5) {\footnotesize#3};
		\node[below,inline text] at (0,0) {\footnotesize#4};
		\node[above,inline text] at (-.6,.5) {\footnotesize#5};
		\node[above,inline text] at (.6,.5) {\footnotesize#6};
		\node[below,inline text] at (1.8,-.5) {\footnotesize#7};
	\end{scope}
}

\def\Rmorphism[#1][#2][#3][#4][#5][#6][#7]{
	\path(#1);
	\pgfgetlastxy{\XCoord}{\YCoord};
	\begin{scope}[xshift=\XCoord,yshift=\YCoord]
		\draw (1.8,-.5)to[out=90,in=-10](.6,.3);
		\draw (-.6,-.5)--(-.6,.5);
		\draw (.6,-.5)--(.6,.5);
		\ifthenelse{\equal{#4}{}}{}
		{
			\draw (-.6,-.35)--(.6,.15);
		};
		\node[below,inline text] at (-.6,-.5) {\footnotesize#2};
		\node[below,inline text] at (.6,-.5) {\footnotesize#3};
		\node[below,inline text] at (0,-.1) {\footnotesize#4};
		\node[above,inline text] at (-.6,.5) {\footnotesize#5};
		\node[above,inline text] at (.6,.5) {\footnotesize#6};
		\node[below,inline text] at (1.8,-.5) {\footnotesize#7};
	\end{scope}
}

\def\Laction[#1][#2][#3]{
		\draw (0,-.5)--(0,.5);
		\draw (-.6,-.5)to[out=90,in=210](0,0);
		\node[below,inline text] at (-.6,-.5) {\footnotesize#1};
		\node[below,inline text] at (0,-.5) {\footnotesize#2};
		\node[above,inline text]  at (0,.5) {\footnotesize#3};
}
\def\Raction[#1][#2][#3]{
	\draw (0,-.5)--(0,.5);
	\draw (.6,-.5)to[out=90,in=-30](0,0);
	\node[below,inline text] at (.6,-.5) {\footnotesize#1};
	\node[below,inline text] at (0,-.5) {\footnotesize#2};
	\node[above,inline text] at (0,.5) {\footnotesize#3};
}

\def\Lassociator[#1][#2][#3][#4]{
	\draw (0,-.5)--(0,.5);
	\draw (-.6,-.5)to[out=90,in=210](0,0);
	\draw (.6,-.5)to[out=90,in=-30](0,.25);
	\node[below,inline text] at (-.6,-.5) {\footnotesize#1};
	\node[below,inline text] at (.6,-.5) {\footnotesize#2};
	\node[below,inline text] at (0,-.5) {\footnotesize#3};
	\node[above,inline text] at (0,.5) {\footnotesize#4};
}
\def\Rassociator[#1][#2][#3][#4]{
	\draw (0,-.5)--(0,.5);
	\draw (-.6,-.5)to[out=90,in=210](0,.25);
	\draw (.6,-.5)to[out=90,in=-30](0,0);
	\node[below,inline text] at (-.6,-.5) {\footnotesize#1};
	\node[below,inline text] at (.6,-.5) {\footnotesize#2};
	\node[below,inline text] at (0,-.5) {\footnotesize#3};
	\node[above,inline text] at (0,.5) {\footnotesize#4};
}

\usepackage{hhline}
\usepackage{etoolbox}

\let\originalleft\left
\let\originalright\right
\renewcommand{\left}{\mathopen{}\mathclose\bgroup\originalleft}
\renewcommand{\right}{\aftergroup\egroup\originalright}

\newcommand{\vvec}[1]
{
	\operatorname{\bf Vec}
	\ifstrequal{#1}{}
	{}
	{\left(#1\right)}
}

\newcommand{\vvectwist}[2]
{
	\operatorname{\bf Vec}^{#2}
	\ifstrequal{#1}{}
	{}
	{\left(#1\right)}
}

\newcommand{\lad}[1]
{
	\operatorname{\bf Lad}
	\ifstrequal{#1}{}
	{}
	{\left(#1\right)}
}

\newcommand{\ann}[3]
{
	\operatorname{\bf Ann}_{\mathcal{#1},\mathcal{#2}}
	\ifstrequal{#3}{}
	{}
	{\left(\mathcal{#3}\right)}
}

\newcommand{\kar}[1]
{
	\operatorname{\bf Kar}
	\ifstrequal{#1}{}
	{}
	{\left(#1\right)}
}

\newcommand{\bpr}[1]
{
	\operatorname{\bf BPR}
	\ifstrequal{#1}{}
	{}
	{\left(#1\right)}
}

\newcommand{\dih}[1]
{
	\operatorname{Dih}
	\ifstrequal{#1}{}
	{}
	{_{#1}}
}

\newcommand{\defect}[5]{
	\ifthenelse{\equal{#3}{}}
	{
		\begin{smallmatrix} #2\hfill	\\ #1\hfill \end{smallmatrix}\!\bigr\vert^{#5}
	}
	{
		\ifthenelse{\equal{#4}{}}
		{
			\begin{smallmatrix} #2\hfill	\\ #1\hfill \end{smallmatrix}\!\bigr\vert_{#3}^{#5}
		}
		{
			\begin{smallmatrix} #2\hfill	\\ #1\hfill \end{smallmatrix}\!\bigr\vert_{(#3,#4)}^{#5}
		}
	}
}

\newcommand{\twodownoneuptube}[6]{
	\node[above,inline text,red] at (0,1) {#4};
	\node[below,inline text,red] at (-.25,-1) {#5};
  \node[below,inline text,red] at (.25,-1) {#6};
	\draw[blue] (-.25,-.5) to [out=90+45,in=-90] (-.5,0) to[out=90,in=180+45] (0,.5) node[left] at (-.5,0) {#1};
	\draw[blue] (.25,-.5) to [out=45,in=-90] (.5,0) to[out=90,in=-45] (0,.75) node[right] at (.5,0) {#2};
	\draw[blue] (-.25,-.65) to [out=45,in=180+45] (.25,-.65) node[below] at (0,-.65) {#3};
\begin{scope}[even odd rule]
	\clip[draw] (0,0) circle (1) (0,0) circle (.33);
	\draw [red,very thick] ({0},{0})--({0},{1});
	\draw [red,very thick] ({-.25},{-1})--({-.25},{0});
	\draw [red,very thick] ({.25},{-1})--({.25},{0});
	\draw (0,0) circle (1) (0,0) circle (.33);
\end{scope}
  }

\newcommand{\onedowntwouptube}[6]{
	\node[below,inline text,red] at (0,-1) {#4};
	\node[above,inline text,red] at (-.25,1) {#5};
  \node[above,inline text,red] at (.25,1) {#6};
	\draw[blue] (-.25,.5) to [out=180+45,in=90] (-.5,0) to [out=-90,in=180-45] (0,-.5) node[left] at (-.5,0) {#1};
	\draw[blue] (.25,.5) to [out=-45,in=90] (.5,0) to [out=-90,in=45] (0,-.75) node[right] at (.5,0) {#2};
	\draw[blue] (-.25,.65) to [out=45,in=180+45] (.25,.75) node[below] at (0,.65) {#3};
	\begin{scope}[even odd rule]
		\clip[draw] (0,0) circle (1) (0,0) circle (.33);
		\draw [red,very thick] ({0},{0})--({0},{-1});
		\draw [red,very thick] ({-.25},{1})--({-.25},{0});
		\draw [red,very thick] ({.25},{1})--({.25},{0});
		\draw (0,0) circle (1) (0,0) circle (.33);
	\end{scope}
  }

\newcommand{\diamondvectorexampleone}[9]{
        \draw[thick,red] (0,0) -- (0,2);
        \draw[thick,orange] (0,2) [out=90+45,in=-90] to (-1,4);
        \draw[thick,blue] (0,2) to [out=45,in=-90] (1,4);
        \draw[thick,nicegreen] (-1,4) to [out=90,in=180+45] (0,6);
        \draw[thick,blue] (1,4) to [out=90,in=-45] (0,6);
        \draw[thick,red] (0,6) -- (0,8);
        \filldraw[fill=black,draw=black] (0,2) circle (3pt);
        \filldraw[fill=black,draw=black] (-1,4) circle (3pt);
        \filldraw[fill=black,draw=black] (1,4) circle (3pt);
        \filldraw[fill=black,draw=black] (0,6) circle (3pt);
        \node at (0,0) [below] {#1};
        \node at (0,8) [above] {#2};
        \node at (-1,4) [left] {#3};
        \node at (1,4) [right] {#4};
        \node at (-0.9,3) [left] {#5};
        \node at (0.9,3) [right] {#6};
        \node at (0.9,5) [right] {#7};
        \node at (-0.9,5) [left] {#8};
        \node at (0,6) [left] {#9};
}

\newcommand{\diamondvectorexampletwo}[8]{
        \draw[thick,red] (0,0) -- (0,2);
        \draw[thick,orange] (0,2) [out=90+45,in=-90] to (-1,4);
        \draw[thick,red] (0,2) to [out=45,in=-90] (1,4);
        \draw[thick,nicegreen] (-1,4) to [out=90,in=180+45] (0,6);
        \draw[thick,blue] (1,4) to [out=90,in=-45] (0,6);
        \draw[thick,blue] (0,6) -- (0,8);
        \filldraw[fill=black,draw=black] (0,2) circle (3pt);
        \filldraw[fill=black,draw=black] (-1,4) circle (3pt);
        \filldraw[fill=black,draw=black] (1,4) circle (3pt);
        \filldraw[fill=black,draw=black] (0,6) circle (3pt);
        \node at (0,0) [below] {#1};
        \node at (0,8) [above] {#2};
        \node at (-1,4) [left] {#3};
        \node at (1,4) [right] {#4};
        \node at (-0.9,3) [left] {#5};
        \node at (0.9,3) [right] {#6};
        \node at (0.9,5) [right] {#7};
        \node at (-0.9,5) [left] {#8};
}

\newcommand{\diamondvectorexamplethree}[8]{
        \draw[thick,red] (0,0) -- (0,2);
        \draw[thick,orange] (0,2) [out=90+45,in=-90] to (-1,4);
        \draw[thick,yellow] (0,2) to [out=45,in=-90] (1,4);
        \draw[thick,nicegreen] (-1,4) to [out=90,in=180+45] (0,6);
        \draw[thick,purple] (1,4) to [out=90,in=-45] (0,6);
        \draw[thick,blue] (0,6) -- (0,8);
        \filldraw[fill=black,draw=black] (0,2) circle (3pt);
        \filldraw[fill=black,draw=black] (-1,4) circle (3pt);
        \filldraw[fill=black,draw=black] (1,4) circle (3pt);
        \filldraw[fill=black,draw=black] (0,6) circle (3pt);
        \node at (0,0) [below] {#1};
        \node at (0,8) [above] {#2};
        \node at (-1,4) [left] {#3};
        \node at (1,4) [right] {#4};
        \node at (-0.9,3) [left] {#5};
        \node at (0.9,3) [right] {#6};
        \node at (0.9,5) [right] {#7};
        \node at (-0.9,5) [left] {#8};
}

\newcommand{\bimoduleassociatorexample}[7]{
          \draw[thick,red] (0,0) -- (0,2);
          \draw[thick,blue] (0,2) -- (0,6);
          \draw[thick,nicegreen] (0,2) -- (2,4);
          \draw[thick,orange] (2,4) -- (2,8);
          \draw[thick,red] (0,6) -- (2,8);
          \draw[thick,red] (2,8) -- (2,10);
          \draw[thick,nicegreen] (2,4) -- (0,6);
          \node[below] at (0,0) {#1};
          \node[left] at (0,4) {#2};
          \node[below] at (1.2,3) {#3};
          \node[below] at (1,5.2) {#4};
          \node[right] at (2,6) {#5};
          \node[left] at (1.5,7.5) {#6};
          \node[above] at (2,10) {#7};
          \filldraw[fill=black,draw=black] (0,2) circle (3pt);
          \filldraw[fill=black,draw=black] (0,6) circle (3pt);
          \filldraw[fill=black,draw=black] (2,4) circle (3pt);
          \filldraw[fill=black,draw=black] (2,8) circle (3pt);
}

\newcommand{\ssbivertex}[4]{
	\draw[thick,red] (0,-.5)--(0,.5) node[above,inline text,pos=1] {$#2$} node[below,inline text,pos=0] {$#1$} node[left, inline text,pos=0.5] {$#3$}  node[right, inline text,pos=0.5] {$#4$};
	\filldraw[black] (0,0) circle (.05);
}

\newcommand{\genericbimoduleassociator}[3]{
\draw[gray,thick,dashed,rounded corners=15pt] (-2,0) rectangle (6,11);
  \draw [thick,red] (0,0)--(0,2);
  \draw [thick,blue] (0,2)--(0,7) node[midway,left] {#1};
  \draw [thick,teal] (0,2)--(4,4);
  \draw [thick,nicegreen] (4,4)--(0,7) node[midway,above,xshift=3pt] {#2};
  \draw [thick,orange] (0,7)--(4,9);
  \draw [thick,purple] (4,4)--(4,9) node[midway,right] {#3};
  \draw [thick,red] (4,9)--(4,11);
  \filldraw[fill=white, draw=black] (0,2) circle (10pt);
  \filldraw[fill=white, draw=black] (4,4) circle (10pt);
  \filldraw[fill=white, draw=black] (0,7) circle (10pt);
  \filldraw[fill=white, draw=black] (4,9) circle (10pt);
}

\newcommand{\vertcompcompoundrep}[5]{
          \draw[thick,red] (0,-0.5)--(0,0) node[below, inline text,pos=0] {$#1$} node[left,inline text] {$#2$};
          \draw[thick,darkgreen] (0,0)--(0,0.5) node[right, inline text, pos=0.5] {$#3$};
          \draw[thick,red] (0,0.5)--(0,1) node[above,inline text] {$#5$} node[left,inline text,pos=0] {$#4$};
          \filldraw[black] (0,0) circle (.05);
          \filldraw[black] (0,.5) circle (.05);
}

\newcommand{\stbivertex}[4]{
	\draw[thick,red] (0,-.5)--(0,0) node[below,inline text,pos=0] {$#1$} node[left, inline text,pos=1] {$#3$}  node[right, inline text,pos=1] {$#4$};
	\draw[thick,darkgreen] (0,0)--(0,.5) node[above,inline text,pos=1] {$#2$};
	\filldraw[black] (0,0) circle (.05);
}

\newcommand{\ssstrivertex}[5]{
	\draw[thick,red] (0,0)--(0,.5) node[above,inline text,pos=1] {$#3$};
	\draw[thick,red] (-.2,-.5)--(0,0) node[below,inline text,pos=0,anchor=north east] {$#1$};
	\draw[thick,red] (.2,-.5)--(0,0) node[below,inline text,pos=0,anchor=north west] {$#2$};
	\filldraw[black] (0,0) circle (.05) node[left,inline text] {$#4$};
}

\newcommand{\ststrivertex}[5]{
	\draw[thick,red] (0,0)--(0,.5) node[above,inline text,pos=1] {$#3$};
	\draw[thick,red] (-.2,-.5)--(0,0) node[below,inline text,pos=0,anchor=north east] {$#1$};
	\draw[thick,darkgreen] (.2,-.5)--(0,0) node[below,inline text,pos=0,anchor=north west] {$#2$};
	\filldraw[black] (0,0) circle (.05) node[left,inline text] {$#4$};
}

\newcommand{\tsstrivertex}[5]{
	\draw[thick,red] (0,0)--(0,.5) node[above,inline text,pos=1] {$#3$};
	\draw[thick,darkgreen] (-.2,-.5)--(0,0) node[below,inline text,pos=0,anchor=north east] {$#1$};
	\draw[thick,red] (.2,-.5)--(0,0) node[below,inline text,pos=0,anchor=north west] {$#2$};
	\filldraw[black] (0,0) circle (.05) node[left,inline text] {$#4$};
}

\newcommand{\rsttrivertex}[5]{
	\draw[thick,red] (0,0)--(0,.5) node[above,inline text,pos=1] {$#3$};
	\draw[thick,darkgreen] (-.2,-.5)--(0,0) node[below,inline text,pos=0,anchor=north east] {$#1$};
	\draw[thick,orange] (.2,-.5)--(0,0) node[below,inline text,pos=0,anchor=north west] {$#2$};
	\filldraw[black] (0,0) circle (.05) node[left,inline text] {$#4$};
}

\newcommand{\ssstrivertexdual}[5]{
	\begin{scope}[yscale=-1]
		\draw[thick,red] (0,0)--(0,.5) node[below,inline text,pos=1] {$#3$};
		\draw[thick,red] (-.2,-.5)--(0,0) node[above,inline text,pos=0,anchor=south east] {$#1$};
		\draw[thick,red] (.2,-.5)--(0,0) node[above,inline text,pos=0,anchor=south west] {$#2$};
		\filldraw[black] (0,0) circle (.05) node[left,inline text] {$#4$};
	\end{scope}
}

\newcommand{\ststrivertexdual}[5]{
	\begin{scope}[yscale=-1]
		\draw[thick,red] (0,0)--(0,.5) node[below,inline text,pos=1] {$#3$};
		\draw[thick,red] (-.2,-.5)--(0,0) node[above,inline text,pos=0,anchor=south east] {$#1$};
		\draw[thick,darkgreen] (.2,-.5)--(0,0) node[above,inline text,pos=0,anchor=south west] {$#2$};
		\filldraw[black] (0,0) circle (.05) node[left,inline text] {$#4$};
	\end{scope}
}

\newcommand{\tsstrivertexdual}[5]{
	\begin{scope}[yscale=-1]
		\draw[thick,red] (0,0)--(0,.5) node[below,inline text,pos=1] {$#3$};
		\draw[thick,darkgreen] (-.2,-.5)--(0,0) node[above,inline text,pos=0,anchor=south east] {$#1$};
		\draw[thick,red] (.2,-.5)--(0,0) node[above,inline text,pos=0,anchor=south west] {$#2$};
		\filldraw[black] (0,0) circle (.05) node[left,inline text] {$#4$};
	\end{scope}
}

\newcommand{\rsttrivertexdual}[5]{
	\begin{scope}[yscale=-1]
	\draw[thick,red] (0,0)--(0,.5) node[below,inline text,pos=1] {$#3$};
	\draw[thick,darkgreen] (-.2,-.5)--(0,0) node[above,inline text,pos=0,anchor=south east] {$#1$};
	\draw[thick,orange] (.2,-.5)--(0,0) node[above,inline text,pos=0,anchor=south west] {$#2$};
	\filldraw[black] (0,0) circle (.05) node[left,inline text] {$#4$};
	\end{scope}
}

\newcommand{\trivalentvertex}[4]{
   \draw (0,0)--(0,1);
   \node[above,inline text] at (0,1) {#1};
   \draw (0,0)--(-0.707,-0.707);
   \node[below,inline text] at (-0.707,-0.707) {#2};
   \draw (0,0)--(0.707,-0.707);
   \node[below,inline text] at (0.707,-0.707) {#3};
   \node[left,inline text] at (0,0) {#4};
}

\newcommand{\vgeneralpantsstp}[6]
{
	\ifthenelse{\equal{#1}{}}{}
	{
		\coordinate (A1) at (0,-1);
		\coordinate (B1) at (-.3,-.6);
		\coordinate (C1) at (0,-.2);
		\node[left,inline text,blue] at (C1) {\footnotesize#1};
		\draw[blue] (A1) to[out=90+45,in=270] (B1) to[out=90,in=180+45] (C1);
	}
	\ifthenelse{\equal{#2}{}}{}
	{
		\coordinate (A1) at (0,-1.1);
		\coordinate (B1) at (.3,-.6);
		\coordinate (C1) at (0,-.1);
		\node[right,inline text,blue] at (C1) {\footnotesize#2};
		\draw[blue] (A1) to[out=90-45,in=270] (B1) to[out=90,in=0-45] (C1);
	}
	\ifthenelse{\equal{#3}{}}{}
	{
		\coordinate (A1) at (0,.2);
		\coordinate (B1) at (-.3,.6);
		\coordinate (C1) at (0,1);
		\node[left,inline text,blue] at (A1) {\footnotesize#3};
		\draw[blue] (A1) to[out=90+45,in=270] (B1) to[out=90,in=180+45] (C1);
	}
	\ifthenelse{\equal{#4}{}}{}
	{
		\coordinate (A1) at (0,.1);
		\coordinate (B1) at (.3,.6);
		\coordinate (C1) at (0,1.1);
		\node[right,inline text,blue] at (A1) {\footnotesize#4};
		\draw[blue] (A1) to[out=90-45,in=270] (B1) to[out=90,in=0-45] (C1);
	}
	\ifthenelse{\equal{#5}{}}{}
	{
		\coordinate (A1) at (0,-1.3);
		\coordinate (B1) at (-.9,0);
		\coordinate (C1) at (0,1.3);
		\node[left,inline text,blue] at (B1) {\footnotesize#5};
		\draw[blue] (A1) to[out=90+45,in=270] (B1) to[out=90,in=180+45] (C1);
	}
	\ifthenelse{\equal{#6}{}}{}
	{
		\coordinate (A1) at (0,-1.4);
		\coordinate (B1) at (.9,0);
		\coordinate (C1) at (0,1.4);
		\node[right,inline text,blue] at (B1) {\footnotesize#6};
		\draw[blue] (A1) to[out=90-45,in=270] (B1) to[out=90,in=0-45] (C1);
	}
	\begin{scope}[even odd rule]
		\clip (0,0) ellipse [x radius=1,y radius=1.5] (0,-.6) circle (.2) (0,.6) circle (.2);
		\draw [red,ultra thick] (0,-1.5)--(0,-.6);
		\draw [orange,ultra thick] (0,-.6)--(0,.6);
		\draw [darkgreen,ultra thick] (0,.6)--(0,1.5);
	\end{scope}
	\draw (0,-.6) circle (.2);
	\draw (0,.6) circle (.2);
	\draw (0,0) ellipse [x radius=1.5,y radius=1.5];
	\node[red,inline text,below] at (0,-1.5) {\footnotesize\ta};
	\node[orange,inline text,above] at (0,.4) {\footnotesize\tb};
	\node[darkgreen,inline text,above] at (0,1.5) {\footnotesize\tc};
}

\newcommand{\vpantsstp}[4]
{
	\ifthenelse{\equal{#1}{}}{}
	{
		\coordinate (A1) at (0,-1);
		\coordinate (B1) at (-.3,-.6);
		\coordinate (C1) at (0,-.2);
		\node[left,inline text,blue] at (B1) {\footnotesize#1};
		\draw[blue] (A1) to[out=90+45,in=270] (B1) to[out=90,in=180+45] (C1);
	}
	\ifthenelse{\equal{#2}{}}{}
	{
		\coordinate (A1) at (0,-1.1);
		\coordinate (B1) at (.3,-.6);
		\coordinate (C1) at (0,-.1);
		\node[right,inline text,blue] at (B1) {\footnotesize#2};
		\draw[blue] (A1) to[out=90-45,in=270] (B1) to[out=90,in=0-45] (C1);
	}
	\ifthenelse{\equal{#3}{}}{}
	{
		\coordinate (A1) at (0,.2);
		\coordinate (B1) at (-.3,.6);
		\coordinate (C1) at (0,1);
		\node[left,inline text,blue] at (B1) {\footnotesize#3};
		\draw[blue] (A1) to[out=90+45,in=270] (B1) to[out=90,in=180+45] (C1);
	}
	\ifthenelse{\equal{#4}{}}{}
	{
		\coordinate (A1) at (0,.1);
		\coordinate (B1) at (.3,.6);
		\coordinate (C1) at (0,1.1);
		\node[right,inline text,blue] at (B1) {\footnotesize#4};
		\draw[blue] (A1) to[out=90-45,in=270] (B1) to[out=90,in=0-45] (C1);
	}
	\begin{scope}[even odd rule]
		\clip (0,0) ellipse [x radius=1,y radius=1.5] (0,-.6) circle (.2) (0,.6) circle (.2);
		\draw [red,ultra thick] (0,-1.5)--(0,-.6);
		\draw [orange,ultra thick] (0,-.6)--(0,.6);
		\draw [darkgreen,ultra thick] (0,.6)--(0,1.5);
	\end{scope}
	\draw (0,-.6) circle (.2);
	\draw (0,.6) circle (.2);
	\draw (0,0) ellipse [x radius=1.5,y radius=1.5];
	\node[red,inline text,below] at (0,-1.5) {\footnotesize\ta};
	\node[orange,inline text,above] at (0,.4) {\footnotesize\tb};
	\node[darkgreen,inline text,above] at (0,1.5) {\footnotesize\tc};
}

\newcommand{\faceOpenNh}[6]
{
        \draw (1., 0.) -- (0.5, 0.866025) node[midway,fill=white,draw=black!30,circle,inner sep=.1] {#1};
        \draw (0.5, 0.866025) -- (-0.5, 0.866025) node[midway,fill=white,draw=black!30,circle,inner sep=.1] {#2};
        \draw (-0.5, 0.866025) -- (-1., 0.) node[midway,fill=white,draw=black!30,circle,inner sep=.1] {#3};
        \draw (-1., 0.) -- (-0.5, -0.866025) node[midway,fill=white,draw=black!30,circle,inner sep=.1] {#4};
        \draw (-0.5, -0.866025) -- (0.5, -0.866025) node[midway,fill=white,draw=black!30,circle,inner sep=.1] {#5};
        \draw(0.5, -0.866025) -- (1,0) node[midway,fill=white,draw=black!30,circle,inner sep=.1] {#6};
        \draw (1., 0.) -- (1.5*1.,1.5* 0.);
        \draw (0.5, 0.866025) -- (1.5*0.5, 1.5*0.866025);
        \draw (-0.5, 0.866025) -- (-1.5*0.5, 1.5*0.866025);
        \draw (-1., 0.) -- (-1.5*1., 0.);
        \draw (-0.5, -0.866025) -- (-1.5*0.5, -1.5*0.866025);
        \draw (0.5, -0.866025) -- (1.5*0.5, -1.5*0.866025);
        \fill (1., 0.) circle (.05);\fill (0.5, 0.866025) circle (.05);\fill (-0.5, 0.866025) circle (.05);
        \fill (-1., 0.) circle (.05);\fill (0.5, -0.866025) circle (.05);\fill (-0.5, -0.866025) circle (.05);
        \draw[white,very thick] (0,0) circle (1.5);
        \draw[gray,dashed,very thick] (0,0) circle (1.5);
}

\newlength{\tabwidth}
\newlength{\tabheight}

\usepackage{soul}

\allowdisplaybreaks
\begin{document}

\title{Computing Defects Associated to Bounded Domain Wall Structures: \\The $\mathbb{Z}/p\mathbb{Z}$ Case}

\author{Jacob C.\ Bridgeman}
\email{jcbridgeman1@gmail.com}
\affiliation{Perimeter Institute for Theoretical Physics, Waterloo, Ontario, Canada}

\author{Daniel Barter}
\email{danielbarter@gmail.com}
\affiliation{Mathematical Sciences Institute, Australian National University, Canberra, Australia}

\date{\today}

\defcitealias{MR2677836}{Etingof \emph{et al.}, Quantum Topol. \textbf{1}, 209 (2010)}

\begin{abstract}
	A domain wall structure consists of a planar graph with faces labeled by fusion categories/topological phases. Edges are labeled by bimodules/domain walls. When the vertices are labeled by point defects we get a compound defect. We present an algorithm, called the domain wall structure algorithm, for computing the compound defect. We apply this algorithm to show that the \emph{bimodule associator}, related to the $O_3$ obstruction of \citepalias{MR2677836}, is trivial for all domain walls of $\vvec{\ZZ{p}}$.
	
	In the language of this paper, the ground states of the Levin-Wen model are compound defects. We use this to define a generalization of the Levin-Wen model with domain walls and point defects. The domain wall structure algorithm can be used to compute the ground states of these generalized Levin-Wen type models.
\end{abstract}

\maketitle

% This line sets the project root file.
% !TEX root = ../defects_domain_wall_structures.tex
% !TEX spellcheck = en_US

%\section{Introduction}

Due to their insensitivity to environmental noise, topological phases have promise as materials for encoding quantum information\cite{MR1951039,Brown2014,Terhal2015,Dennis2002}. By braiding and fusing the emergent quasi-particle excitations, the encoded information can be manipulated in a robust manner. Such protection from the environment is an important requirement for any large scale quantum device. In many phases, especially those most suited to laboratory realization, the quantum computational power is severely limited. It has become clear that the inclusion of defects can improve the materials from this perspective\cite{Dennis2002,0610082,Bombin2007a,Bombin2010,Brown2013a,Pastawski2014,Yoshida2015a,1606.07116,Brown2016,IrisCong1,IrisCong2,PhysRevB.96.195129,Yoshida2017,SETPaper,Brown2018}. A complete understanding of defects, both invertible and noninvertible, is therefore necessary if we are to utilize topological phases to their fullest. In particular, such an understanding should allow for the computation of fusion of general defects. In this paper, we study non-chiral, two-dimensional, long-range-entangled topological phases with general defect structures. 

A defect of a topological phase is a region of positive codimension which differs from the ground state of the underlying bulk phase. Much work has been done on defects in topological phases, for example \onlinecites{Dennis2002,0610082,Bombin2007a,Bombin2010,MR2942952,FUCHS2002353,MR3370609,Kong2013,Brown2013a,MR3063919,Barkeshli2013,Barkeshli2014,Pastawski2014,1508.00573,Yoshida2015a,1606.07116,Brown2016,IrisCong1,IrisCong2,PhysRevB.96.195129,Yoshida2017,PhysRevB.96.125104,SETPaper,Bridgeman2017,Brown2018,1809.00245,1810.10539}. Although a complete classification for defects exists, it is not computational in nature. In this work, we show how the compound defect associated to a defect network\cite{1810.10539} can be computed. Our techniques is not restricted to invertible bimodules or defects.

In \onlinecite{Levin2005}, Levin and Wen (LW) constructed a long range entangled, 2D topological phase of matter associated to any fusion category \cat{C}. When $\cat{C} = \vvec{G}$ for a finite group $G$, this phase agrees with Kitaev's quantum double phase defined in \onlinecite{MR1951039}. For small groups, these Kitaev models are currently of great experimental interest\cite{chow2014implementing,Gambetta1}.
In \onlinecite{MR2942952}, Kitaev and Kong demonstrated that $\cat{C}{-}\cat{D}$ bimodules correspond to domain walls between the corresponding LW phases. 
In \onlinecite{1806.01279}, we showed how to compute the tensor product of $\cat{C}{-}\cat{D}$ bimodules $\cat{M}\otimes_{\cat{D}}\cat{N}$, corresponding to fusing the domain walls in the LW model. Additionally, we gave an explicit physical interpretation of all bimodules for the case $\cat{C}=\cat{D}=\vvec{\ZZ{p}}$ for prime $p$. In \onlinecite{1810.09469}, we extended this work to include \emph{binary interface defects}. We showed how to compute the horizontal fusion (tensor product) and vertical fusion (composition) of these defects. In the case $\cat{C}=\cat{D}=\vvec{\ZZ{p}}$, we provided complete fusion tables and physical interpretations of all binary interface defects.

This paper is a continuation of the work from \onlinecites{1806.01279,1810.09469}. 
In this paper, we present a new procedure, which we call the {\bf domain wall structure algorithm}, that computes the compound point defect associated to a domain wall structure once the holes have been filled in with point defects. We use this algorithm to show that all \emph{domain wall associators} are trivial for $\cat{C}=\cat{D}=\vvec{\ZZ{p}}$. This was not previously known for the noninvertible bimodules over $\vvec{\ZZ{p}}$. As we discuss below, these associators are related to the $O_3$ obstruction\cite{MR2677836}. This obstruction plays an important role in gauging as explained in \onlinecites{1804.01657,MR3555361}. From the condensed matter perspective, it would be interesting to find some non-trivial bimodule associators since it would show how a, potentially complicated, defect could be produced from simple domain walls.

The remainder of this paper is structured as follows. 
In Section~\ref{sec:prelim} we provide some definitions and preliminaries that are required for the remainder of the manuscript. In Section~\ref{sec:domain_wall_structure_algo}, we explain all the data which is required to execute the algorithm and how to execute it by hand. The real virtue of Algorithm~\ref{alg} is that it can be implemented in a computer. The authors have used the algorithm to check the horizontal and vertical fusion tables from \onlinecite{1810.09469} in the computer. 
In Section~\ref{sec:bimodule_associators}, we discuss the domain wall associator. We provide an example calculations and physical interpretations. In Section~\ref{sec:defectLW}, we define a lattice model which generalizes the Levin-Wen Lattice model by introducing domain walls and dimension 0 defects. The ground states of this lattice model can be computed using the {\em domain wall structure algorithm}. We conclude in Section~\ref{sec:conclusions}.

In Appendix~\ref{sec:categoryreps}, we recall some important definitions. We also provide the defining data for the $\vvec{\ZZ{p}}$ bimodules (reproduced from \onlinecite{1806.01279}) and their physical interpretations. Appendix~\ref{sec:extension_theory} briefly discusses the relationship between the bimodule associator structure and \emph{extension theory} of a fusion category.
In Appendix~\ref{sec:representation_tables}, we provide tables defining the irreducible representations of the annular categories for $\vvec{\ZZ{p}}$. 
We provide the complete set of bimodule associator defects in Appendix~\ref{sec:bimodule_associator_tables}. 
In the auxiliary material \footnote{\href{https://arxiv.org/src/1901.08069/anc/}{The auxiliary material can be found by going to the `ancillary files' link in the right hand column of the arXiv page for this paper.}}, we provide a Mathematica notebook that computes the composition of binary interface defects.

% This line sets the project root file.
% !TEX root = ../defects_domain_wall_structures.tex
% !TEX spellcheck = en_US

\section{Preliminaries}\label{sec:prelim}

For definitions of fusion categories and bimodules we refer to \onlinecite{1806.01279} and Appendix~\ref{sec:categoryreps}. 

\begin{definition}[category representation] \label{def:category_representation}
	Let $\cat{A}$ be a category. A representation of $\cat{A}$ is a functor $V : \cat{A} \to \vvec{}$. This functor is comprised of a vector space $V_a$ for each object $a \in \cat{A}$ and a linear map $V_f : V_a \to V_b$ for each morphism $f : a \to b$. The linear maps must satisfy the equations $V_{f \circ g} = V_f V_g$ and $V_{\rm id} = {\rm id}$. A basis for a representation consists of a basis for each vector space $V_a$.	
\end{definition}

\begin{definition}[point defect]\label{def:annular_category}
Consider a sectioned annulus where the faces are labeled by fusion categories and the section dividers are labeled by bimodules, for example
  \begin{align}
    \begin{array}{c}
      \includeTikz{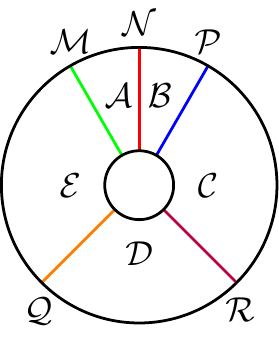}{
        \begin{tikzpicture}[scale=0.7]
          \draw[thick,red] (0,0) -- (0,2);
          \draw[thick,blue] (0,0) -- (1,1.73205);
          \draw[thick,green] (0,0) -- (-1,1.73205);
          \draw[thick,purple] (0,0) -- (1.41421,-1.41421);
          \draw[thick,orange] (0,0) -- (-1.41421,-1.41421);
          \draw[thick] (0,0) circle (2);
          \filldraw[thick,fill=white] (0,0) circle (0.5);
          \node[above] at (0,2) {$\cat{N}$};
          \node[above] at (1,1.73205) {$\cat{P}$};
          \node[above] at (-1,1.73205) {$\cat{M}$};
          \node[below,yshift=-1pt,xshift=1pt] at (1.41421,-1.41421) {$\cat{R}$};
          \node[below,yshift=-1pt,xshift=-1pt] at (-1.41421,-1.41421) {$\cat{Q}$};
          \node at (-0.3,1.3) {$\cat{A}$};
          \node at (0.3,1.3) {$\cat{B}$};
          \node at (1,0) {$\cat{C}$};
          \node at (-1,0) {$\cat{E}$};
          \node at (0,-1) {$\cat{D}$};
        \end{tikzpicture}
      }
    \end{array}.
    \end{align}
Associated to this, we define an {\bf annular category} (also known as a sphere category in \onlinecite{MR2978449}) whose objects are tuples of simple objects from the bimodules and morphisms are string diagrams which can be drawn in the annulus modulo isotopy and local relations
  \begin{align}
    \begin{array}{c}
      \includeTikz{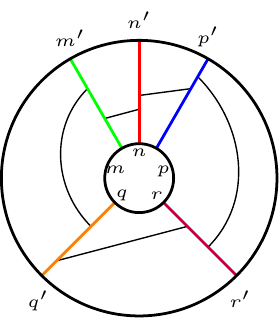}{
        \begin{tikzpicture}[scale=0.7]
          \draw (-0.5,0.866025) -- (0,1);
          \draw (0,1.2) -- (0.75,1.29904);
          \draw (0.85,1.47224) to [out=-45,in=45] (1,-1);
          \draw (-1.2,-1.2) -- (0.7,-0.7);
          \draw (-0.75,1.29904) to [out=180+45,in=90+45] (-0.7,-0.7);
          \draw[thick,red] (0,0) -- (0,2);
          \draw[thick,blue] (0,0) -- (1,1.73205);
          \draw[thick,green] (0,0) -- (-1,1.73205);
          \draw[thick,purple] (0,0) -- (1.41421,-1.41421);
          \draw[thick,orange] (0,0) -- (-1.41421,-1.41421);
          \draw[thick] (0,0) circle (2);
          \filldraw[thick,fill=white] (0,0) circle (0.5);
          \node[above] at (0,2) {\tiny $n'$};
          \node[above] at (1,1.73205) {\tiny $p'$};
          \node[above] at (-1,1.73205) {\tiny$m'$};
          \node[below,yshift=-1pt,xshift=1pt] at (1.41421,-1.41421) {\tiny$r'$};
          \node[below,yshift=-1pt,xshift=-1pt] at (-1.41421,-1.41421) {\tiny$q'$};
          \node[below,xshift=3pt] at (0.2,0.34641) {\tiny$p$};
          \node[below,xshift=-3pt] at (-0.2,0.34641) {\tiny$m$};
          \node at (0,0.37) {\tiny$n$};
          \node at (-0.25,-0.25) {\tiny$q$};
          \node at (0.25,-0.25) {\tiny$r$};
        \end{tikzpicture}
      }
    \end{array}.
  \end{align}
  In \onlinecite{MR2978449}, representations of the annular category are called sphere modules. Physically, these representations parameterize point defects at the domain wall junction, so we shall refer to them as {\bf point defects} or simply {\bf defects} when there is no ambiguity. It is important to not isotope the bimodule strings or rotate the annulus. In order to do this, we would need to specify a consistent set of rigid structures for all the bimodules involved. This can be done \cite{MR3435098}, but is beyond the scope of this work.
\end{definition}

\begin{definition}[domain wall structure] \label{def:domain_wall_structure}
A {\bf domain wall structure} consists of a graph embedded into a disc where the edges don't have critical points. We label the faces of the graph with fusion categories, and the edges of the graph with bimodules between the corresponding fusion categories. For example,
  \begin{align}
    \begin{array}{c}
      \includeTikz{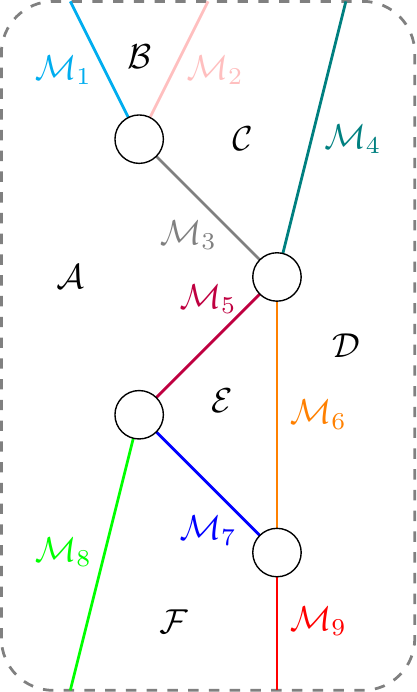}{
        \begin{tikzpicture}[scale=0.7,every node/.style={scale=1}]
          \draw[gray,thick,dashed,rounded corners=15pt] (0,0) rectangle (6,10);
          \draw[red,thick] (4,0) -- (4,2) node[midway,right] {$\cat{M}_9$};
          \draw[blue,thick] (4,2) -- (2,4) node[midway,below,yshift=-2mm] {$\cat{M}_7$};
          \draw[green,thick] (1,0) -- (2,4) node[midway,left] {$\cat{M}_8$};
          \draw[orange,thick] (4,2) -- (4,6) node[midway, right] {$\cat{M}_6$};
          \draw[purple,thick] (2,4) -- (4,6) node[midway, above,yshift=2mm] {$\cat{M}_5$};
          \draw[teal,thick] (4,6) -- (5,10) node[midway, right] {$\cat{M}_4$};
          \draw[gray,thick] (4,6) -- (2,8) node[midway, below,xshift=-2mm] {$\cat{M}_3$};
          \draw[cyan,thick] (2,8) -- (1,10) node[midway,left] {$\cat{M}_1$};
          \draw[pink,thick] (2,8) -- (3,10) node[midway, right] {$\cat{M}_2$};
          \filldraw[fill=white, draw=black] (4,2) circle (10pt);
          \filldraw[fill=white, draw=black] (2,4) circle (10pt);
          \filldraw[fill=white, draw=black] (4,6) circle (10pt);
          \filldraw[fill=white, draw=black] (2,8) circle (10pt);
          \node at (1,6) {$\cat{A}$};
          \node at (2.5,1) {$\cat{F}$};
          \node at (3.2,4.2) {$\cat{E}$};
          \node at (3.5,8) {$\cat{C}$};
          \node at (2,9.2) {$\cat{B}$};
          \node at (5,5) {$\cat{D}$};
        \end{tikzpicture}
        }\end{array}.
  \end{align}
\end{definition}
\begin{definition}[compound defect] \label{def:compount_defect}
	A {\bf compound defect} consists of a domain wall structure, along with an assignment of a point defect (annular category representation) to each vertex. For example,
\begin{align} \label{eq:compound_defect_ref}
	\begin{array}{c}
	\includeTikz{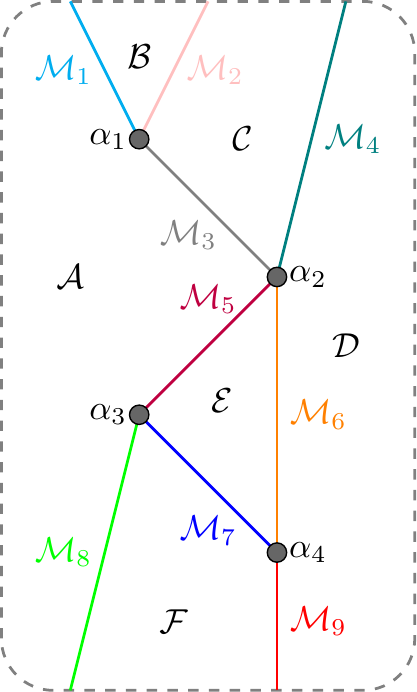}{
		\begin{tikzpicture}[scale=0.7,every node/.style={scale=1}]
		\draw[gray,thick,dashed,rounded corners=15pt] (0,0) rectangle (6,10);
		\draw[red,thick] (4,0) -- (4,2) node[midway,right] {$\cat{M}_9$};
		\draw[blue,thick] (4,2) -- (2,4) node[midway,below,yshift=-2mm] {$\cat{M}_7$};
		\draw[green,thick] (1,0) -- (2,4) node[midway,left] {$\cat{M}_8$};
		\draw[orange,thick] (4,2) -- (4,6) node[midway, right] {$\cat{M}_6$};
		\draw[purple,thick] (2,4) -- (4,6) node[midway, above,yshift=2mm] {$\cat{M}_5$};
		\draw[teal,thick] (4,6) -- (5,10) node[midway, right] {$\cat{M}_4$};
		\draw[gray,thick] (4,6) -- (2,8) node[midway, below,xshift=-2mm] {$\cat{M}_3$};
		\draw[cyan,thick] (2,8) -- (1,10) node[midway,left] {$\cat{M}_1$};
		\draw[pink,thick] (2,8) -- (3,10) node[midway, right] {$\cat{M}_2$};
		\filldraw[fill=dandark, draw=black] (4,2) circle (4pt) node[right]{$\alpha_4$};
		\filldraw[fill=dandark, draw=black] (2,4) circle (4pt) node[left]{$\alpha_3$};
		\filldraw[fill=dandark, draw=black] (4,6) circle (4pt) node[right]{$\alpha_2$};
		\filldraw[fill=dandark, draw=black] (2,8) circle (4pt) node[left]{$\alpha_1$};
		\node at (1,6) {$\cat{A}$};
		\node at (2.5,1) {$\cat{F}$};
		\node at (3.2,4.2) {$\cat{E}$};
		\node at (3.5,8) {$\cat{C}$};
		\node at (2,9.2) {$\cat{B}$};
		\node at (5,5) {$\cat{D}$};
		\end{tikzpicture}
	}\end{array},
\end{align}
with $\alpha_1,\ldots,\alpha_4$ representations of the appropriate annular categories. Diagram~\ref{eq:compound_defect_ref} defines an (generally decomposable) annular category representation. The vectors in this representation are constructed by choosing vectors from the representations $\alpha_1,\ldots,\alpha_4$ subject to consistent labeling of the edges with bimodule objects
\begin{align} \label{eq:compound_defect_vec}
	\begin{array}{c}
	\includeTikz{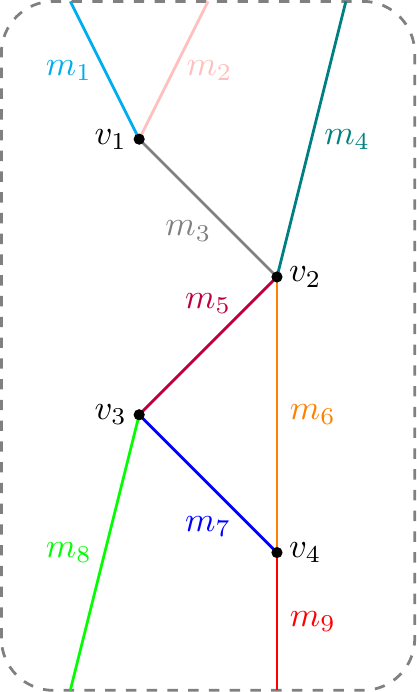}{
		\begin{tikzpicture}[scale=0.7,every node/.style={scale=1}]
		\draw[gray,thick,dashed,rounded corners=15pt] (0,0) rectangle (6,10);
		\draw[red,thick] (4,0) -- (4,2) node[midway,right] {$m_9$};
		\draw[blue,thick] (4,2) -- (2,4) node[midway,below,yshift=-2mm] {$m_7$};
		\draw[green,thick] (1,0) -- (2,4) node[midway,left] {$m_8$};
		\draw[orange,thick] (4,2) -- (4,6) node[midway, right] {$m_6$};
		\draw[purple,thick] (2,4) -- (4,6) node[midway, above,yshift=2mm] {$m_5$};
		\draw[teal,thick] (4,6) -- (5,10) node[midway, right] {$m_4$};
		\draw[gray,thick] (4,6) -- (2,8) node[midway, below,xshift=-2mm] {$m_3$};
		\draw[cyan,thick] (2,8) -- (1,10) node[midway,left] {$m_1$};
		\draw[pink,thick] (2,8) -- (3,10) node[midway, right] {$m_2$};
		\filldraw[fill=black, draw=black] (4,2) circle (2pt) node[right]{$v_4$};
		\filldraw[fill=black, draw=black] (2,4) circle (2pt) node[left]{$v_3$};
		\filldraw[fill=black, draw=black] (4,6) circle (2pt) node[right]{$v_2$};
		\filldraw[fill=black, draw=black] (2,8) circle (2pt) node[left]{$v_1$};
		\end{tikzpicture}
	}\end{array}.
\end{align}
The annular category action on a vector in a compound defect is %computed as follows
\begin{align}
\begin{array}{c}
\includeTikz{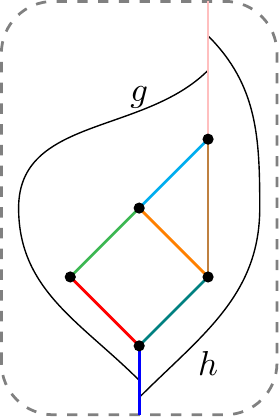}
	{
	\begin{tikzpicture}[scale=0.7,every node/.style={scale=1}]
		\draw[gray,thick,dashed,rounded corners=15pt] (-2,-3) rectangle (2,3);
		\draw (0,-2.5)to[out=135,in=270](-1.75,0)to[out=90,in=225] (1,2);\node at(0,1.6) {$g$};
		\draw (0,-2.75)to[out=45,in=270](1.75,0)to[out=90,in=-45] (1,2.5);\node at(1,-2.25) {$h$};
		\draw[blue,thick] (0,-3) -- (0,-2);
		\draw[red,thick] (0,-2) -- (-1,-1);
		\draw[nicegreen,thick] (-1,-1)--(0,0);
		\draw[teal,thick] (0,-2)--(1,-1);
		\draw[orange,thick] (1,-1)--(0,0);
		\draw[brown,thick] (1,-1)--(1,1);
		\draw[cyan,thick] (0,0)--(1,1);
		\draw[pink,thick] (1,1)--(1,3);
		\filldraw[fill=black, draw=black] (0,-2) circle (2pt);
		\filldraw[fill=black, draw=black] (-1,-1) circle (2pt);
		\filldraw[fill=black, draw=black] (1,-1) circle (2pt);
		\filldraw[fill=black, draw=black] (0,0) circle (2pt);
		\filldraw[fill=black, draw=black] (1,1) circle (2pt);
	\end{tikzpicture}
	}
\end{array}
&=
\begin{array}{c}
\includeTikz{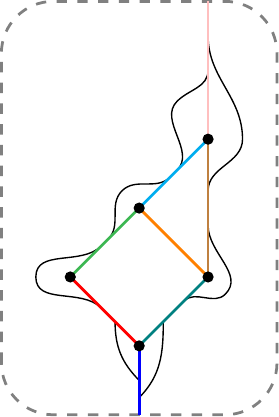}
{
	\begin{tikzpicture}[scale=0.7,every node/.style={scale=1}]
	\draw[gray,thick,dashed,rounded corners=15pt] (-2,-3) rectangle (2,3);
	\draw (0,-2.5)to[out=135,in=270](-.35,-1.65);
	\draw (-.55,-1.45)to[out=135,in=270](-1.5,-1)to[out=90,in=180+45](-.55,-.55);
	\draw (-.45,-.45)to[out=45,in=180+45](-.25,.25)to[out=45,in=180+45] (.45,.45);
	\draw (.55,.55) to[out=45,in=180+45] (.55,1.55) to[out=45,in=270] (1,2) ;
	\draw (0,-2.75)to[out=45,in=270](.35,-1.65);
	\draw (.65,-1.35) to[out=45,in=180+45](1.25,-1.25)to[out=45,in=270](1,-.25);
	\draw (1,.25)to[out=90,in=270] (1.5,1)to[out=90,in=270](1,2.5);
	\draw[blue,thick] (0,-3) -- (0,-2);
	\draw[red,thick] (0,-2) -- (-1,-1);
	\draw[nicegreen,thick] (-1,-1)--(0,0);
	\draw[teal,thick] (0,-2)--(1,-1);
	\draw[orange,thick] (1,-1)--(0,0);
	\draw[brown,thick] (1,-1)--(1,1);
	\draw[cyan,thick] (0,0)--(1,1);
	\draw[pink,thick] (1,1)--(1,3);
	\filldraw[fill=black, draw=black] (0,-2) circle (2pt);
	\filldraw[fill=black, draw=black] (-1,-1) circle (2pt);
	\filldraw[fill=black, draw=black] (1,-1) circle (2pt);
	\filldraw[fill=black, draw=black] (0,0) circle (2pt);
	\filldraw[fill=black, draw=black] (1,1) circle (2pt);
	\end{tikzpicture}
}
\end{array}.
\end{align}
There is also a bubble action for each internal cavity
\begin{align}
\begin{array}{c}
\includeTikz{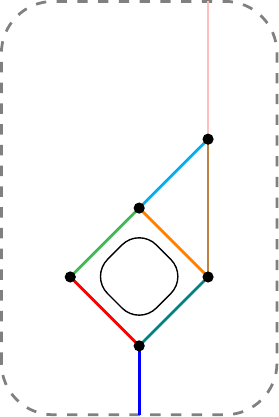}
	{
	\begin{tikzpicture}[scale=0.7,every node/.style={scale=1}]
		\draw[gray,thick,dashed,rounded corners=15pt] (-2,-3) rectangle (2,3);
		\draw[blue,thick] (0,-3) -- (0,-2);
		\draw[red,thick] (0,-2) -- (-1,-1);
		\draw[nicegreen,thick] (-1,-1)--(0,0);
		\draw[teal,thick] (0,-2)--(1,-1);
		\draw[orange,thick] (1,-1)--(0,0);
		\draw[brown,thick] (1,-1)--(1,1);
		\draw[cyan,thick] (0,0)--(1,1);
		\draw[pink,thick] (1,1)--(1,3);
    \draw[rounded corners=7pt,rotate around = {45:(0,-1.70)}] (0,-1.70) rectangle (1,-0.70);
		\filldraw[fill=black, draw=black] (0,-2) circle (2pt);
		\filldraw[fill=black, draw=black] (-1,-1) circle (2pt);
		\filldraw[fill=black, draw=black] (1,-1) circle (2pt);
		\filldraw[fill=black, draw=black] (0,0) circle (2pt);
		\filldraw[fill=black, draw=black] (1,1) circle (2pt);
	\end{tikzpicture}
	}
\end{array}
&=
\begin{array}{c}
\includeTikz{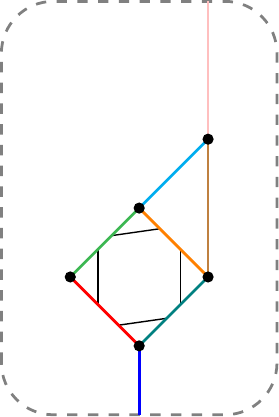}
{
	\begin{tikzpicture}[scale=0.7,every node/.style={scale=1}]
	\draw[gray,thick,dashed,rounded corners=15pt] (-2,-3) rectangle (2,3);
  \draw (-0.6,-1.4) -- (-0.6,-0.6);
  \draw (0.6,-1.4) -- (0.6,-0.6);
  \draw (-0.4,-0.4) -- (0.3,-0.3);
  \draw (-0.3,-1.7) -- (0.4,-1.6);
	\draw[blue,thick] (0,-3) -- (0,-2);
	\draw[red,thick] (0,-2) -- (-1,-1);
	\draw[nicegreen,thick] (-1,-1)--(0,0);
	\draw[teal,thick] (0,-2)--(1,-1);
	\draw[orange,thick] (1,-1)--(0,0);
	\draw[brown,thick] (1,-1)--(1,1);
	\draw[cyan,thick] (0,0)--(1,1);
	\draw[pink,thick] (1,1)--(1,3);
	\filldraw[fill=black, draw=black] (0,-2) circle (2pt);
	\filldraw[fill=black, draw=black] (-1,-1) circle (2pt);
	\filldraw[fill=black, draw=black] (1,-1) circle (2pt);
	\filldraw[fill=black, draw=black] (0,0) circle (2pt);
	\filldraw[fill=black, draw=black] (1,1) circle (2pt);
	\end{tikzpicture}
}
\end{array}.
\end{align}
We must quotient away the bubble actions for every internal cavity because bubbles internal to a cavity should evaluate to the dimension of their labeling object. For ${\bf Vec}(\mathbb{Z}/p\mathbb{Z})$ every simple object has dimension 1. This bubble action appears in the definition of the Levin-Wen Hamiltonian from \onlinecite{Levin2005}, and the version presented in Section~\ref{sec:defectLW}.
\end{definition}

\begin{algorithm}[domain wall structure algorithm] \label{alg}
	
	The \emph{domain wall structure algorithm} decomposes a compound defect into a sum of simple representations.
	
The main steps in the domain wall structure algorithm are as follows:
\begin{enumerate}
\item Construct a compound defect by filling the holes in the domain wall structure with vectors from the corresponding annular category representations, subject to the labels on the internal edges agreeing.
\item Quotient out the bubble action for each internal cavity
\item Compute all relevant idempotent actions on the quotient representation. This lets us decompose the quotient representation into simple annular category representations.
\end{enumerate}
\end{algorithm}

% This line sets the project root file.
% !TEX root = ../defects_domain_wall_structures.tex
% !TEX spellcheck = en_US

\section{The Domain Wall Structure Algorithm} \label{sec:domain_wall_structure_algo}

The goal of this section is to explain how to compute the compound defect. We shall demonstrate how the computation works using $\vvec{\ZZ{p}}$ as our central example, but everything we describe works in much more generality. 

All the annular categories ${\bf Ann}$ of interest in this paper are semi-simple, so we can describe their representations as functors ${\bf Ann} \to \vvec{}$ or as indecomposable idempotent endomorphisms in ${\bf Ann}$. In \onlinecite{1810.09469}, we exclusively used the idempotent description. In this paper, we shall use both ways of presenting a representation of ${\bf Ann}$. 

Given an indecomposable idempotent $i:a\to a$ in ${\bf Ann}$, the corresponding functor $V:{\bf Ann} \to \vvec{}$ is defined by $V_x={\bf Ann}(a,x) \circ i$. This vector space is the image of the projection
\begin{align}
  {\bf Ann}(a,x) &\to {\bf Ann}(a,x) \\
  f &\mapsto f \circ i.
\end{align}
The vector space $V_x$ is nontrivial exactly when there is a nontrivial morphism $a\to x$ of the form $f\circ i$.

All of the functors corresponding to binary interface defects described in \onlinecite{1810.09469} have been tabulated in Appendix~\ref{sec:representation_tables}. The vector space in which the tabulated vectors live can be read off from the string labels.
We refer to \onlinecite{1806.01279} for the bimodule definitions for $\vvec{\ZZ{p}}$. Definitions of idempotents corresponding to all 2-string annular categories can be found in \onlinecite{1810.09469}. We will use the notation (introduced in \onlinecite{1810.09469}) $\defect{M}{N}{x}{}{}$ for a defect called `$x$' occurring at the interface of a bimodule $M$ and a bimodule $N$.

We shall now demonstrate how an entry of the representation tables (Appendix~\ref{sec:representation_tables}) are computed. 
\begin{example}[Constructing irreducible representations]
  Consider the defect $\defect{R}{F_r}{x}{}{}$ which was defined in \onlinecite{1810.09469} by the idempotent
\begin{align}
\defect{R}{F_r}{x}{}{}=\frac{1}{p}\sum_k\omega^{kx}
\begin{array}{c}
	\includeTikz{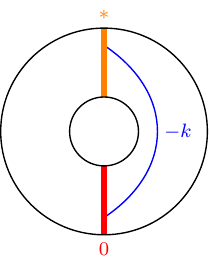}
	{
		\begin{tikzpicture}[scale=.7,every node/.style={scale=.7}]
		\annparamst{$0$}{$*$}{}{};
		\annst{}{$-k$};
		\end{tikzpicture}
	}
\end{array}.
\end{align}
This idempotent serves two purposes. Firstly, it labels an irreducible representation of ${\bf Ann}_{R,F_r}$. Secondly, the idempotent projects onto the representation it labels. If $f : (0;*) \to (m;*)$, then $f \mapsto f \circ \defect{F_r}{R}{x}{}{}$ is an endomorphism of ${\bf Ann}_{R,F_r}((0;*),(m;*))$. We choose the following basis for the image of this endomorphism
\begin{align}
\begin{array}{c}\includeTikz{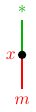}{
					\begin{tikzpicture}[scale=.7,every node/.style={scale=.5}]
					\stbivertex{m}{*}{x}{};
					\end{tikzpicture}}
\end{array}
:=
\frac{1}{p}\sum_k\omega^{k(x+rm)}
	\begin{array}{c}
				\includeTikz{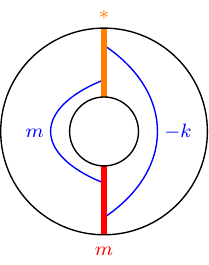}
				{
					\begin{tikzpicture}[scale=.7,every node/.style={scale=.7}]
					\annparamst{$m$}{$*$}{}{};
					\annst{$m$}{$-k$};
					\end{tikzpicture}
				}
	\end{array}.
\end{align}
Acting by a general morphism
\begin{align}
  \begin{array}{c}
    \includeTikz{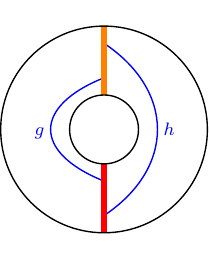}{
    \begin{tikzpicture}[scale=.7,every node/.style={scale=.65}]
	\annparamss{}{}{}{};
	\annst{$g$}{$h$};
\end{tikzpicture}}
\end{array}
\end{align}
on the basis vectors gives
\begin{align}
 \omega^{h(x+r(g+m))}
  \begin{array}{c}\includeTikz{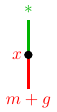}{
					\begin{tikzpicture}[scale=.7,every node/.style={scale=.5}]
					\stbivertex{m+g}{*}{x}{};
					\end{tikzpicture}}
\end{array},
\end{align}
as tabulated in Appendix~\ref{sec:representation_tables}.
\end{example}

\subsection{Vertical defect fusion}
The simplest case of the domain wall structure algorithm is vertical defect fusion, corresponding to the domain wall structure
\begin{align}
  \begin{array}{c}
    \includeTikz{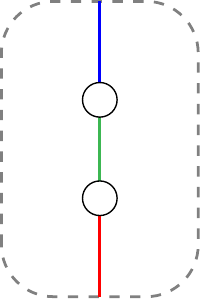}{
      \begin{tikzpicture}[scale=0.5]
        \draw[gray,thick,dashed,rounded corners=15pt] (0,0) rectangle (4,6);
        \draw[thick,red] (2,0) -- (2,2);
        \draw[thick,nicegreen] (2,2) -- (2,4);
        \draw[thick,blue] (2,4) -- (2,6);
        \filldraw[fill=white, draw=black] (2,2) circle (10pt);
        \filldraw[fill=white, draw=black] (2,4) circle (10pt);
        \end{tikzpicture}
      }
    \end{array}.\label{eqn:composition_structure}
\end{align}
In \onlinecite{1810.09469}, we computed these vertical defect fusions for all compatible pairs of defects in the $\vvec{\ZZ{p}}$ model. These vertical fusions can also be computed using the domain wall structure algorithm. 
%The basic idea is that we fill all the holes in the domain wall structure with vectors from the corresponding representations, subject to consistent labeling of the internal edges. 
Given a pair of point defects $\alpha_1,\alpha_2$ (equivalently representations of 2-string annular categories), the compound defect is formed by filling the holes in Eqn.~\ref{eqn:composition_structure} with these defects
\begin{align} V =
\begin{array}{c}
\includeTikz{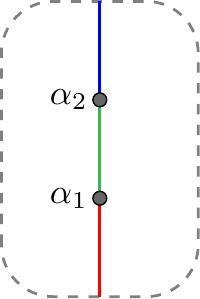}{
	\begin{tikzpicture}[scale=0.5]
	\draw[gray,thick,dashed,rounded corners=15pt] (0,0) rectangle (4,6);
	\draw[thick,red] (2,0) -- (2,2);
	\draw[thick,nicegreen] (2,2) -- (2,4);
	\draw[thick,blue] (2,4) -- (2,6);
	\filldraw[fill=dandark, draw=black] (2,2) circle (4pt) node[left] {$\alpha_1$};
	\filldraw[fill=dandark, draw=black] (2,4) circle (4pt) node[left] {$\alpha_2$};
	\end{tikzpicture}
}
\end{array}.\label{eqn:composition_structure_filled}
\end{align}
This forms a (possibly reducible) annular category representation. A vector in this representation looks like
\begin{align}
\begin{array}{c}
\includeTikz{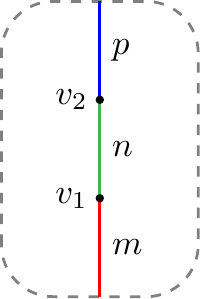}{
	\begin{tikzpicture}[scale=0.5]
	\draw[gray,thick,dashed,rounded corners=15pt] (0,0) rectangle (4,6);
	\draw[thick,red] (2,0) -- (2,2);
	\draw[thick,nicegreen] (2,2) -- (2,4);
	\draw[thick,blue] (2,4) -- (2,6);
  \node[right] at (2,1) {$m$};
  \node[right] at (2,3) {$n$};
  \node[right] at (2,5) {$p$};
	\filldraw[fill=black, draw=black] (2,2) circle (2pt) node[left] {$v_1$};
	\filldraw[fill=black, draw=black] (2,4) circle (2pt) node[left] {$v_2$};
	\end{tikzpicture}
}
\end{array} \in V_{(m,p)}.
\end{align}
If $\alpha$ is a binary interface defect and $i_{\alpha} : (m_{\alpha},n_{\alpha}) \to (m_{\alpha},n_{\alpha})$ is the corresponding idempotent from \cite{1810.09469}, then we have
\begin{align} \label{eq:isotypic}
  V \cong \bigoplus_{\alpha} \dim (i_{\alpha}V_{(m_{\alpha},n_{\alpha})}) \cdot \alpha.
\end{align}
In representation theory, this is called an isotypic decomposition. The general theory of isotypic decompositions is explained in chapter 4 of \onlinecite{MR2522486}. We use Eqn.~\ref{eq:isotypic} to decompose $V$ into irreducible representations.

\begin{example}
  Consider the vertical defect fusion $\defect{R}{F_r}{x}{}{} \circ \defect{F_r}{R}{z}{}{}$. We begin by building the compound representation of the annular category ${\bf Ann}_{R,R}$ from our chosen defects. It has basis vectors of the form
  \begin{align}
    \begin{array}{c}
	\includeTikz{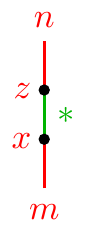}{
        \begin{tikzpicture}
\vertcompcompoundrep{m}{x}{*}{z}{n}
          \end{tikzpicture}
        }
      \end{array}.
    \end{align}
 In \onlinecite{1810.09469}, we define the idempotent
\begin{align}
	\defect{R}{R}{\alpha}{\zeta}{}=\frac{1}{p}\sum_\gamma\omega^{\gamma\zeta}
	\begin{array}{c}
		\includeTikz{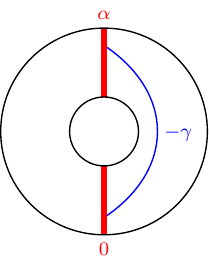}
		{
			\begin{tikzpicture}[scale=.7,every node/.style={scale=.7}]
			\annparamss{$0$}{$\alpha$}{}{};
			\annss{}{$-\gamma$};
			\end{tikzpicture}
		}
	\end{array}.
\end{align}
Recall that this idempotent projects onto the irreducible representation labeled by $(\alpha,\zeta)$. Applying this idempotent to the basis vector
\begin{align}
	\begin{array}{c}
	\includeTikz{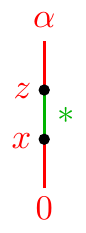}{
        \begin{tikzpicture}
\vertcompcompoundrep{0}{x}{*}{z}{\alpha}
          \end{tikzpicture}
        }
      \end{array}
\end{align}
  gives
  \begin{align}
    \frac{1}{p}\sum_{\gamma}\omega^{\gamma(\zeta-z-x+r\alpha)}
    \begin{array}{c}
      \includeTikz{vertical_composition_rep_2}{
        \begin{tikzpicture}
\vertcompcompoundrep{0}{x}{*}{z}{\alpha}
          \end{tikzpicture}
        }
      \end{array}
  \end{align}
  which is zero unless $\zeta = x + z - r\alpha$. Therefore we have
\begin{align}
\defect{R}{F_r}{x}{}{} \circ \defect{F_r}{R}{z}{}{} \cong \oplus_{\alpha} \defect{R}{R}{\alpha}{x+z-r\alpha}{}.
\end{align}
 \end{example}

\subsection{Horizontal defect fusion}
If we only use annular categories with two bimodule strings, the domain wall structure algorithm only computes vertical composition of defects. To compute more interesting compound defects, we need to include annular categories with three or more bimodule strings. In \onlinecite{1806.01279}, we computed the Brauer-Picard ring for the fusion category $\vvec{\ZZ{p}}$. More precisely, for all pairs of $\vvec{\ZZ{p}}$ bimodules $\cat{M},\cat{N}$, we computed an explicit isomorphism $\cat{M} \otimes_{\ZZ{p}} \cat{N} \cong \oplus_i \cat{P}_i$. These explicit isomorphisms are recorded in the inflation tables in \onlinecite{1810.09469}. If we take the identity (under vertical fusion) defect on $M$ and inflate the top or bottom part, we get an idempotent in a three string annular category. The corresponding representations play the role of bimodule trivalent vertices. These representations have been tabulated in Appendix~\ref{sec:representation_tables}. Now we demonstrate how to compute an entry of this table.

\begin{example}
  Consider the trivial defect on the $X_x$ domain wall
  \begin{align}
  \defect{X_x}{X_x}{0}{0}{}=\frac{1}{p}\sum_g
	\begin{array}{c}
		\includeTikz{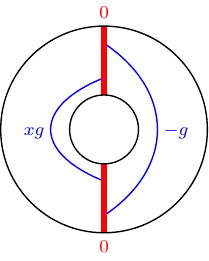}{
			\begin{tikzpicture}[scale=.7,every node/.style={scale=.65}]
				\annparamss{$0$}{$0$}{}{};
				\annss{$xg$}{$-g$};
			\end{tikzpicture}
		}
	\end{array}.
	\end{align}
          If we inflate the top half of this idempotent along the isomorphism $X_x \cong F_q \otimes_{\ZZ{p}} F_r$ where $x = q^{-1}r$, then we get the idempotent
 \begin{align}
   \frac{1}{p^2} \sum_{g,k}
   \begin{array}{c}
     \includeTikz{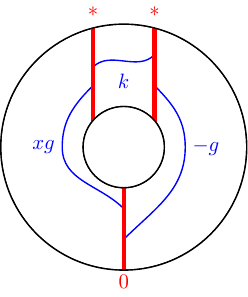}{
       \begin{tikzpicture}[scale=1.25,every node/.style={scale=.6}]
         \onedowntwouptube{$xg$}{$-g$}{$k$}{$0$}{$*$}{$*$}
       \end{tikzpicture}
     }
     \end{array}.
 \end{align}
Composing morphisms $(0;*,*) \to (m;*,*)$ on the outside gives us a linear endomorphism of ${\bf Ann}_{F_q,F_r,X_x}((0;*,*),(m;*,*))$. We choose the following basis for the image of this endomorphism:
\begin{align}
	 \begin{array}{c}\includeTikz{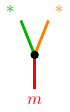}{
	\begin{tikzpicture}[scale=.7,every node/.style={scale=.5}]
	\rsttrivertexdual{*}{*}{m}{};
	\end{tikzpicture}
	 }\end{array} := \frac{1}{p^2} \sum_{g,k} \omega^{qmk}
	 \begin{array}{c}
	   \includeTikz{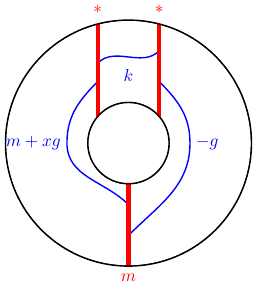}{
	     \begin{tikzpicture}[scale=1.25,every node/.style={scale=.5}]
	       \onedowntwouptube{$m+xg$}{$-g$}{$k$}{$m$}{$*$}{$*$}
	     \end{tikzpicture}
	   }
	   \end{array}
\end{align}
This forms the basis for our representation. Applying
\begin{align}
\begin{array}{c}
  \includeTikz{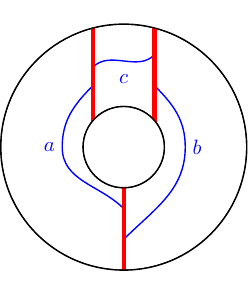}
{
	\begin{tikzpicture}[scale=1.25,every node/.style={scale=.6}]
\onedowntwouptube{$a$}{$b$}{$c$}{}{}{}
	\end{tikzpicture}	
}
\end{array}.
\end{align}
and making the substitutions $g \to g + b, k \to k - c$ gives
\begin{align}
\omega^{-c (q (a+m + xb))}\begin{array}{c}\includeTikz{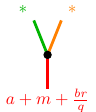}{
\begin{tikzpicture}[scale=.7,every node/.style={scale=.5}]
\rsttrivertexdual{*}{*}{a+m+\frac{b r}{q}}{};
\end{tikzpicture}
}\end{array},
\end{align}
as recorded in Table~\ref{tab:dualtrivertextable}.
\end{example}
Now that we have a collection of 2 and 3 bimodule string annular category representations at our disposal, we can discuss some more complicated domain wall structures and compute the corresponding compound defects. Of particular interest is the domain wall structure
\begin{align}
  \begin{array}{c}
    \includeTikz{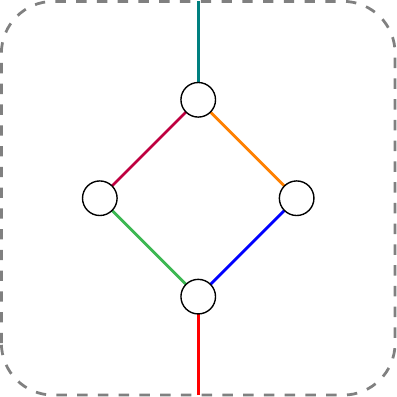}{
      \begin{tikzpicture}[scale=0.5]
        \draw[gray,thick,dashed,rounded corners=15pt] (0,0) rectangle (8,8);
        \draw[thick,red] (4,0) -- (4,2);
        \draw[thick,nicegreen] (4,2) -- (2,4);
        \draw[thick,blue] (4,2) -- (6,4);
        \draw[thick,purple] (2,4) -- (4,6);
        \draw[thick,orange] (6,4) -- (4,6);
        \draw[thick,teal] (4,6) -- (4,8);
        \filldraw[fill=white, draw=black] (4,2) circle (10pt);
        \filldraw[fill=white, draw=black] (2,4) circle (10pt);
        \filldraw[fill=white, draw=black] (6,4) circle (10pt);
        \filldraw[fill=white, draw=black] (4,6) circle (10pt);
        \end{tikzpicture}
      }
    \end{array}.
\end{align}
This domain wall structure corresponds to horizontal defect fusion. In the $\vvec{\ZZ{p}}$ case, we computed all possible horizontal defect fusions in \onlinecite{1810.09469}. In the following example, we demonstrate how to compute horizontal defect fusion using the domain wall structure algorithm. This example is the first time we encounter the internal cavity {\em bubble action}, which we need to trivialize to get the correct answer.

\begin{example} \label{ex:hoz_fusion_1}
  Consider the horizontal fusion $\defect{F_q}{R}{x}{}{} \otimes \defect{L}{L}{c}{z}{}$. Using the trivalent vertices corresponding to the isomorphisms $R \otimes_{\ZZ{p}} L \cong p \cdot T$ and $F_q \otimes_{\ZZ{p}} L \cong T$, we can construct a (reducible) representation of the category ${\bf Ann}_{T,T}$. It has the basis
\begin{align}
  \begin{array}{c}
    \includeTikz{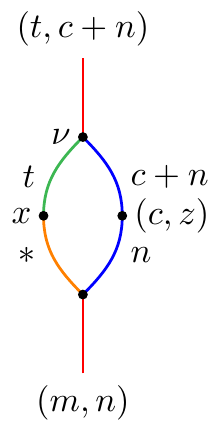}{
      \begin{tikzpicture}[scale = 0.4]
        \diamondvectorexampleone{$(m,n)$}{$(t,c+n)$}{$x$}{$(c,z)$}{$*$}{$n$}{$c+n$}{$t$}{$\nu$};
       \end{tikzpicture}
      }
    \end{array}.
\end{align}
This representation is too large. It has a $\ZZ{p}$ action by introducing a bubble into the middle cavity. In order to get a physically relevant representation, we need to quotient away this action to construct the representation of interest. Acting by a $g$ bubble multiplies the above vector by $\omega^{g(x+z-\nu-q(t+m))}$
\begin{align}
\begin{array}{c}
\includeTikz{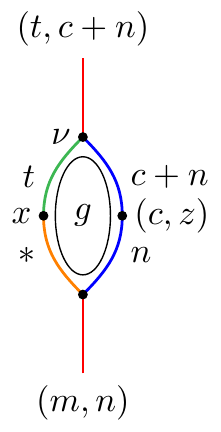}{
	\begin{tikzpicture}[scale = 0.4]
	\diamondvectorexampleone{$(m,n)$}{$(t,c+n)$}{$x$}{$(c,z)$}{$*$}{$n$}{$c+n$}{$t$}{$\nu$};
	\draw (0,4) ellipse (.7 and 1.5);\node at (0,4) {$g$};
	\end{tikzpicture}
}
\end{array}
&=
\begin{array}{c}
\includeTikz{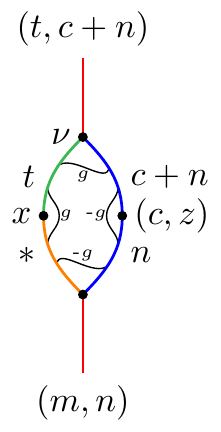}{
	\begin{tikzpicture}[scale = 0.4]
	\draw (-.9,3.25) to[out=90,in=270](-.6,4)to[out=90,in=270](-.9,4.75);\node[] at (-.6,4) {\tiny $\,\,\,g$};
	\draw (.9,3.25) to[out=90,in=270](.6,4)to[out=90,in=270](.9,4.75);\node[] at (.5,4) {\tiny -$g\,\,\,$};
	\draw (-.66,5.25)to[out=45,in=270](.66,5.25); \node at (0,5) {\tiny$g$};
	\draw (-.66,2.75)to[out=90,in=225](.66,2.75); \node at (0,3) {\tiny-$g$};
	\diamondvectorexampleone{$(m,n)$}{$(t,c+n)$}{$x$}{$(c,z)$}{$*$}{$n$}{$c+n$}{$t$}{$\nu$};
	\end{tikzpicture}
}
\end{array}
=\omega^{gx}\omega^{g(z-qt)}
\begin{array}{c}
\includeTikz{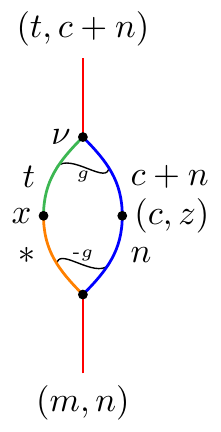}{
	\begin{tikzpicture}[scale = 0.4]
	\draw (-.66,5.25)to[out=45,in=270](.66,5.25); \node at (0,5) {\tiny$g$};
	\draw (-.66,2.75)to[out=90,in=225](.66,2.75); \node at (0,3) {\tiny-$g$};
	\diamondvectorexampleone{$(m,n)$}{$(t,c+n)$}{$x$}{$(c,z)$}{$*$}{$n$}{$c+n$}{$t$}{$\nu$};
	\end{tikzpicture}
}
\end{array}
=\omega^{g(x+z-qt)}\omega^{gqm}\omega^{-g\nu}
\begin{array}{c}
\includeTikz{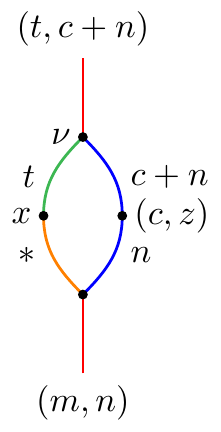}{
	\begin{tikzpicture}[scale = 0.4]
	\diamondvectorexampleone{$(m,n)$}{$(t,c+n)$}{$x$}{$(c,z)$}{$*$}{$n$}{$c+n$}{$t$}{$\nu$};
	\end{tikzpicture}
}
\end{array}.
\end{align}
Therefore, unless $t = q^{-1}(x+z-\nu) + m$, the vector projects onto zero in the quotient. After taking the quotient, the idempotent $\defect{T}{T}{\alpha}{\beta}{}$ acts as zero unless $\alpha = q^{-1}(x+z-\nu)$ and $\beta = c$. This is exactly the horizontal fusion outcome $\defect{F_q}{R}{x}{}{} \otimes \defect{L}{L}{c}{z}{} \cong \defect{T}{T}{q^{-1}(x+z-\nu)}{c}{\nu}$ which was computed in \onlinecite{1810.09469}.
  \end{example}

\begin{example}
Consider the horizontal fusion $\defect{X_k}{X_l}{}{}{} \otimes \defect{F_0}{R}{z}{}{}$. As in Example~\ref{ex:hoz_fusion_1}, we construct a representation of the category ${\bf Ann}_{F_0,T}$. It has the basis
\begin{align}
  \begin{array}{c}
    \includeTikz{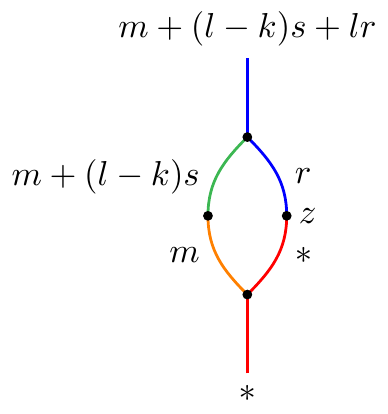}{
      \begin{tikzpicture}[scale = 0.4]
        \diamondvectorexampletwo{$*$}{$m+(l-k)s+lr$}{}{$z$}{$m$}{$*$}{$r$}{$m+(l-k)s$};
       \end{tikzpicture}
      }
    \end{array}.
\end{align}
Acting by a bubble labeled with $r$ in the internal cavity sends this vector to
\begin{align}
  \begin{array}{c}
    \includeTikz{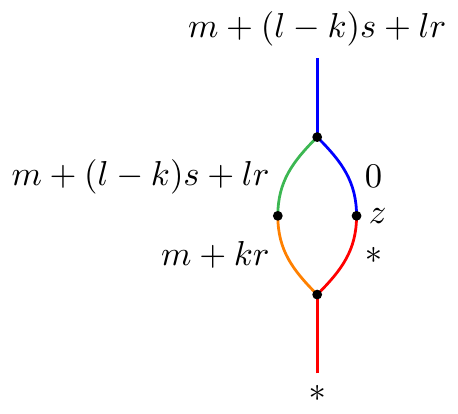}{
      \begin{tikzpicture}[scale = 0.4]
        \diamondvectorexampletwo{$*$}{$m+(l-k)s+lr$}{}{$z$}{$m+kr$}{$*$}{$0$}{$m+(l-k)s+lr$};
       \end{tikzpicture}
      }
    \end{array}.
\end{align}
Therefore, if we relabel $t^\prime=m+(l-k)s+lr$ and $m^\prime=m+kr$ we have the following basis when we quotient away the bubble action
\begin{align}
  \begin{array}{c}
    \includeTikz{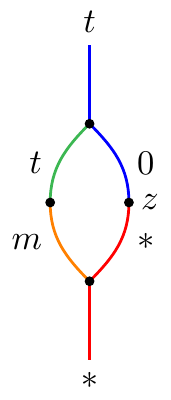}{
      \begin{tikzpicture}[scale = 0.4]
        \diamondvectorexampletwo{$*$}{$t$}{}{$z$}{$m$}{$*$}{$0$}{$t$};
       \end{tikzpicture}
      }
    \end{array}.
\end{align}
If we want to act by the idempotent $\defect{F_0}{R}{\zeta}{}{}$, we must have $t=0$.  The result of applying the projection is nonzero if and only if $\zeta = z$. Since $m$ is arbitrary, we have
$\defect{X_k}{X_l}{}{}{} \otimes \defect{F_0}{R}{z}{}{} \cong p \cdot \defect{F_0}{R}{z}{}{}$.
\end{example}

\begin{example}
  Consider the horizontal fusion $\defect{F_0}{F_r}{}{}{} \otimes \defect{T}{F_t}{}{}{}$. First we construct the compound representation of ${\bf Ann}_{L,X_{r^{-1}t}}$
\begin{align}
  \begin{array}{c}
    \includeTikz{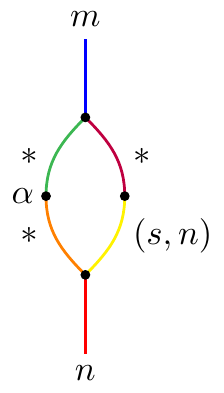}{
      \begin{tikzpicture}[scale = 0.4]
        \diamondvectorexamplethree{$n$}{$m$}{$\alpha$}{}{$*$}{$(s,n)$}{$*$}{$*$};
       \end{tikzpicture}
      }
    \end{array}.
\end{align}
Acting by an $s$ bubble sends this vector to 
\begin{align} \omega^{tsn+sr(\alpha-m)}
  \begin{array}{c}
    \includeTikz{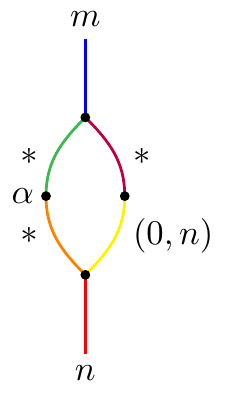}{
      \begin{tikzpicture}[scale = 0.4]
        \diamondvectorexamplethree{$n$}{$m$}{$\alpha$}{}{$*$}{$(0,n)$}{$*$}{$*$};
       \end{tikzpicture}
      }
    \end{array}
\end{align}
which forms a basis for the quotient. To apply $\defect{L}{X_{r^{-1}t}}{}{}{}$, we must have $m=n=0$. Since $\alpha \in \ZZ{p}$ (the representation corresponding to $\defect{F_0}{F_r}{}{}{}$) is $p$-dimensional, we have $\defect{F_0}{F_r}{}{}{} \otimes \defect{T}{F_t}{}{}{} \cong p \cdot \defect{L}{X_{r^{-1}t}}{}{}{}$ 
\end{example}

% This line sets the project root file.
% !TEX root = ../defects_domain_wall_structures.tex
\section{Bimodule Associators} \label{sec:bimodule_associators}

Now that we have seen that we must quotient away the bubble actions corresponding to internal cavities, we have seen everything needed to compute the compound defects corresponding to arbitrarily complex domain wall structures. Another interesting example is the following compound defect.
\begin{align}
[M,N,P]:=
\begin{array}{c}
\includeTikz{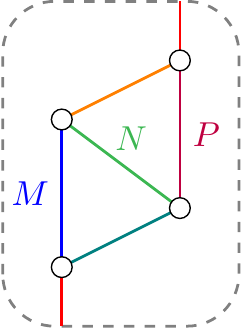}
{
	\begin{tikzpicture}[scale=0.3]
	\genericbimoduleassociator{$M$}{$N$}{$P$}
	\end{tikzpicture}
}
\end{array}\label{eqn:bimod_associator_def}.
\end{align}
We shall call this domain wall structure the bimodule associator for the triple $M,N,P$. If this defect projects onto nontrivial point defects, it indicates an obstruction to defining an extension (as described in Def.~\ref{def:BPR}). For $M,N,P$ invertible, this is closely related to the $O_3$ obstruction of \onlinecite{MR2677836} being nontrivial (see Appendix~\ref{sec:extension_theory}). From a physics viewpoint, this obstruction means we cannot \emph{gauge} the defects\cite{Barkeshli2014,1804.01657,MR3555361}.

We now provide an example calculation of a bimodule associator. The full set of associators (all trivial) can be found in Table~\ref{tab:bimod_associator_table}.

\begin{example}
	Let us compute the bimodule associator $[F_q,L,X_l]$. First, we can construct the following representation out of our trivalent vertices
	\begin{align}
	\begin{array}{c}
	\includeTikz{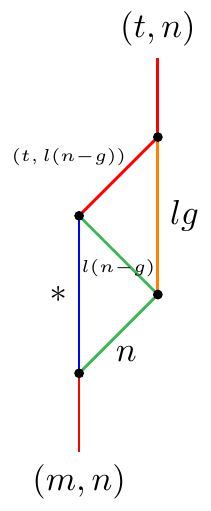}{
		\begin{tikzpicture}[scale=0.4]
		\bimoduleassociatorexample{$(m,n)$}{$*$}{$n$}{\tiny$l(n{-}g)$}{$lg$}{\tiny
		 $(t,l(n{-}g))$}{$(t,n)$}
		\end{tikzpicture}
	}
	\end{array}.
	\end{align}
	We need to quotient out the bubble actions from both of the cavities. Acting by an $lg$ bubble in the top cavity sends the vector to
	\begin{align}
	\begin{array}{c}
	\includeTikz{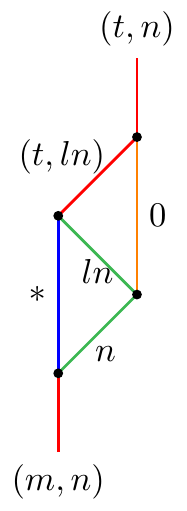}{
		\begin{tikzpicture}[scale=0.4]
		\bimoduleassociatorexample{$(m,n)$}{$*$}{$n$}{$ln$}{$0$}{$(t,ln)$}{$(t,n)$}
		\end{tikzpicture}
	}
	\end{array}.
	\end{align}
	Acting by a $h$ bubble in the bottom cavity multiplies this vector by $\omega^{hq(m-t)}$. So the vector is projected to zero in the quotient unless $m=t$. Therefore we have the following basis for the quotient
	\begin{align}
	\begin{array}{c}
	\includeTikz{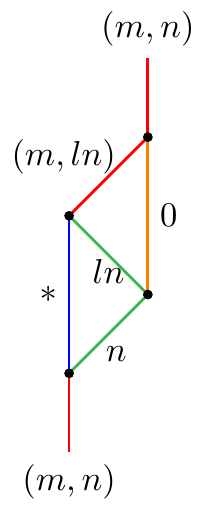}{
		\begin{tikzpicture}[scale=0.4]
		\bimoduleassociatorexample{$(m,n)$}{$*$}{$n$}{$ln$}{$0$}{$(m,ln)$}{$(m,n)$}
		\end{tikzpicture}
	}
	\end{array}.
	\end{align}
	To apply $\defect{T}{T}{\alpha}{\beta}{}$, we must have $m=n=\alpha=\beta=0$. Therefore the compound defect is $\defect{T}{T}{0}{0}{}$.
\end{example}

\begin{example}
We can also compute bimodule associators using the physical interpretations of the bimodules from \onlinecite{1806.01279} (Table~\ref{tab:bimodinterp}). The parameters $\mu$ and $\nu$ in our 3-string annular category representations physically correspond to the presence of a non-condensable anyon at the corner. Recall that the rough boundary condenses the $e$ anyons and the smooth boundary condenses the $m$ anyons.
\begin{align}
\left[R,F_y,R\right]\to
	\begin{array}{c}
		\includeTikz{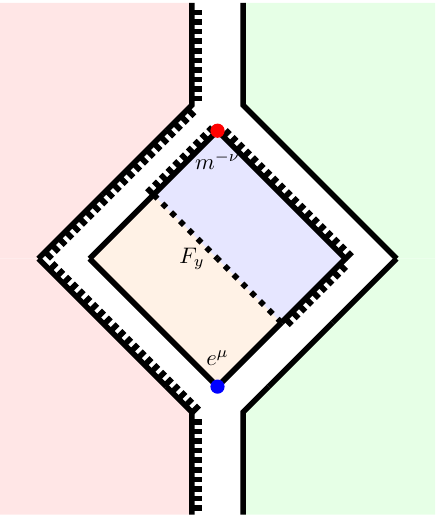}
		{
			\begin{tikzpicture}[scale=.65,every node/.style={scale=.65}]
			\def\dx{.4};
			%			\draw (0,0)--(1,1);\draw (1,1)--(2,0);\draw(1,1)--(2,2);\draw(2,2)--(4,0);\draw(2,2)--(2,3);
			\fill[red!10](-1,0)--(-\dx,0)--(2-\dx,2)--(2,2+\dx)--(2,4)--(-1,4)--cycle;
			\fill[red!10](-1,0)--(-\dx,0)--(2-\dx,-2)--(2,-2-\dx)--(2,-4)--(-1,-4)--cycle;
			\fill[orange!10] (\dx,0)--(1+\dx,1)--(3+\dx,-1)--(2+\dx,-2)--cycle;
			\begin{scope}[shift={(4+2*\dx,0)},scale=-1]
			\fill[green!10](-1,0)--(-\dx,0)--(2-\dx,2)--(2,2+\dx)--(2,4)--(-1,4)--cycle;
			\fill[green!10](-1,0)--(-\dx,0)--(2-\dx,-2)--(2,-2-\dx)--(2,-4)--(-1,-4)--cycle;
			\fill[blue!10] (\dx,0)--(1+\dx,1)--(3+\dx,-1)--(2+\dx,-2)--cycle;
			\end{scope}
			\begin{scope}
			\draw[ultra thick,shift={(-\dx,0)}] (0,0)--(1,1);
			\draw[ultra thick,shift={(\dx,0)}] (0,0)--(1,1);
			\draw[ultra thick,dotted,shift={(\dx,0)}] (1,1)--(2,0);
			\draw[ultra thick,shift={(-\dx,0)}](1,1)--(2,2);
			\draw[ultra thick,shift={(\dx,0)}](1,1)--(2,2);
			\draw[ultra thick,shift={(\dx,0)}](2,2)--(4,0);
			\draw[ultra thick,shift={(3*\dx,0)}](2,2)--(4,0);
			\draw[ultra thick,shift={(-\dx,0)}](2,2)--(2+\dx,2+\dx)--(2+\dx,4);
			\draw[ultra thick,shift={(3*\dx,0)}](2,2)--(2-\dx,2+\dx)--(2-\dx,4);
			\end{scope}
			\begin{scope}[shift={(4+2*\dx,0)},scale=-1]
			\draw[ultra thick,shift={(-\dx,0)}] (0,0)--(1,1);
			\draw[ultra thick,shift={(\dx,0)}] (0,0)--(1,1);
			\draw[ultra thick,dotted,shift={(\dx,0)}] (1,1)--(2,0);
			\draw[ultra thick,shift={(-\dx,0)}](1,1)--(2,2);
			\draw[ultra thick,shift={(\dx,0)}](1,1)--(2,2);
			\draw[ultra thick,shift={(\dx,0)}](2,2)--(4,0);
			\draw[ultra thick,shift={(3*\dx,0)}](2,2)--(4,0);
			\draw[ultra thick,shift={(-\dx,0)}](2,2)--(2+\dx,2+\dx)--(2+\dx,4);
			\draw[ultra thick,shift={(3*\dx,0)}](2,2)--(2-\dx,2+\dx)--(2-\dx,4);
			\end{scope}
			\foreach \x in {0,...,9}
			{
				\draw[ultra thick,shift={(2,-3)}] (0,-.9+.15*\x)--(.15,-.9+.15*\x);
				\draw[ultra thick,shift={(2,3+\dx)}] (0,-.9+.15*\x)--(.15,-.9+.15*\x);
			};
			\foreach \x in {1,...,22}
			{
				\draw[ultra thick,shift={(-\dx,0)},rotate=-45] (0,.15*\x)--(.15,.15*\x);
				\draw[,yscale=-1,ultra thick,shift={(-\dx,0)},rotate=-45] (0,.15*\x+.1)--(.15,.15*\x+.1);
			};
			\foreach \x in {0,...,9}
			{
				\draw[ultra thick,shift={(1+\dx-.1,1.1)},rotate=-45] (0,.15*\x)--(.15,.15*\x);				
			};
			\foreach \x in {0,...,-18}
			{
				\draw[scale=-1,ultra thick,shift={(-4-\dx-.1,-.1)},rotate=45] (0,.15*\x)--(.15,.15*\x);
			};
			\foreach \x in {-10,...,-18}
			{
				\draw[ultra thick,shift={(3+\dx-.1-2*\dx,1.1-2-3*\dx)},scale=-1,,rotate=-45] (0,.15*\x)--(.15,.15*\x);
			};
			\filldraw[red] (2+\dx,2) circle (.1);
			\node[below,black] at (2+\dx,2-\dx/2) {$m^{-\nu}$};
			\filldraw[blue] (2+\dx,-2) circle (.1);
			\node[above,black] at (2+\dx,-2+\dx/2) {$e^{\mu}$};
			\node at (2,0) {$F_y$};
			\end{tikzpicture}
		}
	\end{array}
\to&
\begin{array}{c}
	\includeTikz{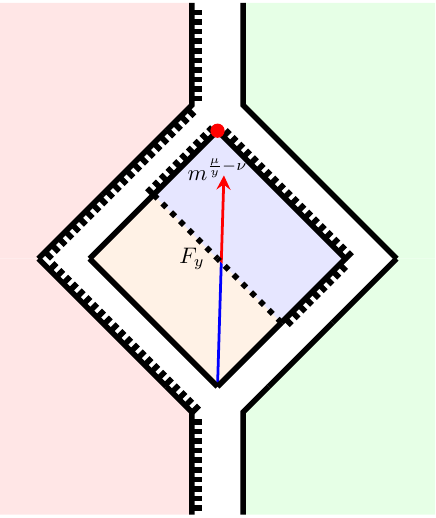}
	{
		\begin{tikzpicture}[scale=.65,every node/.style={scale=.65}]
		\def\dx{.4};
		%			\draw (0,0)--(1,1);\draw (1,1)--(2,0);\draw(1,1)--(2,2);\draw(2,2)--(4,0);\draw(2,2)--(2,3);
		\fill[red!10](-1,0)--(-\dx,0)--(2-\dx,2)--(2,2+\dx)--(2,4)--(-1,4)--cycle;
		\fill[red!10](-1,0)--(-\dx,0)--(2-\dx,-2)--(2,-2-\dx)--(2,-4)--(-1,-4)--cycle;
		\fill[orange!10] (\dx,0)--(1+\dx,1)--(3+\dx,-1)--(2+\dx,-2)--cycle;
		\begin{scope}[shift={(4+2*\dx,0)},scale=-1]
		\fill[green!10](-1,0)--(-\dx,0)--(2-\dx,2)--(2,2+\dx)--(2,4)--(-1,4)--cycle;
		\fill[green!10](-1,0)--(-\dx,0)--(2-\dx,-2)--(2,-2-\dx)--(2,-4)--(-1,-4)--cycle;
		\fill[blue!10] (\dx,0)--(1+\dx,1)--(3+\dx,-1)--(2+\dx,-2)--cycle;
		\end{scope}
				\begin{scope}
		\clip (\dx,0)--(1+\dx,1)--(3+\dx,-1)--(2+\dx,-2)--cycle;
		\draw[thick,blue,-stealth] (2+\dx,-2)->(2.5,1.3);
		\end{scope}
		\begin{scope}
		\clip (1+\dx,1)--(2+\dx,2)--(4+\dx,0)--(3+\dx,-1)--cycle;
		\draw[thick,red,-stealth] (2+\dx,-2)->(2.5,1.3);
		\end{scope}
		\begin{scope}
		\draw[ultra thick,shift={(-\dx,0)}] (0,0)--(1,1);
		\draw[ultra thick,shift={(\dx,0)}] (0,0)--(1,1);
		\draw[ultra thick,dotted,shift={(\dx,0)}] (1,1)--(2,0);
		\draw[ultra thick,shift={(-\dx,0)}](1,1)--(2,2);
		\draw[ultra thick,shift={(\dx,0)}](1,1)--(2,2);
		\draw[ultra thick,shift={(\dx,0)}](2,2)--(4,0);
		\draw[ultra thick,shift={(3*\dx,0)}](2,2)--(4,0);
		\draw[ultra thick,shift={(-\dx,0)}](2,2)--(2+\dx,2+\dx)--(2+\dx,4);
		\draw[ultra thick,shift={(3*\dx,0)}](2,2)--(2-\dx,2+\dx)--(2-\dx,4);
		\end{scope}
		\begin{scope}[shift={(4+2*\dx,0)},scale=-1]
		\draw[ultra thick,shift={(-\dx,0)}] (0,0)--(1,1);
		\draw[ultra thick,shift={(\dx,0)}] (0,0)--(1,1);
		\draw[ultra thick,dotted,shift={(\dx,0)}] (1,1)--(2,0);
		\draw[ultra thick,shift={(-\dx,0)}](1,1)--(2,2);
		\draw[ultra thick,shift={(\dx,0)}](1,1)--(2,2);
		\draw[ultra thick,shift={(\dx,0)}](2,2)--(4,0);
		\draw[ultra thick,shift={(3*\dx,0)}](2,2)--(4,0);
		\draw[ultra thick,shift={(-\dx,0)}](2,2)--(2+\dx,2+\dx)--(2+\dx,4);
		\draw[ultra thick,shift={(3*\dx,0)}](2,2)--(2-\dx,2+\dx)--(2-\dx,4);
		\end{scope}
		\foreach \x in {0,...,9}
		{
			\draw[ultra thick,shift={(2,-3)}] (0,-.9+.15*\x)--(.15,-.9+.15*\x);
			\draw[ultra thick,shift={(2,3+\dx)}] (0,-.9+.15*\x)--(.15,-.9+.15*\x);
		};
		\foreach \x in {1,...,22}
		{
			\draw[ultra thick,shift={(-\dx,0)},rotate=-45] (0,.15*\x)--(.15,.15*\x);
			\draw[,yscale=-1,ultra thick,shift={(-\dx,0)},rotate=-45] (0,.15*\x+.1)--(.15,.15*\x+.1);
		};
		\foreach \x in {0,...,9}
		{
			\draw[ultra thick,shift={(1+\dx-.1,1.1)},rotate=-45] (0,.15*\x)--(.15,.15*\x);				
		};
		\foreach \x in {0,...,-18}
		{
			\draw[scale=-1,ultra thick,shift={(-4-\dx-.1,-.1)},rotate=45] (0,.15*\x)--(.15,.15*\x);
		};
		\foreach \x in {-10,...,-18}
		{
			\draw[ultra thick,shift={(3+\dx-.1-2*\dx,1.1-2-3*\dx)},scale=-1,,rotate=-45] (0,.15*\x)--(.15,.15*\x);
		};
		\filldraw[red] (2+\dx,2) circle (.1);
		\node[below,black] at (2+\dx,2-\dx/1.3) {$m^{\frac{\mu}{y}{-}\nu}$};
		\node at (2,0) {$F_y$};
		\end{tikzpicture}
	}
\end{array}
\to
\begin{array}{c}
	\includeTikz{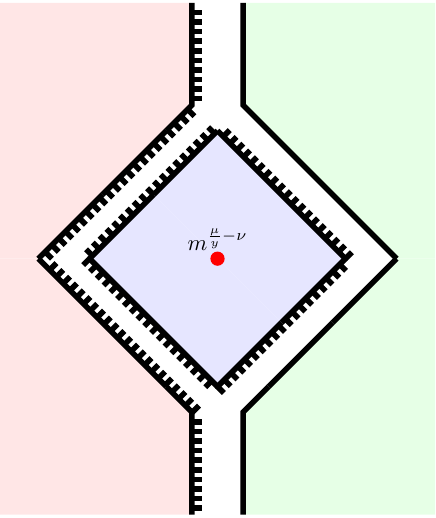}
	{
		\begin{tikzpicture}[scale=.65,every node/.style={scale=.65}]
		\def\dx{.4};
		%			\draw (0,0)--(1,1);\draw (1,1)--(2,0);\draw(1,1)--(2,2);\draw(2,2)--(4,0);\draw(2,2)--(2,3);
		\fill[red!10](-1,0)--(-\dx,0)--(2-\dx,2)--(2,2+\dx)--(2,4)--(-1,4)--cycle;
		\fill[red!10](-1,0)--(-\dx,0)--(2-\dx,-2)--(2,-2-\dx)--(2,-4)--(-1,-4)--cycle;
		\fill[blue!10] (\dx,0)--(1+\dx,1)--(3+\dx,-1)--(2+\dx,-2)--cycle;
		\begin{scope}[shift={(4+2*\dx,0)},scale=-1]
		\fill[green!10](-1,0)--(-\dx,0)--(2-\dx,2)--(2,2+\dx)--(2,4)--(-1,4)--cycle;
		\fill[green!10](-1,0)--(-\dx,0)--(2-\dx,-2)--(2,-2-\dx)--(2,-4)--(-1,-4)--cycle;
		\fill[blue!10] (\dx,0)--(1+\dx,1)--(3+\dx,-1)--(2+\dx,-2)--cycle;
		\end{scope}
		\begin{scope}
		\draw[ultra thick,shift={(-\dx,0)}] (0,0)--(1,1);
		\draw[ultra thick,shift={(\dx,0)}] (0,0)--(1,1);
%		\draw[ultra thick,dotted,shift={(\dx,0)}] (1,1)--(2,0);
		\draw[ultra thick,shift={(-\dx,0)}](1,1)--(2,2);
		\draw[ultra thick,shift={(\dx,0)}](1,1)--(2,2);
		\draw[ultra thick,shift={(\dx,0)}](2,2)--(4,0);
		\draw[ultra thick,shift={(3*\dx,0)}](2,2)--(4,0);
		\draw[ultra thick,shift={(-\dx,0)}](2,2)--(2+\dx,2+\dx)--(2+\dx,4);
		\draw[ultra thick,shift={(3*\dx,0)}](2,2)--(2-\dx,2+\dx)--(2-\dx,4);
		\end{scope}
		\begin{scope}[shift={(4+2*\dx,0)},scale=-1]
		\draw[ultra thick,shift={(-\dx,0)}] (0,0)--(1,1);
		\draw[ultra thick,shift={(\dx,0)}] (0,0)--(1,1);
%		\draw[ultra thick,dotted,shift={(\dx,0)}] (1,1)--(2,0);
		\draw[ultra thick,shift={(-\dx,0)}](1,1)--(2,2);
		\draw[ultra thick,shift={(\dx,0)}](1,1)--(2,2);
		\draw[ultra thick,shift={(\dx,0)}](2,2)--(4,0);
		\draw[ultra thick,shift={(3*\dx,0)}](2,2)--(4,0);
		\draw[ultra thick,shift={(-\dx,0)}](2,2)--(2+\dx,2+\dx)--(2+\dx,4);
		\draw[ultra thick,shift={(3*\dx,0)}](2,2)--(2-\dx,2+\dx)--(2-\dx,4);
		\end{scope}
		\foreach \x in {0,...,9}
		{
			\draw[ultra thick,shift={(2,-3)}] (0,-.9+.15*\x)--(.15,-.9+.15*\x);
			\draw[ultra thick,shift={(2,3+\dx)}] (0,-.9+.15*\x)--(.15,-.9+.15*\x);
		};
		\foreach \x in {1,...,22}
		{
			\draw[ultra thick,shift={(-\dx,0)},rotate=-45] (0,.15*\x)--(.15,.15*\x);
			\draw[,yscale=-1,ultra thick,shift={(-\dx,0)},rotate=-45] (0,.15*\x+.1)--(.15,.15*\x+.1);
		};
		\foreach \x in {-9,...,9}
		{
			\draw[ultra thick,shift={(1+\dx-.1,1.1)},rotate=-45] (0,.15*\x)--(.15,.15*\x);				
		};
		\foreach \x in {0,...,-18}
		{
			\draw[scale=-1,ultra thick,shift={(-4-\dx-.1,-.1)},rotate=45] (0,.15*\x)--(.15,.15*\x);
			\draw[scale=1,ultra thick,shift={(\dx-.1,-.1)},rotate=45] (0,.15*\x)--(.15,.15*\x);
		};
		\foreach \x in {-0,...,-18}
		{
			\draw[ultra thick,shift={(3+\dx-.1-2*\dx,1.1-2-3*\dx)},scale=-1,,rotate=-45] (0,.15*\x)--(.15,.15*\x);
		};
		\filldraw[red] (2+\dx,0) circle (.1);
		\node[above,black] at (2+\dx,0) {$m^{\frac{\mu}{y}{-}\nu}$};
		\end{tikzpicture}
	}
\end{array}
\\
\begin{array}{c}
	\includeTikz{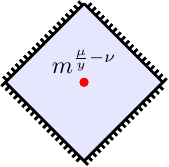}
	{
		\begin{tikzpicture}[scale=.4,every node/.style={scale=.7}]
		\def\dx{.4};;
		\fill[blue!10] (\dx,0)--(1+\dx,1)--(3+\dx,-1)--(2+\dx,-2)--cycle;
		\begin{scope}[shift={(4+2*\dx,0)},scale=-1]
		\fill[blue!10] (\dx,0)--(1+\dx,1)--(3+\dx,-1)--(2+\dx,-2)--cycle;
		\end{scope}
		\begin{scope}
		\draw[thick,shift={(\dx,0)}] (0,0)--(1,1);
		\draw[thick,shift={(\dx,0)}](1,1)--(2,2);
		\draw[thick,shift={(\dx,0)}](2,2)--(4,0);
		\end{scope}
		\begin{scope}[shift={(4+2*\dx,0)},scale=-1]
		\draw[thick,shift={(\dx,0)}] (0,0)--(1,1);
		\draw[thick,shift={(\dx,0)}](1,1)--(2,2);
		\draw[thick,shift={(\dx,0)}](2,2)--(4,0);
		\end{scope}
		\foreach \x in {-9,...,9}
		{
			\draw[thick,shift={(1+\dx-.1,1.1)},rotate=-45] (0,.15*\x)--(.15,.15*\x);				
		};
		\foreach \x in {0,...,-18}
		{
			\draw[scale=-1,thick,shift={(-4-\dx-.1,-.1)},rotate=45] (0,.15*\x)--(.15,.15*\x);
			\draw[scale=1, thick,shift={(\dx-.1,-.1)},rotate=45] (0,.15*\x)--(.15,.15*\x);
		};
		\foreach \x in {-0,...,-18}
		{
			\draw[thick,shift={(3+\dx-.1-2*\dx,1.1-2-3*\dx)},scale=-1,,rotate=-45] (0,.15*\x)--(.15,.15*\x);
		};
		\filldraw[red] (2+\dx,0) circle (.1);
		\node[above,black] at (2+\dx,0) {$m^{\frac{\mu}{y}{-}\nu}$};
		\end{tikzpicture}
	}
\end{array}
&=\delta_{\mu}^{y\nu},
\end{align}
since the internal disc must contain $0$ anyons that cannot be fused into the boundary. Therefore this associator is $\delta^{\nu y}_{\mu}\defect{R}{R}{0}{0}{\mu,\nu}{}$.
\end{example}

% This line sets the project root file.
% !TEX root = ../defects_domain_wall_structures.tex
% !TEX spellcheck = en_US

\section{A Defect Levin-Wen Lattice Model}\label{sec:defectLW}

A particularly important application of the domain wall structure algorithm is finding the ground states of Levin-Wen Hamiltonians\cite{Levin2005}, which are compound defects. This requires a slight modification of the usual definition of these models.
In this section, we show how a Levin-Wen type model\cite{Levin2005}, complete with defects, can be constructed using annular category representations.

For simplicity, we restrict to the following lattice
\begin{align}
  \begin{array}{c}
    \includeTikz{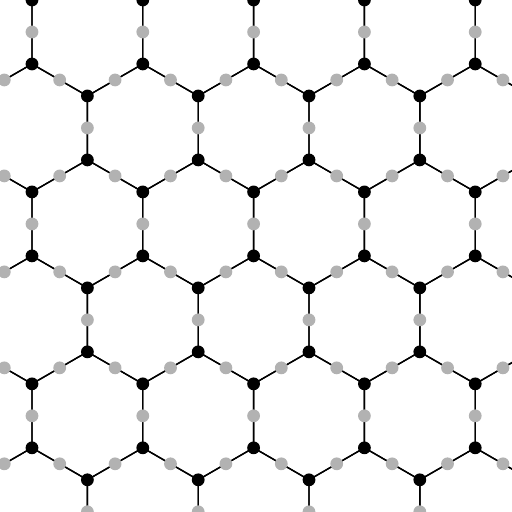}{
\begin{tikzpicture}[scale=.65,every node/.style={scale=.85}]
\clip (-.5,-1) rectangle (7.5,7);
\pgfmathsetmacro\dx{sin(60)};
\pgfmathsetmacro\dy{cos(60)};
;
\foreach \x in {-5,...,10}{
	\foreach \y in {-5,...,10}{
		\pgfmathsetmacro\my{mod(\y,2)};
		\begin{scope}[shift={(2*\dx*\x+\dx*\my,\dy*\y+\y)}]
		\draw (90:0)--(90:1) (90+120:0)--(90+120:1) (90-120:0)--(90-120:1) ;
		\draw[shift={(\dx,-\dy)},rotate=180] (90:0)--(90:1) (90+120:0)--(90+120:1) (90-120:0)--(90-120:1) ;
		\end{scope}
	}
}
\foreach \x in {-5,...,10}{
	\foreach \y in {-5,...,10}{
		\pgfmathsetmacro\my{mod(\y,2)};
		\begin{scope}[shift={(2*\dx*\x+\dx*\my,\dy*\y+\y)}]
		\fill (90:0) circle (.1);
		\fill (90:1) circle (.1);
		\fill[black!30] (90:.5) circle (.1);
		\fill[black!30] (90-120:.5) circle (.1);
		\fill[black!30] (90+120:.5) circle (.1);
		\end{scope}
	}
}
\end{tikzpicture}
      }
  \end{array}.
\end{align}
The model described here has degrees of freedom both on edges and vertices, indicated by gray and black dots respectively.

Each face of the lattice is labeled with a fusion category $\cat{C}_f$, edges are labeled with bimodules $M_e$ between the appropriate $\cat{C}_f$. Each vertex is labeled with an (irreducible) representation $V^z$ of the corresponding annular category. Recall that this representation consists of a set of vector spaces $V^z_{m_1,m_2,m_3}$, where the $m_i$ are the edge labels. A partial list of irreducible 3-bimodule annular category representations are provided in Tables~\ref{tab:trivertextable} and \ref{tab:dualtrivertextable}. All other irreducible representations can be obtained by composing with the 2-bimodule vertex representations from Table~\ref{tab:bivertextable} in the following way
\begin{align}
  \begin{array}{c}
    \includeTikz{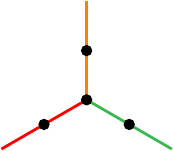}{
      \begin{tikzpicture}
        \draw[thick,red] (0,0) -- (-0.866025, -0.5);
        \draw[thick,nicegreen] (0,0) -- (0.866025, -0.5);
        \draw[thick,orange] (0,0) -- (0,1);
        \draw[fill=black] (0,0) circle (0.05);
        \draw[fill=black] (-0.433013, -0.25) circle (0.05);
        \draw[fill=black] (0.433013, -0.25) circle (0.05);
        \draw[fill=black] (0, 0.5) circle (0.05);
      \end{tikzpicture}
      }
  \end{array} \qquad
\begin{array}{c}
    \includeTikz{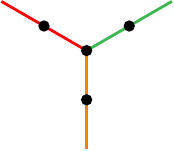}{
      \begin{tikzpicture}[yscale=-1]
        \draw[thick,red] (0,0) -- (-0.866025, -0.5);
        \draw[thick,nicegreen] (0,0) -- (0.866025, -0.5);
        \draw[thick,orange] (0,0) -- (0,1);
        \draw[fill=black] (0,0) circle (0.05);
        \draw[fill=black] (-0.433013, -0.25) circle (0.05);
        \draw[fill=black] (0.433013, -0.25) circle (0.05);
        \draw[fill=black] (0, 0.5) circle (0.05);
      \end{tikzpicture}
      }
  \end{array}.
\end{align}

A basis for the edge Hilbert space is given by the set of objects $m\in M_e$. A basis for the vertex degree of freedom is given by a basis for the assigned representation.

If $\psi$ is a basis state on the lattice, in a neighborhood of the vertex $z$ it looks like
\begin{align}
\begin{array}{c}
\includeTikz{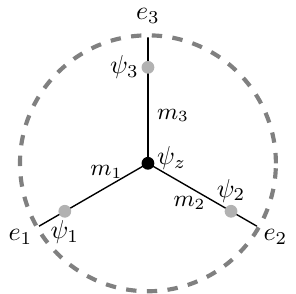}{
	\begin{tikzpicture}[scale=.65,every node/.style={scale=.75}]
	\draw (0,0) -- (-30:2) node [pos=1.15] {$e_2$};
	\draw (0,0) -- (210:2) node [pos=1.15] {$e_1$};
	\draw (0,0) -- (0,2) node [pos=1.15] {$e_3$};
	\draw[white,very thick] (0,0) circle (2);
	\draw[gray,very thick,dashed] (0,0) circle (2);	
	\fill [black] (0,0) circle (.1);
	\node[right] at (0,.075) { $\psi_z$};
	\fill [black!30] (210:1.5) circle (.1) node[below,text=black] { $\psi_1$};
	\fill [black!30] (-30:1.5) circle (.1) node[above,text=black] { $\psi_2$};
	\fill [black!30] (90:1.5) circle (.1) node[left,text=black] { $\psi_3$};
	\node[right] at (90:.75) {\small $m_3$};
	\node[below] at (-30:.75) {\small $m_2$};
	\node[above] at (210:.75) {\small $m_1$};
	\end{tikzpicture}
}
\end{array},
\end{align}
where $\psi_z \in V^z_{m_1,m_2,m_3}$, $\psi_i,m_i\in\ob{M_{e_i}}$. The Hamiltonian for this model consists of two parts. The vertex operator $H_{z,e_i}$ projects onto the states where $\psi_i = m_i$. We define
\begin{align}
  H_z = H_{z,e_1} + H_{z,e_2} + H_{z,e_3}.
\end{align}
 In a neighborhood of the face $f$, the state $\psi$ looks like
\begin{align}
  \begin{array}{c}
    \includeTikz{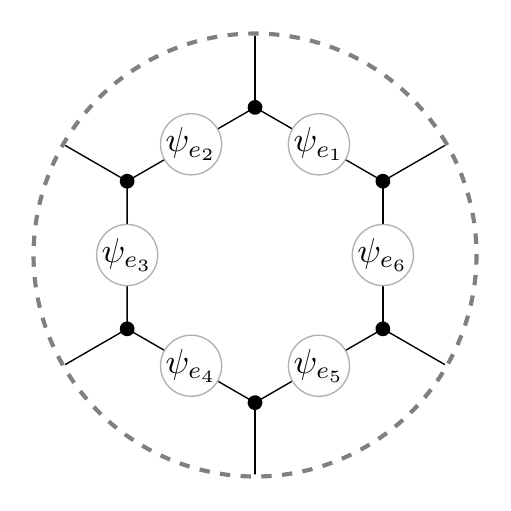}{
      \begin{tikzpicture}[scale=1.5,rotate=30]
        \faceOpenNh{$\psi_{e_1}$}{$\psi_{e_2}$}{$\psi_{e_3}$}{$\psi_{e_4}$}{$\psi_{e_5}$}{$\psi_{e_6}$}
      \end{tikzpicture}
      }
    \end{array}.
\end{align}
Fix $g \in \mathbb{Z}/p\mathbb{Z}$. The operator acts $H_{f,g}$ maps this to:
\begin{align}
  \begin{array}{c}
    \includeTikz{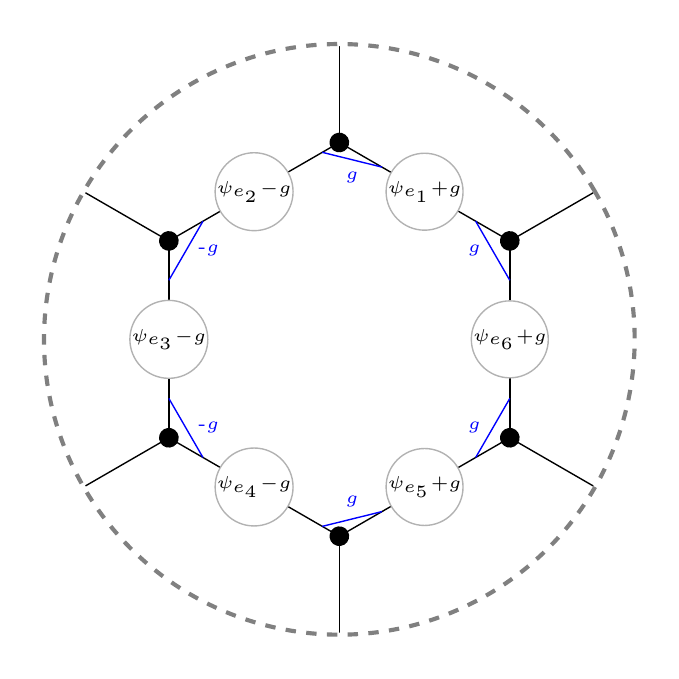}{
      \begin{tikzpicture}[scale=2,rotate=30]
        \faceOpenNh{\tiny $\psi_{e_1}{+}g$}{\tiny $\psi_{e_2}{-}g$}{\tiny $\psi_{e_3}{-}g$}{\tiny $\psi_{e_4}{-}g$}{\tiny $\psi_{e_5}{+g}$}{\tiny $\psi_{e_6}{+}g$}
        \draw[blue] (0.9, -0.173205) -- (0.9, 0.173205) node[midway,left] {\tiny $g$};
        \draw[blue] (-0.3, 0.866025) -- (-0.6, 0.69282) node[midway,right] {\tiny{-}$g$};
        \draw[blue] (-0.9, 0.173205) -- (-0.9, -0.173205) node[midway,right] {\tiny{-}$g$};
        \draw[blue] (0.3, -0.866025) -- (0.6, -0.69282) node[midway,left] {\tiny $g$};
        \draw[blue] (0.625, 0.649519) -- (0.4, 0.866025) node[midway,below] {\tiny $g$};
        \draw[blue] (-0.55, -0.779423) -- (-0.25, -0.866025) node[midway,above] {\tiny $g$};
      \end{tikzpicture}
      }
    \end{array}
\end{align}
Since the local degrees of freedom which live on the vertices are vectors in annular category representations, they can absorb the blue strings (as shown in Tables~\ref{tab:trivertextable} and \ref{tab:dualtrivertextable}). We define
\begin{align}
  H_f = \frac{1}{p} \sum_g H_{f,g}.
\end{align}
The Hamiltonian for the model is
\begin{align}
  \sum_z (1-H_z) + \sum_f (1-H_f).
\end{align}
where $z$ ranges over vertices and $f$ ranges over faces.

% This line sets the project root file.
% !TEX root = ../defects_domain_wall_structures.tex
% !TEX spellcheck = en_US

\subsection{Ground state configurations}

For clarity, we now show how the `string-net' ground state arises from this construction. Restrict to the bimodule $X_1$ for $\ZZ{p}$. From Tables~\ref{tab:trivertextable} and \ref{tab:dualtrivertextable}, all vector spaces arising in the annular category representation at a trivertex are one dimensional. We can therefore label the vertex degrees of freedom by the edge labels. A basis for the vertex Hilbert space is given by
\begin{align}
\left\{
\begin{array}{c}
\includeTikz{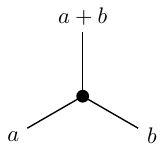}{
	\begin{tikzpicture}[scale=.65,every node/.style={scale=.65}]
		\draw(90:0)--(90:1) node[pos=1.25] {$a+b$};
		\draw(90+120:0)--(90+120:1) node[pos=1.25] {$a$};
		\draw(90-120:0)--(90-120:1) node[pos=1.25] {$b$};
		\fill (0,0) circle (.1);
	\end{tikzpicture}}
\end{array}\right\}_{a,b\in\ZZ{p}}.\label{eqn:bulkbasis}
\end{align}
The edge terms in the Hamiltonian force the edge labels to be consistent. The face terms fluctuate between configurations of edge labelings. For example, in the $\ZZ{2}$ case
\begin{align}
\begin{array}{c}
\includeTikz{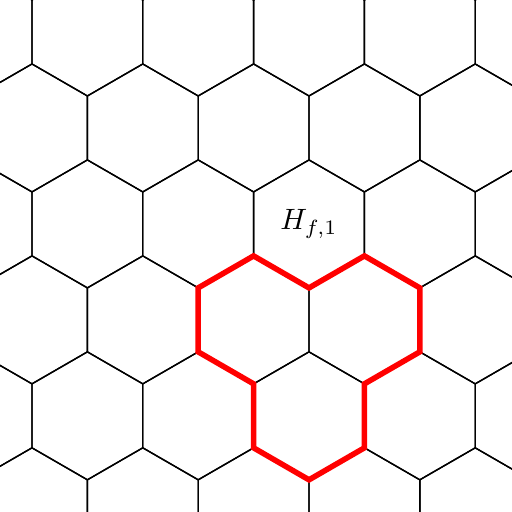}{
	\begin{tikzpicture}[scale=.65,every node/.style={scale=.85}]
	\clip (-.5,-1) rectangle (7.5,7);
	\pgfmathsetmacro\dx{sin(60)};
	\pgfmathsetmacro\dy{cos(60)};
	;
	\foreach \x in {-5,...,10}{
		\foreach \y in {-5,...,10}{
			\pgfmathsetmacro\my{mod(\y,2)};
	\begin{scope}[shift={(2*\dx*\x+\dx*\my,\dy*\y+\y)}]
		\draw (90:0)--(90:1) (90+120:0)--(90+120:1) (90-120:0)--(90-120:1) ;
		\draw[shift={(\dx,-\dy)},rotate=180] (90:0)--(90:1) (90+120:0)--(90+120:1) (90-120:0)--(90-120:1) ;
	\end{scope}
	}
	}
\def\x{3};\def\y{0};\pgfmathsetmacro\my{mod(\y,2)};
	\draw[ultra thick,red] (5*\dx,-\dy)--(6*\dx,0)--(6*\dx,2*\dy)--(7*\dx,3*\dy)--(7*\dx,5*\dy)--(6*\dx,6*\dy)--(6*\dx,6*\dy)--(5*\dx,5*\dy)--(4*\dx,6*\dy)--(3*\dx,5*\dy)--(3*\dx,3*\dy)--(4*\dx,2*\dy)--(4*\dx,0)--cycle;
	\node[] at (5*\dx,7*\dy) {$H_{f,1}$};
	\end{tikzpicture}
}
\end{array}
=
\begin{array}{c}
\includeTikz{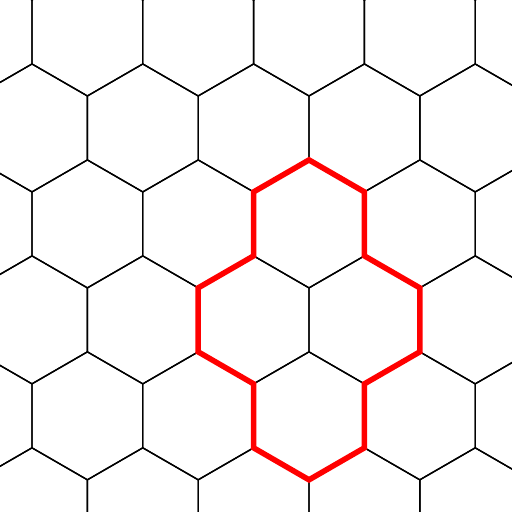}{
	\begin{tikzpicture}[scale=.65,every node/.style={scale=.85}]
	\clip (-.5,-1) rectangle (7.5,7);
	\pgfmathsetmacro\dx{sin(60)};
	\pgfmathsetmacro\dy{cos(60)};
	;
	\foreach \x in {-5,...,10}{
		\foreach \y in {-5,...,10}{
			\pgfmathsetmacro\my{mod(\y,2)};
			\begin{scope}[shift={(2*\dx*\x+\dx*\my,\dy*\y+\y)}]
			\draw (90:0)--(90:1) (90+120:0)--(90+120:1) (90-120:0)--(90-120:1) ;
			\draw[shift={(\dx,-\dy)},rotate=180] (90:0)--(90:1) (90+120:0)--(90+120:1) (90-120:0)--(90-120:1) ;
			\end{scope}
		}
	}
	\def\x{3};\def\y{0};\pgfmathsetmacro\my{mod(\y,2)};
	\draw[ultra thick,red] (5*\dx,-\dy)--(6*\dx,0)--(6*\dx,2*\dy)--(7*\dx,3*\dy)--(7*\dx,5*\dy)--(6*\dx,6*\dy)--(6*\dx,8*\dy)--(5*\dx,9*\dy)--(4*\dx,8*\dy)--(4*\dx,8*\dy)--(4*\dx,6*\dy)--(3*\dx,5*\dy)--(3*\dx,3*\dy)--(4*\dx,2*\dy)--(4*\dx,0)--cycle;
	\end{tikzpicture}
}
\end{array}
.
\end{align}

Excitations in this model are created by mismatches along the edges, for example the configuration
\begin{align}
\begin{array}{c}
\includeTikz{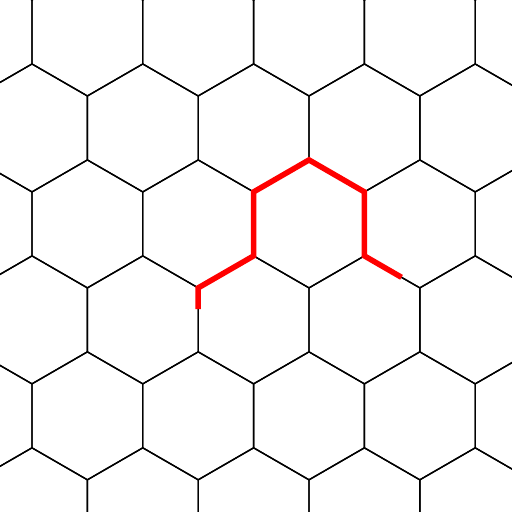}{
	\begin{tikzpicture}[scale=.65,every node/.style={scale=.85}]
	\clip (-.5,-1) rectangle (7.5,7);
	\pgfmathsetmacro\dx{sin(60)};
	\pgfmathsetmacro\dy{cos(60)};
	;
	\foreach \x in {-5,...,10}{
		\foreach \y in {-5,...,10}{
			\pgfmathsetmacro\my{mod(\y,2)};
			\begin{scope}[shift={(2*\dx*\x+\dx*\my,\dy*\y+\y)}]
			\draw (90:0)--(90:1) (90+120:0)--(90+120:1) (90-120:0)--(90-120:1) ;
			\draw[shift={(\dx,-\dy)},rotate=180] (90:0)--(90:1) (90+120:0)--(90+120:1) (90-120:0)--(90-120:1) ;
			\end{scope}
		}
	}
	\def\x{3};\def\y{0};\pgfmathsetmacro\my{mod(\y,2)};
	\draw[ultra thick,red] (7*\dx-.33*\dx,5*\dy+.33*\dy)--(6*\dx,6*\dy)--(6*\dx,8*\dy)--(5*\dx,9*\dy)--(4*\dx,8*\dy)--(4*\dx,8*\dy)--(4*\dx,6*\dy)--(3*\dx,5*\dy)--(3*\dx,5*\dy-.66*\dy);
	\end{tikzpicture}
}
\end{array}
\end{align}
violates two vertex terms, corresponding to a pair of particles (anyons). By inserting phases along the string, plaquette terms can also be violated.

Ground states with defect lines and points can also be described in this model. Consider the following arrangement of defect lines
\begin{align}
\begin{array}{c}
\includeTikz{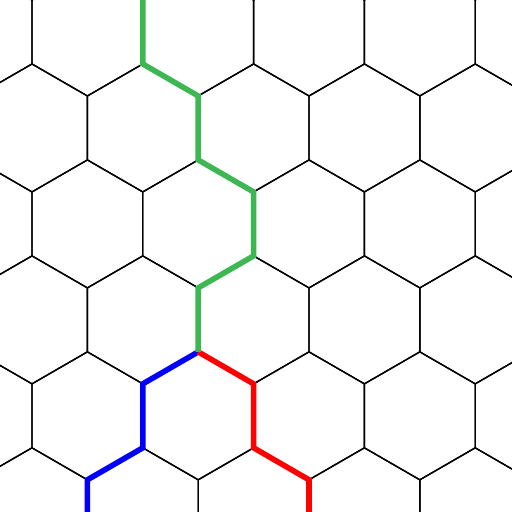}{
	\begin{tikzpicture}[scale=.65,every node/.style={scale=.85}]
	\clip (-.5,-1) rectangle (7.5,7);
	\pgfmathsetmacro\dx{sin(60)};
	\pgfmathsetmacro\dy{cos(60)};
	;
	\foreach \x in {-5,...,10}{
		\foreach \y in {-5,...,10}{
			\pgfmathsetmacro\my{mod(\y,2)};
			\begin{scope}[shift={(2*\dx*\x+\dx*\my,\dy*\y+\y)}]
			\draw (90:0)--(90:1) (90+120:0)--(90+120:1) (90-120:0)--(90-120:1) ;
			\draw[shift={(\dx,-\dy)},rotate=180] (90:0)--(90:1) (90+120:0)--(90+120:1) (90-120:0)--(90-120:1) ;
			\end{scope}
		}
	}
	\def\x{3};\def\y{0};\pgfmathsetmacro\my{mod(\y,2)};
%	\draw[ultra thick,blue] (7*\dx-.33*\dx,5*\dy+.33*\dy)--(6*\dx,6*\dy)--(6*\dx,8*\dy)--(5*\dx,9*\dy)--(4*\dx,8*\dy)--(4*\dx,8*\dy)--(4*\dx,6*\dy)--(3*\dx,5*\dy)--(3*\dx,5*\dy-.66*\dy);
\draw[ultra thick,blue] (1*\dx,-3*\dy)--(\dx,-1*\dy)--(2*\dx,0*\dy)--(2*\dx,2*\dy)--(3*\dx,3*\dy);
\draw[ultra thick,red] (5*\dx,-3*\dy)--(5*\dx,-1*\dy)--(4*\dx,0*\dy)--(4*\dx,2*\dy)--(3*\dx,3*\dy);
\draw[ultra thick,nicegreen] (3*\dx,3*\dy)--(3*\dx,5*\dy)--(4*\dx,6*\dy)--(4*\dx,8*\dy)--(3*\dx,9*\dy)--(3*\dx,11*\dy)--(2*\dx,12*\dy)--(2*\dx,14*\dy);
	\end{tikzpicture}
}
\end{array}
\leftrightarrow
\begin{array}{c}
\includeTikz{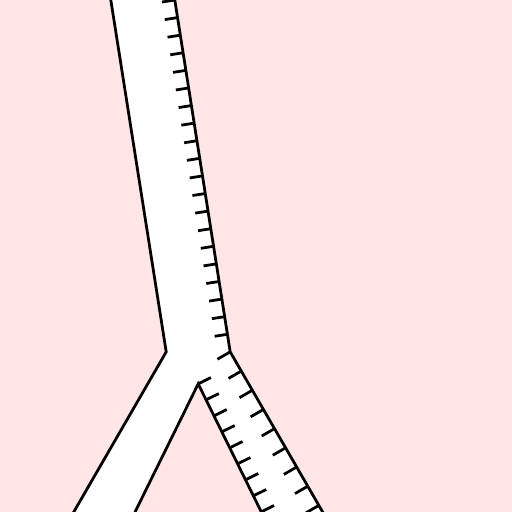}{
	\begin{tikzpicture}[scale=.65,every node/.style={scale=.85}]
		\pgfmathsetmacro\dx{sin(60)};
	\pgfmathsetmacro\dy{cos(60)};
	\clip (-.5,-1) rectangle (7.5,7);
%	\fill[red!10] (-.5,-1) rectangle (7.5,7);
	\fill[red!10] (1*\dx-.5,-3*\dy)--(3*\dx-.5,3*\dy)--(2*\dx-.5,14*\dy)--(-.5,7)--(-.5,-1)--cycle;
	\fill[red!10] (5*\dx+.5,-3*\dy)--(3*\dx+.5,3*\dy)--(2*\dx+.5,14*\dy)--(7.5,7)--(7.5,-1)--cycle;
	\fill[red!10] (1*\dx+.5,-3*\dy)--(3*\dx,2*\dy)--(5*\dx-.5,-3*\dy);
	\draw[thick, black] (1*\dx-.5,-3*\dy)--(3*\dx-.5,3*\dy)--(2*\dx-.5,14*\dy);
	\draw[thick, black] (5*\dx+.5,-3*\dy)--(3*\dx+.5,3*\dy)--(2*\dx+.5,14*\dy);
	\draw[thick,black] (1*\dx+.5,-3*\dy)--(3*\dx,2*\dy)--(5*\dx-.5,-3*\dy);
	\foreach \t in {0,...,10}{\draw[thick,shift={($(3*\dx,2*\dy)+(2*\dx*\t/10-\t/20,-5*\dy*\t/10)$)}] (0,0)--(.2,.2*4*\dx/\dy/10-.2/\dy/10);}
	\foreach \t in {0,...,10}{\draw[thick,shift={($(5*\dx+.5,-3*\dy)+(-2*\dx*\t/10,6*\dy*\t/10)$)}] (0,0)--(-.2,-.2*\dx/\dy/3);}
	\foreach \t in {1,...,20}{\draw[thick,shift={($(3*\dx+.5,3*\dy)+(-\dx*\t/20,11*\dy*\t/20)$)}] (0,0)--(-.2,-.2*\dx/\dy/11);}
	\end{tikzpicture}
}
\end{array},
\end{align}
where the blue, red, green lines are the bimodules $F_0$, $T$  and $L$ respectively. The remaining edges are labeled by the $X_1$ bimodule, corresponding to the bulk theory. Recall from \onlinecite{1806.01279} (reproduced in Table~\ref{tab:bimodinterp}) that these domain walls correspond to boundaries in the theory. At each vertex, we assign the representation from Tables~\ref{tab:trivertextable} and \ref{tab:dualtrivertextable}.

In addition to the vertex basis from Eqn.~\ref{eqn:bulkbasis}, we now require bases for the new vertices (see Tables~\ref{tab:trivertextable} and \ref{tab:dualtrivertextable}). As before, excitations can be created by violating either vertex or face terms, both in the bulk and in the vicinity of the boundary.
% This line sets the project root file.
% !TEX root = ../defects_domain_wall_structures.tex
% !TEX spellcheck = en_US

\section{Remarks}\label{sec:conclusions}
In this work, we have described a framework for computing the compound defect associated to a domain wall structure. The algorithm described is agnostic to the invertibility of the bimodules and point defects forming the structure. Using this algorithm, we have shown how the fusion (both vertical and horizontal) of defects are expressed as domain wall structures, and how the results of \onlinecite{1806.01279} can be replicated in this new, computer-friendly manner. Additionally, we have applied our algorithm to show that the \emph{domain wall associators} for all bimodules over $\vvec{\ZZ{p}}$ are trivial.

Although we have specialized to $\vvec{\ZZ{p}}$ for this work, the ideas described here are not restricted to this class of fusion categories. 
Due to the large number of fault tolerant gates that can be implemented, the category $\vvec{\ZZ{2}}\times\vvec{\ZZ{2}}$, called the color code in quantum computing, is of particular interest\cite{ColorCode,Yoshida2015a}. The large number of bimodules of this model (270) make a computer-implementable method, such as that outlined here, necessary to study the defects. We also expect these techniques to be useful for $\vvec{G}$ when $G$ is not abelian, and other non-abelian fusion categories.

We have shown how the domain wall associators can be computed in this framework. These associators are closely related to the $O_3$ obstruction of \onlinecite{MR2677836}. When this obstruction vanishes (as is the case for $\vvec{\ZZ{p}}$), a further obstruction, called $O_4$ in \onlinecite{MR2677836}, can arise. This obstruction is related to natural isomorphisms of defects. It would be extremely useful if the techniques developed in this work can be extended to include this data.

\acknowledgments
This work is supported by the Australian Research Council (ARC) via Discovery Project ``Subfactors and symmetries'' DP140100732 and Discovery Project ``Low dimensional categories'' DP16010347. This research was supported in part by Perimeter Institute for Theoretical Physics. Research at Perimeter Institute is supported by the Government of Canada through the Department of Innovation, Science and Economic Development Canada and by the Province of Ontario through the Ministry of Economic Development, Job Creation and Trade. We thank Corey Jones, Cain Edie-Michell and Scott Morrison for explaining to us how fusion category extension theory works.

\vspace{5mm}

%------------------------------------------------------------------------------------------------------------%
\bibliographystyle{apsrev_jacob}
\bibliography{refs}
%------------------------------------------------------------------------------------------------------------%

\onecolumngrid
\appendix

% This line sets the project root file.
% !TEX root = ../defects_domain_wall_structures.tex
% !TEX spellcheck = en_US
\clearpage
\section{A gentle review of bimodule categories over a fusion category}\label{sec:categoryreps}

To aid the unfamiliar reader, in this appendix we provide a sketch of the definitions of some of the mathematical structures used in this paper. Whilst these definitions are sufficient for studying the $\vvec{\ZZ{p}}$ case, we note that they must be refined if the algorithms described here are to be used for more general fusion categories. In particular, care must be taken when the fusion categories have nontrivial Frobenius-Schur indicators. In the general case, all the strings must be oriented and we need to be more careful with rigid structures. More complete definitions can be found in \onlinecite{MR3674995}.

The particle types in a long range entangled (2+1)D topological phase are labeled by objects in a \emph{modular tensor category} \cite{MR2640343}. The Drinfeld center $Z(\mathcal{C})$ of a fusion category $\mathcal{C}$ is a modular tensor category. The fusion category $\mathcal{C}$ encodes a local Hamiltonian which can be used to construct a gapped commuting projector Hamiltonian with anyonic excitations $Z(\mathcal{C})$ as described by Levin and Wen in \onlinecite{Levin2005}. It is important to not think of the objects in $\mathcal{C}$ as anyons. The objects in $\mathcal{C}$ are local degrees of freedom in the Levin-Wen lattice model. 

Throughout this paper, all categories are enriched over $\vvec{}$, the category of $\mathbb{C}$-vector spaces. This means that the Hom spaces are vector spaces and composition is bilinear. In particular, the endomorphisms of any object form an algebra.

\begin{definition}[Fusion category]
A {\em tensor category} $\mathcal{C}$ is a category $\mathcal{C}$ equipped with a functor $- \otimes - : \mathcal{C} \otimes \mathcal{C} \to \mathcal{C}$, a natural isomorphism $ (- \otimes -) \otimes - \cong - \otimes (- \otimes -)$ called the associator and a special object $1 \in \mathcal{C}$ which satisfy the pentagon equation and unit equations respectively. These can be found on Page~22 of \onlinecite{MR3242743}. If $\mathcal{C}$ is semi-simple, using the string diagram notation as explained in \onlinecite{MR3674995}, a vector in $\mathcal{C}(a \otimes b, c)$ can be represented by a trivalent vertex: 
\begin{align}
\begin{array}{c}
  \includeTikz{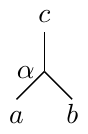}{
    \begin{tikzpicture}[scale=0.4, every node/.style={scale=0.8}]
      \trivalentvertex{$c$}{$a$}{$b$}{$\alpha$};
    \end{tikzpicture}
  }
\end{array}.
\end{align}

If we choose bases for all the vector spaces $\mathcal{C}(a \otimes b,c)$, then the associator can be represented as a tensor F, where
  \begin{align}
    \begin{array}{c}
	    \includeTikz{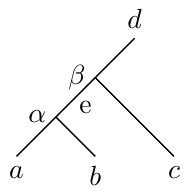}
	    {
	    	\begin{tikzpicture}[scale=0.4, every node/.style={scale=0.8}]
		    	\draw(0,0)--(3,3) node[below] at (0,0) {$a$} node[above] at (3,3) {$d$};
		    	\draw(2,0)--(1,1) node[below] at (2,0) {$b$};
		    	\draw(4,0)--(2,2) node[below] at (4,0) {$c$};
		    	\node[] at (1.75,1.25) {e};
		    	\node[left] at(1,1){$\alpha$};
		    	\node[left] at(2,2){$\beta$};
	    	\end{tikzpicture}
	    }
    \end{array}
    =\sum_{(f,\mu,\nu)}
    \biggl[F^{abc}_{d}\biggr]_{(e,\alpha,\beta)(f,\mu,\nu)}
    \begin{array}{c}
    \includeTikz{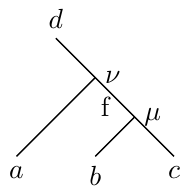}
    {
    	\begin{tikzpicture}[scale=0.4, every node/.style={scale=0.8}]
    	\draw(0,0)--(-3,3) node[below] at (0,0) {$c$} node[above] at (-3,3) {$d$};
    	\draw(-2,0)--(-1,1) node[below] at (-2,0) {$b$};
    	\draw(-4,0)--(-2,2) node[below] at (-4,0) {$a$};
    	\node at (-1.75,1.25) {f};
    	\node[right] at(-1,1){$\mu$};
    	\node[right] at(-2,2){$\nu$};
    	\end{tikzpicture}
    }
    \end{array}.
  \end{align}
  where $\alpha \in \mathcal{C}(a \otimes b,e), \beta \in \mathcal{C}(e \otimes c,d), \mu \in \mathcal{C}(b \otimes c,f)$ and $\nu \in \mathcal{C}(a \otimes f,d)$ are basis vectors.

A {\em fusion category} is a semi-simple rigid tensor category $\mathcal{C}$ with a finite number of simple objects and a simple unit. An object is {\em simple} if there are no non-trivial sub-objects. Semi-simple means that every object in $\mathcal{C}$ is a direct sum of simple objects. Rigid is a technical condition defined on Page~40 of \onlinecite{MR3242743} which implies that objects in $\mathcal{C}$ have duals and they behave like duals in the category of vector spaces. When doing fusion category computations, it is customary to fix a set of simple objects and express all other objects as direct sums of the chosen simple ones.

A subtle but important issue is the interaction between rigid structures and the graphical calculus used to describe morphisms in a fusion category. In order for the graphical calculus to fully capture the rigid structure, all strings need to be oriented. If the Frobenius-Schur indicators associated to the rigid structure are trivial, then the orientations can be ignored. This is case for the fusion category $\vvec{G}$, which is why we can ignore the string orientations.
\end{definition}

\begin{example}[$\vvec{G}$]
  Let $G$ be a finite group. The fusion category $\vvec{G}$ has simple objects the elements $g \in G$. The tensor product is $g \otimes h \cong gh$. Trivalent vertices
  \begin{align}
\begin{array}{c}
  \includeTikz{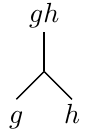}{
    \begin{tikzpicture}[scale=0.4, every node/.style={scale=0.8}]
      \trivalentvertex{$gh$}{$g$}{$h$}{};
    \end{tikzpicture}
  }
      \end{array}
  \end{align}
  can be chosen so that
\begin{align}
    \begin{array}{c}
	    \includeTikz{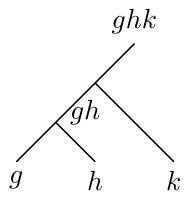}
	    {
	    	\begin{tikzpicture}[scale=0.4, every node/.style={scale=0.8}]
		    	\draw(0,0)--(3,3) node[below] at (0,0) {$g$} node[above] at (3,3) {$ghk$};
		    	\draw(2,0)--(1,1) node[below] at (2,0) {$h$};
		    	\draw(4,0)--(2,2) node[below] at (4,0) {$k$};
		    	\node[] at (1.75,1.25) {$gh$};
		    	\node[left] at(1,1){};
		    	\node[left] at(2,2){};
	    	\end{tikzpicture}
	    }
    \end{array}
    =
    \begin{array}{c}
    \includeTikz{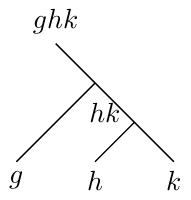}
    {
    	\begin{tikzpicture}[scale=0.4, every node/.style={scale=0.8}]
    	\draw(0,0)--(-3,3) node[below] at (0,0) {$k$} node[above] at (-3,3) {$ghk$};
    	\draw(-2,0)--(-1,1) node[below] at (-2,0) {$h$};
    	\draw(-4,0)--(-2,2) node[below] at (-4,0) {$g$};
    	\node at (-1.75,1.25) {$hk$};
    	\node[right] at(-1,1){};
    	\node[right] at(-2,2){};
    	\end{tikzpicture}
    }
    \end{array}
\end{align}
\end{example}

\begin{definition}[Bimodule category]
  Let $\mathcal{C},\,\mathcal{D}$ be fusion categories. A {\em bimodule category} $\mathcal{C} \curvearrowright \mathcal{M} \curvearrowleft \mathcal{D}$ is a semi-simple category equipped with functors $ - \vartriangleright - : \mathcal{C} \otimes \mathcal{M} \to \mathcal{M}$ and $ - \vartriangleleft - : \mathcal{M} \otimes \mathcal{D} \to \mathcal{M}$ and three natural isomorphisms
  \begin{align}
    &- \vartriangleright (- \vartriangleright -) \cong (- \otimes - ) \vartriangleright - \\
    & (- \vartriangleleft -) \vartriangleleft - \cong - \vartriangleleft (- \otimes -) \\
    & - \vartriangleright (- \vartriangleleft -) \cong (- \vartriangleright -) \vartriangleleft -.
  \end{align}
  If we choose bases for the Hom spaces $\mathcal{M}(c \vartriangleright m, n)$ and $\mathcal{M}(m \vartriangleleft d,n)$, then these natural isomorphisms can be represented as tensors: 
  
  \begin{align}
  \begin{array}{c}
	  \includeTikz{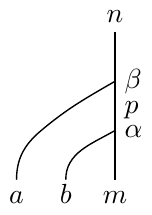}
	  {
	  	\begin{tikzpicture}[scale=0.5, every node/.style={scale=0.8}]
	  	\draw(0,0)--(0,3);
	  	\draw(-1,0)to[out=90,in=210](0,1);
	  	\draw(-2,0)to[out=90,in=220](-1.5,1)to[out=40,in=210](0,2);
	  	\node[below,inline text] at (-2,0) {$a$};
	  	\node[below,inline text] at (-1,0) {$b$};
	  	\node[below,inline text] at (0,0) {$m$};
	  	\node[above,inline text] at (0,3) {$n$};
	  	\node[right,inline text] at (0,1) {$\alpha$};
	  	\node[right,inline text] at (0,1.5) {$p$};
	  	\node[right,inline text] at (0,2) {$\beta$};
	  	\end{tikzpicture}
	  }
  \end{array}
  &=\sum_{(q,\mu,\nu)}
  \biggl[L^{abm}_{n}\biggr]_{(p,\alpha,\beta)(q,\mu,\nu)}
  \begin{array}{c}
  \includeTikz{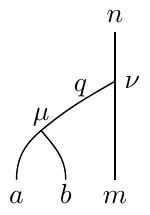}
  {
  	\begin{tikzpicture}[scale=0.5, every node/.style={scale=0.8}]
	  	\draw(0,0)--(0,3);
	  	\draw(-1,0)to[out=90,in=-50](-1.5,1);
	  	\draw(-2,0)to[out=90,in=220](-1.5,1)to[out=40,in=210](0,2);
	  	\node[below,inline text] at (-2,0) {$a$};
	  	\node[below,inline text] at (-1,0) {$b$};
	  	\node[below,inline text] at (0,0) {$m$};
	  	\node[above,inline text] at (0,3) {$n$};
	  	\node[above,inline text] at (-1.5,1) {$\mu$};
	  	\node[above,inline text] at (-.7,1.6) {$q$};
	  	\node[right,inline text] at (0,2) {$\nu$};
  	\end{tikzpicture}
  }
  \end{array}.\label{eqn:leftassociatorA} \\
  \begin{array}{c}
  \includeTikz{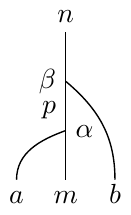}
  {
  	\begin{tikzpicture}[scale=0.5, every node/.style={scale=0.8}]
  	\draw(0,0)--(0,3);
  	\draw(1,0)to[out=90,in=-40](0,2);
  	\draw(-1,0)to[out=90,in=200](0,1);
  	\node[below,inline text] at (-1,0) {$a$};
  	\node[below,inline text] at (1,0) {$b$};
  	\node[below,inline text] at (0,0) {$m$};
  	\node[above,inline text] at (0,3) {$n$};
  	\node[right,inline text] at (0,1) {$\alpha$};
  	\node[left,inline text] at (0,1.5) {$p$};
  	\node[left,inline text] at (0,2) {$\beta$};
  	\end{tikzpicture}
  }
  \end{array}
  &=\sum_{(q,\mu,\nu)}
  \biggl[C^{abm}_{n}\biggr]_{(p,\alpha,\beta)(q,\mu,\nu)}
  \begin{array}{c}
  \includeTikz{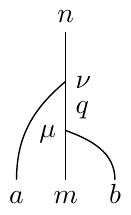}
  {
  	\begin{tikzpicture}[scale=0.5, every node/.style={scale=0.8}]
  	\draw(0,0)--(0,3);
  	\draw(1,0)to[out=90,in=-20](0,1);
  	\draw(-1,0)to[out=90,in=220](0,2);
  	\node[below,inline text] at (-1,0) {$a$};
  	\node[below,inline text] at (1,0) {$b$};
  	\node[below,inline text] at (0,0) {$m$};
  	\node[above,inline text] at (0,3) {$n$};
  	\node[left,inline text] at (0,1) {$\mu$};
  	\node[right,inline text] at (0,1.5) {$q$};
  	\node[right,inline text] at (0,2) {$\nu$};
  	\end{tikzpicture}
  }
  \end{array},\\
  \begin{array}{c}
  \includeTikz{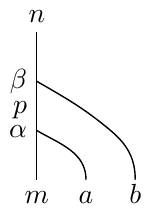}
  {
  	\begin{tikzpicture}[xscale=-0.5,yscale=.5, every node/.style={scale=0.8}]
  	\draw(0,0)--(0,3);
  	\draw(-1,0)to[out=90,in=210](0,1);
  	\draw(-2,0)to[out=90,in=220](-1.5,1)to[out=40,in=210](0,2);
  	\node[below,inline text] at (-2,0) {$b$};
  	\node[below,inline text] at (-1,0) {$a$};
  	\node[below,inline text] at (0,0) {$m$};
  	\node[above,inline text] at (0,3) {$n$};
  	\node[left,inline text] at (0,1) {$\alpha$};
  	\node[left,inline text] at (0,1.5) {$p$};
  	\node[left,inline text] at (0,2) {$\beta$};
  	\end{tikzpicture}
  }
  \end{array}
  &=\sum_{(q,\mu,\nu)}
  \biggl[R^{mab}_{n}\biggr]_{(p,\alpha,\beta)(q,\mu,\nu)}
  \begin{array}{c}
  \includeTikz{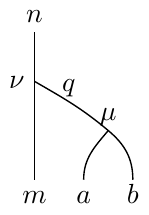}
  {
  	\begin{tikzpicture}[xscale=-0.5,yscale=.5, every node/.style={scale=0.8}]
  	\draw(0,0)--(0,3);
  	\draw(-1,0)to[out=90,in=-50](-1.5,1);
  	\draw(-2,0)to[out=90,in=220](-1.5,1)to[out=40,in=210](0,2);
  	\node[below,inline text] at (-2,0) {$b$};
  	\node[below,inline text] at (-1,0) {$a$};
  	\node[below,inline text] at (0,0) {$m$};
  	\node[above,inline text] at (0,3) {$n$};
  	\node[above,inline text] at (-1.5,1) {$\mu$};
  	\node[above,inline text] at (-.7,1.6) {$q$};
  	\node[left,inline text] at (0,2) {$\nu$};
  	\end{tikzpicture}
  }
  \end{array}.
  \end{align}
\end{definition}

\begin{example}[$\vvec{\mathbb{Z}/p\mathbb{Z}}$-bimodules]
  In \onlinecite{1806.01279}, following \onlinecite{MR3242743}, we gave a complete list of the simple $\vvec{\mathbb{Z}/p\mathbb{Z}}$-bimodules (reproduced in Table~\ref{tab:zpdata}). In this example, we can always gauge away the left and right associators, so we only tabulate the center associator in Table~\ref{tab:zpdata}. In \onlinecite{1806.01279}, we computed the physical interpretations for each of these bimodules. This data is reproduced in Table~\ref{tab:bimodinterp}.
 \end{example}

\begin{example}[$\vvec{\mathbb{Z}/p\mathbb{Z}}$-bimodule functors]
As explained by Morrison and Walker in Chapter 6 of \cite{MR2978449}, bimodule functors correspond to representations of annular categories from definition \ref{def:annular_category}. We use the notation $\defect{N}{M}{\bullet}{}{}$ to refer to primitive idempotents in $\vvec{\mathbb{Z}/p\mathbb{Z}}$-bimodule annular categories. All the primitive idempotents are reproduced in table \ref{tab:idempotents}.
\end{example}

\begin{table}[h]
%	\resizebox{\linewidth}{!}{
		%	\begin{tabular}{|c  |  c  |  c  |  c | c | c | c |}
		\begin{tabular}{!{\vrule width 1pt}>{\columncolor[gray]{.9}[\tabcolsep]}c!{\vrule width 1pt}c !{\color[gray]{.8}\vrule} c !{\color[gray]{.8}\vrule} c!{\vrule width 1pt}}
			\toprule[1pt]
			\rowcolor[gray]{.9}[\tabcolsep]Bimodule label & Left action & Right action & Associator \\
			\toprule[1pt]
			$T$ &
			$
			\begin{array}{c}
			\includeTikz{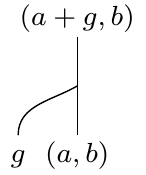}
			{
				\begin{tikzpicture}
				\Laction[$g$][$(a,b)$][$(a+g,b)$];
				\end{tikzpicture}
			}
			\end{array}
			$
			&
			$
			\begin{array}{c}
			\includeTikz{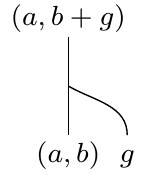}
			{
				\begin{tikzpicture}
				\Raction[$g$][$(a,b)$][$(a,b+g)$];
				\end{tikzpicture}
			}
			\end{array}
			$&
			$
			\begin{array}{c}
			\includeTikz{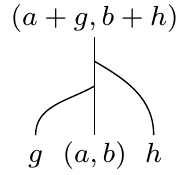}
			{
				\begin{tikzpicture}
				\Lassociator[$g$][$h$][$(a,b)$][$(a+g,b+h)$];
				\end{tikzpicture}
			}
			\end{array}
			=
			\begin{array}{c}
			\includeTikz{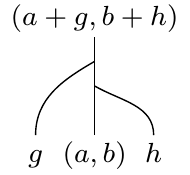}
			{
				\begin{tikzpicture}
				\Rassociator[$g$][$h$][$(a,b)$][$(a+g,b+h)$];
				\end{tikzpicture}
			}
			\end{array}
			$
			\\
			\greycline{2-3}
			$L$ &
			$
			\begin{array}{c}
			\includeTikz{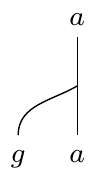}
			{
				\begin{tikzpicture}
				\Laction[$g$][$a$][$a$];
				\end{tikzpicture}
			}
			\end{array}
			$
			&
			$
			\begin{array}{c}
			\includeTikz{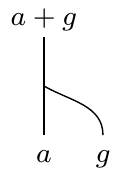}
			{
				\begin{tikzpicture}
				\Raction[$g$][$a$][$a+g$];
				\end{tikzpicture}
			}
			\end{array}
			$&
			$
			\begin{array}{c}
			\includeTikz{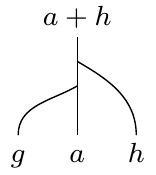}
			{
				\begin{tikzpicture}
				\Lassociator[$g$][$h$][$a$][$a+h$];
				\end{tikzpicture}
			}
			\end{array}
			=
			\begin{array}{c}
			\includeTikz{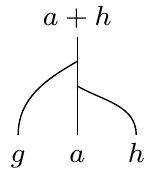}
			{
				\begin{tikzpicture}
				\Rassociator[$g$][$h$][$a$][$a+h$];
				\end{tikzpicture}
			}
			\end{array}
			$
			\\
			\greycline{2-3} 
			$R$ &
			$
			\begin{array}{c}
			\includeTikz{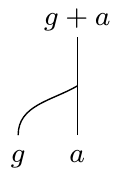}
			{
				\begin{tikzpicture}
				\Laction[$g$][$a$][$g+a$];
				\end{tikzpicture}
			}
			\end{array}
			$
			&
			$
			\begin{array}{c}
			\includeTikz{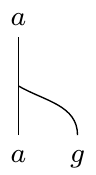}
			{
				\begin{tikzpicture}
				\Raction[$g$][$a$][$a$];
				\end{tikzpicture}
			}
			\end{array}
			$&
			$
			\begin{array}{c}
			\includeTikz{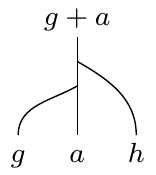}
			{
				\begin{tikzpicture}
				\Lassociator[$g$][$h$][$a$][$g+a$];
				\end{tikzpicture}
			}
			\end{array}
			=
			\begin{array}{c}
			\includeTikz{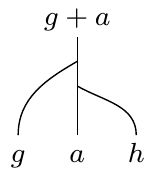}
			{
				\begin{tikzpicture}
				\Rassociator[$g$][$h$][$a$][$g+a$];
				\end{tikzpicture}
			}
			\end{array}
			$
			\\
			\greycline{2-3} 
			$F_0$ &
			$
			\begin{array}{c}
			\includeTikz{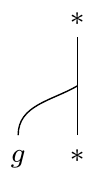}
			{
				\begin{tikzpicture}
				\Laction[$g$][$*$][$*$];
				\end{tikzpicture}
			}
			\end{array}
			$
			&
			$
			\begin{array}{c}
			\includeTikz{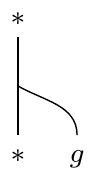}
			{
				\begin{tikzpicture}
				\Raction[$g$][$*$][$*$];
				\end{tikzpicture}
			}
			\end{array}
			$&
			$
			\begin{array}{c}
			\includeTikz{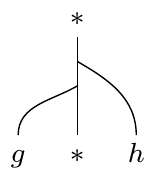}
			{
				\begin{tikzpicture}
				\Lassociator[$g$][$h$][$*$][$*$];
				\end{tikzpicture}
			}
			\end{array}
			=
			\begin{array}{c}
			\includeTikz{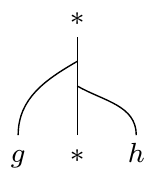}
			{
				\begin{tikzpicture}
				\Rassociator[$g$][$h$][$*$][$*$];
				\end{tikzpicture}
			}
			\end{array}
			$
			\\		
			\toprule[1pt]
			$X_k$ &
			$
			\begin{array}{c}
			\includeTikz{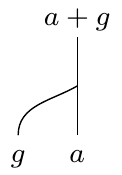}
			{
				\begin{tikzpicture}
				\Laction[$g$][$a$][$a+g$];
				\end{tikzpicture}
			}
			\end{array}
			$
			&
			$
			\begin{array}{c}
			\includeTikz{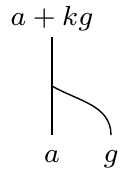}
			{
				\begin{tikzpicture}
				\Raction[$g$][$a$][$a+kg$];
				\end{tikzpicture}
			}
			\end{array}
			$&
			$
			\begin{array}{c}
			\includeTikz{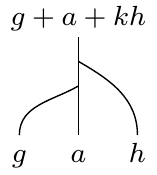}
			{
				\begin{tikzpicture}
				\Lassociator[$g$][$h$][$a$][$g+a+kh$];
				\end{tikzpicture}
			}
			\end{array}
			=
			\begin{array}{c}
			\includeTikz{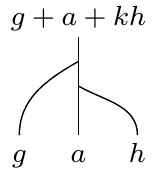}
			{
				\begin{tikzpicture}
				\Rassociator[$g$][$h$][$a$][$g+a+kh$];
				\end{tikzpicture}
			}
			\end{array}
			$
			\\
			\greycline{2-3}
			$F_q$&
			$
			\begin{array}{c}
			\includeTikz{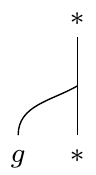}
			{
				\begin{tikzpicture}
				\Laction[$g$][$*$][$*$];
				\end{tikzpicture}
			}
			\end{array}
			$
			&
			$
			\begin{array}{c}
			\includeTikz{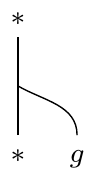}
			{
				\begin{tikzpicture}
				\Raction[$g$][$*$][$*$];
				\end{tikzpicture}
			}
			\end{array}
			$&
			$
			\begin{array}{c}
			\includeTikz{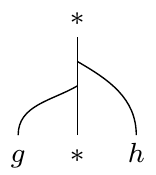}
			{
				\begin{tikzpicture}
				\Lassociator[$g$][$h$][$*$][$*$];
				\end{tikzpicture}
			}
			\end{array}
			=
			e^{\frac{2\pi i}{p} q g h}
			\begin{array}{c}
			\includeTikz{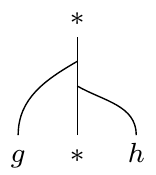}
			{
				\begin{tikzpicture}
				\Rassociator[$g$][$h$][$*$][$*$];
				\end{tikzpicture}
			}
			\end{array}
			$
			\\
			\hline
		\end{tabular}
%	}
	\caption{Data for all $\vvec{\ZZ{p}}-\vvec{\ZZ{p}}$ bimodules. $q \in H^2(\ZZ{p},U(1)) \cong \ZZ{p}$. Bimodules below the thick line are invertible. Reproduced from \onlinecite{1806.01279}.}\label{tab:zpdata}
\end{table}

\begin{table}[h]
	\begin{tabular}{!{\vrule width 1pt}c!{\vrule width 1pt}c|c!{\vrule width 1pt}}
		\toprule[1pt]
		\rowcolor[gray]{.9}[\tabcolsep]Bimodule label & Domain wall & Action on particles\\
		\toprule[1pt]
		$T$&$\begin{array}{c}\includeTikz{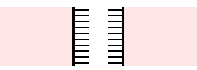}{
			\begin{tikzpicture}[yscale=.3]
			\draw[white](0,-1.1)--(0,1.21);
			\fill[red!10](-1,-1) rectangle (-.25,1);
			\fill[red!10](1,-1) rectangle (.25,1);
			\draw[thick](-.25,-1)--(-.25,1);
			\draw[thick](.25,-1)--(.25,1);
			\foreach \x in {0,...,9}{\draw (.1,-.9+.2*\x)--(.25,-.9+.2*\x);};
			\foreach \x in {0,...,9}{\draw (-.1,-.9+.2*\x)--(-.25,-.9+.2*\x);};
			\end{tikzpicture}}
		\end{array}$&Condenses $e$ on both sides\\
		%			\hline
		$L$&$\begin{array}{c}\includeTikz{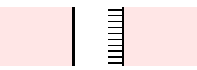}{\begin{tikzpicture}[yscale=.3]
			\draw[white](0,-1.1)--(0,1.21);
			\fill[red!10](-1,-1) rectangle (-.25,1);
			\fill[red!10](1,-1) rectangle (.25,1);
			\draw[thick](-.25,-1)--(-.25,1);
			\draw[thick](.25,-1)--(.25,1);
			\foreach \x in {0,...,9}{\draw (.1,-.9+.2*\x)--(.25,-.9+.2*\x);};
			\end{tikzpicture}}\end{array}$&Condenses $m$ on left and $e$ on right\\
		%			\hline
		$R$&$\begin{array}{c}\includeTikz{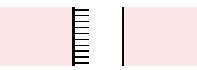}{\begin{tikzpicture}[yscale=.3]
			\draw[white](0,-1.1)--(0,1.21);
			\fill[red!10](-1,-1) rectangle (-.25,1);
			\fill[red!10](1,-1) rectangle (.25,1);
			\draw[thick](-.25,-1)--(-.25,1);
			\draw[thick](.25,-1)--(.25,1);
			\foreach \x in {0,...,9}{\draw (-.1,-.9+.2*\x)--(-.25,-.9+.2*\x);};
			\end{tikzpicture}}\end{array}$&Condenses $e$ on left and $m$ on right\\
		%			\hline
		$F_0$&$\begin{array}{c}\includeTikz{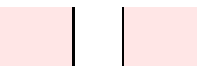}{\begin{tikzpicture}[yscale=.3]
			\draw[white](0,-1.1)--(0,1.21);
			\fill[red!10](-1,-1) rectangle (-.25,1);
			\fill[red!10](1,-1) rectangle (.25,1);
			\draw[thick](-.25,-1)--(-.25,1);
			\draw[thick](.25,-1)--(.25,1);
			\end{tikzpicture}}\end{array}$&Condenses $m$ on both sides\\
		\toprule[1pt]
		$X_k$&$\begin{array}{c}\includeTikz{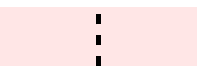}{\begin{tikzpicture}[yscale=.3]
			\draw[white](0,-1.1)--(0,1.21);
			\fill[red!10](-1,-1) rectangle (1,1);
			\draw[ultra thick,dashed] (0,-1)--(0,1);
			\end{tikzpicture}}\end{array}$&$X_k:e^am^b\mapsto e^{ka}m^{k^{-1} b}$\\
		%			\hline
		$F_{q}=F_1 X_q$&$\begin{array}{c}\includeTikz{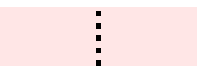}{\begin{tikzpicture}[yscale=.3]
			\draw[white](0,-1.1)--(0,1.21);
			\fill[red!10](-1,-1) rectangle (1,1);
			\draw[ultra thick,dotted] (0,-1)--(0,1);\end{tikzpicture}}\end{array}$&
		$F_1:e^am^b\mapsto e^{b}m^{a}$\\
		\toprule[1pt]
	\end{tabular}
	\caption{Domain walls on the lattice corresponding to bimodules. Reproduced from \onlinecite{1806.01279}.}\label{tab:bimodinterp}
\end{table}

\clearpage
\vspace*{20mm}
	\begin{table}[h]
\begin{minipage}[t][.8\textheight][c]{\textwidth}
	\rotatebox[]{90}{
\resizebox{.98\textheight}{!}{
			\begin{tabular}{!{\vrule width 1pt}>{\columncolor[gray]{.9}[\tabcolsep]}c!{\vrule width 1pt}c!{\color[gray]{.8}\vrule}c!{\color[gray]{.8}\vrule}c!{\color[gray]{.8}\vrule}c!{\color[gray]{.8}\vrule}c!{\color[gray]{.8}\vrule}c!{\vrule width 1pt}}
				\toprule[1pt]
				\rowcolor[gray]{.9}[\tabcolsep]&$T$&$L$&$R$&$F_0$&$X_l$&$F_r$\\
				\toprule[1pt]
				$T$&$\defect{T}{T}{a}{b}{}=
				\begin{array}{c}
				\includeTikz{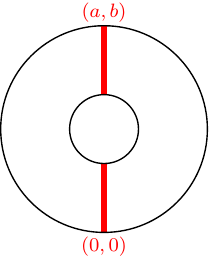}
				{
					\begin{tikzpicture}[scale=.7,every node/.style={scale=.7}]
					\annparamss{$(0,0)$}{$(a,b)$}{}{};
					\annss{}{};
					\end{tikzpicture}
				}
				\end{array}$&$\defect{T}{L}{a}{}{}=
				\begin{array}{c}
				\includeTikz{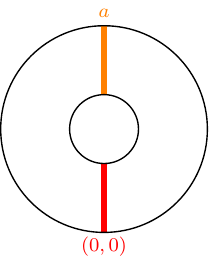}
				{
					\begin{tikzpicture}[scale=.7,every node/.style={scale=.7}]
					\annparamst{$(0,0)$}{$a$}{}{};
					\annst{}{};
					\end{tikzpicture}
				}
				\end{array}$&$\defect{T}{R}{a}{}{}=
				\begin{array}{c}
				\includeTikz{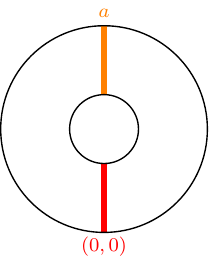}
				{
					\begin{tikzpicture}[scale=.7,every node/.style={scale=.7}]
					\annparamst{$(0,0)$}{$a$}{}{};
					\annst{}{};
					\end{tikzpicture}
				}
				\end{array}$&$\defect{T}{F_0}{}{}{}=
				\begin{array}{c}
				\includeTikz{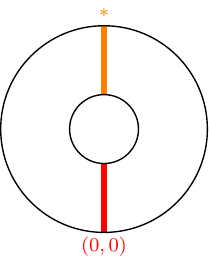}
				{
					\begin{tikzpicture}[scale=.7,every node/.style={scale=.7}]
					\annparamst{$(0,0)$}{$*$}{}{};
					\annst{}{};
					\end{tikzpicture}
				}
				\end{array}$&$\defect{T}{X_l}{a}{}{}=
				\begin{array}{c}
				\includeTikz{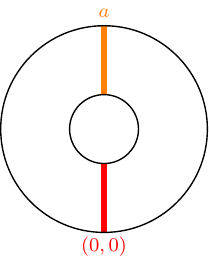}
				{
					\begin{tikzpicture}[scale=.7,every node/.style={scale=.7}]
					\annparamst{$(0,0)$}{$a$}{}{};
					\annst{}{};
					\end{tikzpicture}
				}
				\end{array}$&$\defect{T}{F_r}{}{}{}=
				\begin{array}{c}
				\includeTikz{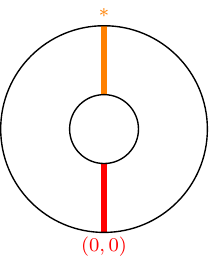}
				{
					\begin{tikzpicture}[scale=.7,every node/.style={scale=.7}]
					\annparamst{$(0,0)$}{$*$}{}{};
					\annst{}{};
					\end{tikzpicture}
				}
				\end{array}$\\
\greycline{2-7}
				$L$&
				$\defect{L}{T}{a}{}{}=
				\begin{array}{c}
				\includeTikz{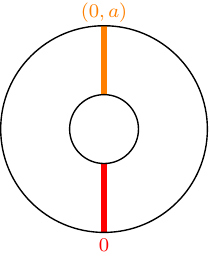}
				{
					\begin{tikzpicture}[scale=.7,every node/.style={scale=.7}]
					\annparamst{$0$}{$(0,a)$}{}{};
					\annst{}{};
					\end{tikzpicture}
				}
				\end{array}$
				&$\defect{L}{L}{a}{x}{}=\frac{1}{p}\sum_{g}\omega^{gx}
				\begin{array}{c}
				\includeTikz{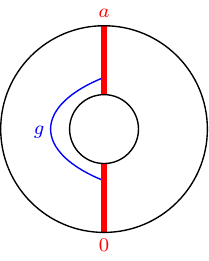}
				{
					\begin{tikzpicture}[scale=.7,every node/.style={scale=.7}]
					\annparamss{$0$}{$a$}{}{};
					\annss{$g$}{};
					\end{tikzpicture}
				}
				\end{array}$&$\defect{L}{R}{}{}{}=
				\begin{array}{c}
				\includeTikz{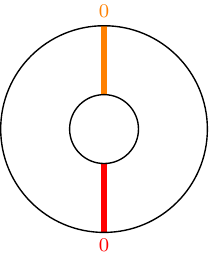}
				{
					\begin{tikzpicture}[scale=.7,every node/.style={scale=.7}]
					\annparamst{$0$}{$0$}{}{};
					\annst{}{};
					\end{tikzpicture}
				}
				\end{array}$&$\defect{L}{F_0}{x}{}{}=\frac{1}{p}\sum_g \omega^{gx}
				\begin{array}{c}
				\includeTikz{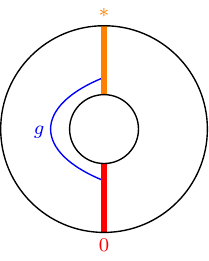}
				{
					\begin{tikzpicture}[scale=.7,every node/.style={scale=.7}]
					\annparamst{$0$}{$*$}{}{};
					\annst{$g$}{};
					\end{tikzpicture}
				}
				\end{array}$&$\defect{L}{X_l}{}{}{}=
				\begin{array}{c}
				\includeTikz{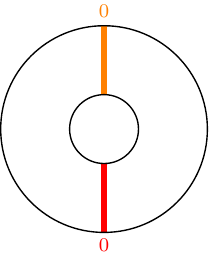}
				{
					\begin{tikzpicture}[scale=.7,every node/.style={scale=.7}]
					\annparamst{$0$}{$0$}{}{};
					\annst{}{};
					\end{tikzpicture}
				}
				\end{array}$&$\defect{L}{F_r}{x}{}{}=\frac{1}{p}\sum_g \omega^{gx}
				\begin{array}{c}
				\includeTikz{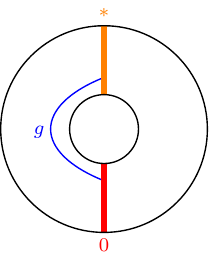}
				{
					\begin{tikzpicture}[scale=.7,every node/.style={scale=.7}]
					\annparamst{$0$}{$*$}{}{};
					\annst{$g$}{};
					\end{tikzpicture}
				}
				\end{array}$\\
\greyhline
				$R$&$\defect{R}{T}{a}{}{}=
				\begin{array}{c}
				\includeTikz{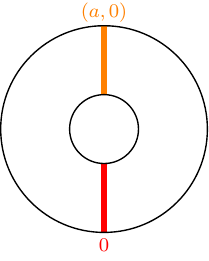}
				{
					\begin{tikzpicture}[scale=.7,every node/.style={scale=.7}]
					\annparamst{$0$}{$(a,0)$}{}{};
					\annst{}{};
					\end{tikzpicture}
				}
				\end{array}$&$\defect{R}{L}{}{}{}=
				\begin{array}{c}
				\includeTikz{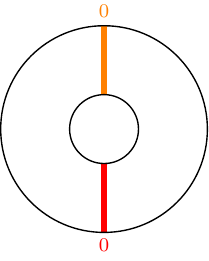}
				{
					\begin{tikzpicture}[scale=.7,every node/.style={scale=.7}]
					\annparamst{$0$}{$0$}{}{};
					\annst{}{};
					\end{tikzpicture}
				}
				\end{array}$&$\defect{R}{R}{a}{x}{}=\frac{1}{p}\sum_g\omega^{gx}
				\begin{array}{c}
				\includeTikz{RR_idempotent}
				{
					\begin{tikzpicture}[scale=.7,every node/.style={scale=.7}]
					\annparamss{$0$}{$a$}{}{};
					\annss{}{$-g$};
					\end{tikzpicture}
				}
				\end{array}$&$\defect{R}{F_0}{x}{}{}=\frac{1}{p}\sum_g\omega^{gx}
				\begin{array}{c}
				\includeTikz{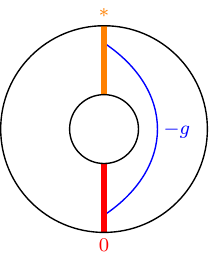}
				{
					\begin{tikzpicture}[scale=.7,every node/.style={scale=.7}]
					\annparamst{$0$}{$*$}{}{};
					\annst{}{$-g$};
					\end{tikzpicture}
				}
				\end{array}$&$\defect{R}{X_l}{}{}{}=
				\begin{array}{c}
				\includeTikz{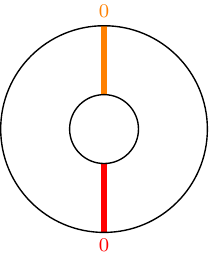}
				{
					\begin{tikzpicture}[scale=.7,every node/.style={scale=.7}]
					\annparamst{$0$}{$0$}{}{};
					\annst{}{};
					\end{tikzpicture}
				}
				\end{array}$&$\defect{R}{F_r}{x}{}{}=\frac{1}{p}\sum_g\omega^{gx}
				\begin{array}{c}
				\includeTikz{RFr_idempotent}
				{
					\begin{tikzpicture}[scale=.7,every node/.style={scale=.7}]
					\annparamst{$0$}{$*$}{}{};
					\annst{}{$-g$};
					\end{tikzpicture}
				}
				\end{array}$\\
\greycline{2-7}
				$F_0$&$\defect{F_0}{T}{}{}{}=
				\begin{array}{c}
				\includeTikz{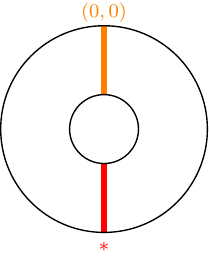}
				{
					\begin{tikzpicture}[scale=.7,every node/.style={scale=.7}]
					\annparamst{$*$}{$(0,0)$}{}{};
					\annst{}{};
					\end{tikzpicture}
				}
				\end{array}$&$\defect{F_0}{L}{x}{}{}=\frac{1}{p}\sum_g\omega^{gx}
				\begin{array}{c}
				\includeTikz{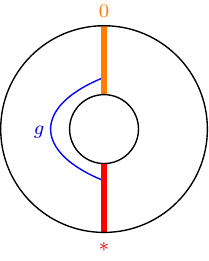}
				{
					\begin{tikzpicture}[scale=.7,every node/.style={scale=.7}]
					\annparamst{$*$}{$0$}{}{};
					\annst{$g$}{};
					\end{tikzpicture}
				}
				\end{array}$&$\defect{F_0}{R}{x}{}{}=\frac{1}{p}\sum_g\omega^{gx}
				\begin{array}{c}
				\includeTikz{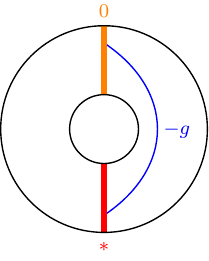}
				{
					\begin{tikzpicture}[scale=.7,every node/.style={scale=.7}]
					\annparamst{$*$}{$0$}{}{};
					\annst{}{$-g$};
					\end{tikzpicture}
				}
				\end{array}$&$\defect{F_0}{F_0}{x}{y}{}=\frac{1}{p^2}\sum_{g,h}\omega^{gx+hy}
				\begin{array}{c}
				\includeTikz{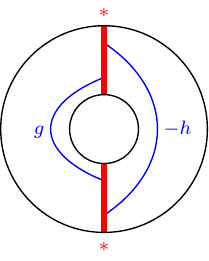}
				{
					\begin{tikzpicture}[scale=.7,every node/.style={scale=.7}]
					\annparamss{$*$}{$*$}{}{};
					\annss{$g$}{$-h$};
					\end{tikzpicture}
				}
				\end{array}$&$\defect{F_0}{X_l}{x}{}{}=\frac{1}{p}\sum_{g}\omega^{gx}
				\begin{array}{c}
				\includeTikz{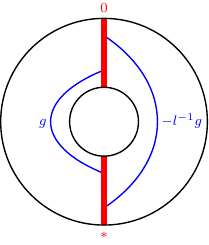}
				{
					\begin{tikzpicture}[scale=.7,every node/.style={scale=.5}]
					\annparamss{$*$}{$0$}{}{};
					\annss{$g$}{$-l^{-1}g$};
					\end{tikzpicture}
				}
				\end{array}$&$\defect{F_0}{F_r}{}{}{}=\frac{1}{p}\sum_{g}
				\begin{array}{c}
				\includeTikz{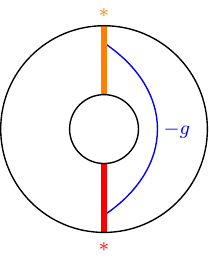}
				{
					\begin{tikzpicture}[scale=.7,every node/.style={scale=.7}]
					\annparamst{$*$}{$*$}{}{};
					\annst{}{$-g$};
					\end{tikzpicture}
				}
				\end{array}$\\
\greycline{2-7}
				$X_k$&$\defect{X_k}{T}{a}{}{}=
				\begin{array}{c}
				\includeTikz{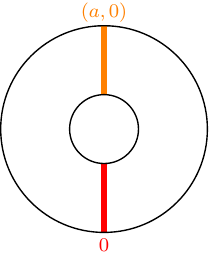}
				{
					\begin{tikzpicture}[scale=.7,every node/.style={scale=.7}]
					\annparamst{$0$}{$(a,0)$}{}{};
					\annst{}{};
					\end{tikzpicture}
				}
				\end{array}$&$\defect{X_k}{L}{}{}{}=
				\begin{array}{c}
				\includeTikz{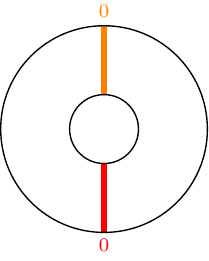}
				{
					\begin{tikzpicture}[scale=.7,every node/.style={scale=.7}]
					\annparamst{$0$}{$0$}{}{};
					\annst{}{};
					\end{tikzpicture}
				}
				\end{array}$&$\defect{X_k}{R}{}{}{}=
				\begin{array}{c}
				\includeTikz{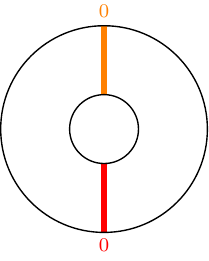}
				{
					\begin{tikzpicture}[scale=.7,every node/.style={scale=.7}]
					\annparamst{$0$}{$0$}{}{};
					\annst{}{};
					\end{tikzpicture}
				}
				\end{array}$&$\defect{X_k}{F_0}{x}{}{}=\frac{1}{p}\sum_g\omega^{gx}
				\begin{array}{c}
				\includeTikz{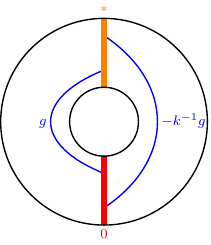}
				{
					\begin{tikzpicture}[scale=.7,every node/.style={scale=.5}]
					\annparamst{$0$}{$*$}{}{};
					\annst{$g$}{$-k^{-1}g$};
					\end{tikzpicture}
				}
				\end{array}$&
				\begin{tabular}{c}
					$\defect{X_k}{X_k}{a}{x}{}=\frac{1}{p}\sum_g\omega^{gx}
					\begin{array}{c}
					\includeTikz{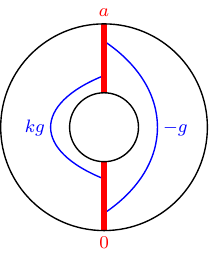}
					{
						\begin{tikzpicture}[scale=.7,every node/.style={scale=.65}]
						\annparamss{$0$}{$a$}{}{};
						\annss{$kg$}{$-g$};
						\end{tikzpicture}
					}
					\end{array}$
					\\
					\greyhline
					$\defect{X_k}{X_l}{}{}{}=
					\begin{array}{c}
					\includeTikz{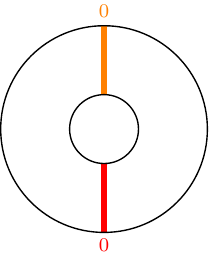}
					{
						\begin{tikzpicture}[scale=.7,every node/.style={scale=.7}]
						\annparamst{$0$}{$0$}{}{};
						\annst{}{};
						\end{tikzpicture}
					}
					\end{array}$
				\end{tabular}
				&	
				$\defect{X_k}{F_r}{x}{}{}=\frac{1}{p}\sum_g\Theta_{x,kr}(g)
				\begin{array}{c}
				\includeTikz{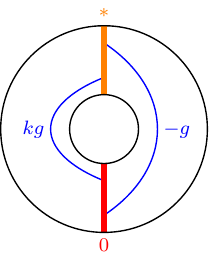}
				{
					\begin{tikzpicture}[scale=.7,every node/.style={scale=.7}]
					\annparamst{$0$}{$*$}{}{};
					\annst{$kg$}{$-g$};
					\end{tikzpicture}
				}
				\end{array}$
				\\
\greycline{2-7}
				$F_q$&$\defect{F_q}{T}{}{}{}=
				\begin{array}{c}
				\includeTikz{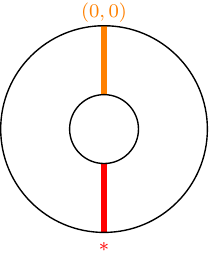}
				{
					\begin{tikzpicture}[scale=.7,every node/.style={scale=.7}]
					\annparamst{$*$}{$(0,0)$}{}{};
					\annst{}{};
					\end{tikzpicture}
				}
				\end{array}$&$\defect{F_q}{L}{x}{}{}=\frac{1}{p}\sum_g\omega^{gx}
				\begin{array}{c}
				\includeTikz{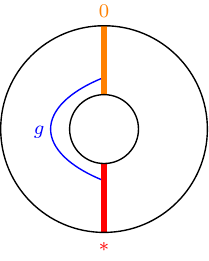}
				{
					\begin{tikzpicture}[scale=.7,every node/.style={scale=.7}]
					\annparamst{$*$}{$0$}{}{};
					\annst{$g$}{};
					\end{tikzpicture}
				}
				\end{array}$&$\defect{F_q}{R}{x}{}{}=\frac{1}{p}\sum_g\omega^{gx}
				\begin{array}{c}
				\includeTikz{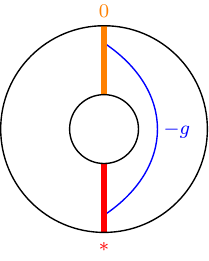}
				{
					\begin{tikzpicture}[scale=.7,every node/.style={scale=.7}]
					\annparamst{$*$}{$0$}{}{};
					\annst{}{$-g$};
					\end{tikzpicture}
				}
				\end{array}$&$\defect{F_q}{F_0}{}{}{}=\frac{1}{p}\sum_{g}
				\begin{array}{c}
				\includeTikz{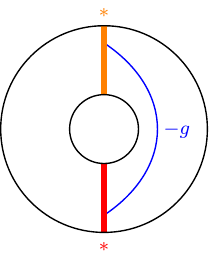}
				{
					\begin{tikzpicture}[scale=.7,every node/.style={scale=.7}]
					\annparamst{$*$}{$*$}{}{};
					\annst{}{$-g$};
					\end{tikzpicture}
				}
				\end{array}$&
				$\defect{F_q}{X_l}{x}{}{}=\frac{1}{p}\sum_g\Theta_{x,-ql}(g)
				\begin{array}{c}
				\includeTikz{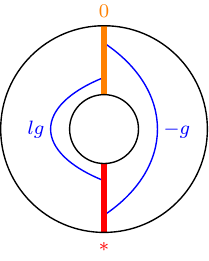}
				{
					\begin{tikzpicture}[scale=.7,every node/.style={scale=.7}]
					\annparamst{$*$}{$0$}{}{};
					\annst{$lg$}{$-g$};
					\end{tikzpicture}
				}
				\end{array}$
				&
				\begin{tabular}{c}
					$\defect{F_q}{F_q}{x}{y}{}=\frac{1}{p^2}\sum_{g,h}\omega^{gx+hy}
					\begin{array}{c}
					\includeTikz{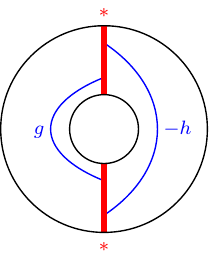}
					{
						\begin{tikzpicture}[scale=.7,every node/.style={scale=.7}]
						\annparamss{$*$}{$*$}{}{};
						\annss{$g$}{$-h$};
						\end{tikzpicture}
					}
					\end{array}$\\
					\greyhline
					$\defect{F_q}{F_r}{}{}{}=\frac{1}{p}\sum_{g}
					\begin{array}{c}
					\includeTikz{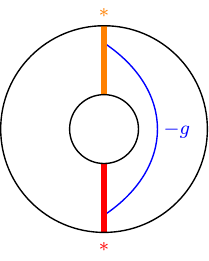}
					{
						\begin{tikzpicture}[scale=.7,every node/.style={scale=.7}]
						\annparamst{$*$}{$*$}{}{};
						\annst{}{$-g$};
						\end{tikzpicture}
					}
					\end{array}$
				\end{tabular}
				\\
				\toprule[1pt]
			\end{tabular}
		}
	}
	\caption{Indecomposable idempotents for 2-string annuli of all domain walls, corresponding to defects. For $p=2$, $\Theta_{x,a}(g)=(-1)^{gx} i^{ag}$, whilst for odd $p$ $\Theta_{x,a}(g)=\omega^{gx+ag^2 2^{-1}}$, where $2^{-1}$ is the modular inverse of 2.  Reproduced from \onlinecite{1810.09469}.}\label{tab:idempotents}
\end{minipage}
		%}
	\end{table}
%\end{minipage}

% This line sets the project root file.
% !TEX root = ../defects_domain_wall_structures.tex
% !TEX spellcheck = en_US
\clearpage
\section{Relationship to Extension Theory} \label{sec:extension_theory}

The computations described in this paper are closely related to extension 
theory as described in \onlinecite{MR2677836}.
\begin{definition} \label{def:BPR}
	We define the Brauer-Picard 3-category ${\bf BPR}$ as follows: Objects are fusion categories, 1-morphisms are bimodules, 2-morphisms are bimodule functors and 3-morphisms are natural transformations.
\end{definition}
Let $M$ be a finite monoid and $\cat{M}$ the tensor category ${\bf Vec}(M)$ of $M$-graded vector spaces considered as a 3-category with a single object $*$ and only identity 3-morphisms. Then {\bf extension data} is exactly a 3-functor $\cat{M} \to {\bf BPR}$. Such a 3-functor contains the following data:
\begin{itemize}
\item A fusion category $* \mapsto \cat{C}$.
\item A $\cat{C}{-}\cat{C}$ bimodule $g \mapsto M_g$ for each element $g \in M$.
\item Annular category representations
  \begin{align}
  \begin{array}{c}
    \includeTikz{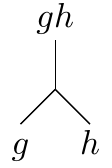}{
      \begin{tikzpicture}[scale=0.5]
        \trivalentvertex{$gh$}{$g$}{$h$}{}
      \end{tikzpicture}
      }
  \end{array} \mapsto \Lambda_{g,h} \in {\bf Rep}\left(
  \begin{array}{c}
    \includeTikz{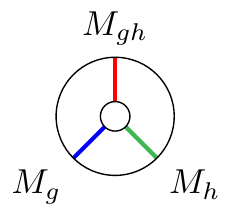}{
      \begin{tikzpicture}[scale=0.3]
        \draw[very thick,red] (0,0) -- (0,2);
        \draw[very thick,blue] (0,0) -- (-1.4142,-1.4142);
        \draw[very thick,nicegreen] (0,0) -- (1.4142,-1.4142);
        \filldraw[fill=white] (0,0) circle (0.5);
        \draw (0,0) circle (2);
        \node[above] at (0,2) {$M_{gh}$};
        \node[below left] at (-1.4142,-1.4142) {$M_g$};
        \node[below right] at (1.4142,-1.4142) {$M_h$};
        \end{tikzpicture}
    }
    \end{array}
  \right),\\
  \begin{array}{c}
    \includeTikz{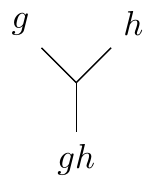}{
      \begin{tikzpicture}[scale=0.5]
        \draw (0,0) -- (0,-1);
        \draw (0,0) -- (-0.7071,0.7071);
        \draw (0,0) -- (0.7071,0.7071);
        \node[below] at (0,-1) {$gh$};
        \node[above left] at (-0.7071,0.7071) {$g$};
        \node[above right] at (0.7071,0.7071) {$h$}; 
      \end{tikzpicture}
      }
  \end{array} \mapsto V_{g,h} \in {\bf Rep}\left(
  \begin{array}{c}
    \includeTikz{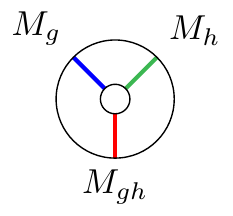}{
      \begin{tikzpicture}[scale=0.3,yscale=-1]
        \draw[very thick,red] (0,0) -- (0,2);
        \draw[very thick,blue] (0,0) -- (-1.4142,-1.4142);
        \draw[very thick,nicegreen] (0,0) -- (1.4142,-1.4142);
        \filldraw[fill=white] (0,0) circle (0.5);
        \draw (0,0) circle (2);
        \node[below] at (0,2) {$M_{gh}$};
        \node[above left] at (-1.4142,-1.4142) {$M_g$};
        \node[above right] at (1.4142,-1.4142) {$M_h$};
        \end{tikzpicture}
    }
    \end{array}
  \right).
  \end{align}\label{eqn:annrep}
\end{itemize}
In order for the 3-functor to be defined at this level, the following diagrams must map to the identity defect
  \begin{align}
  \begin{array}{c}
    \includeTikz{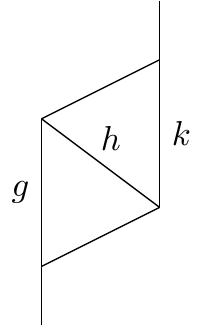}{
      \begin{tikzpicture}[scale=0.3]
        \draw (0,0)--(0,2);
  \draw (0,2)--(0,7) node[midway,left] {$g$};
  \draw (0,2)--(4,4);
  \draw (4,4)--(0,7) node[midway,above,xshift=3pt] {$h$};
  \draw (0,7)--(4,9);
  \draw (4,4)--(4,9) node[midway,right] {$k$};
  \draw (4,9)--(4,11);
      \end{tikzpicture}
    }
  \end{array}
  &\mapsto
\begin{array}{c}
    \includeTikz{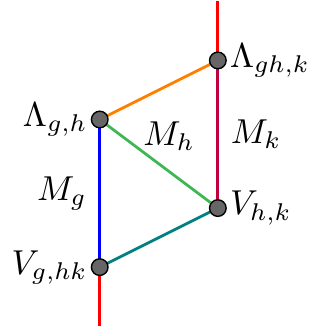}{
      \begin{tikzpicture}[scale=0.3]
\draw [thick,red] (0,0)--(0,2);
  \draw [thick,blue] (0,2)--(0,7) node[black,midway,left] {$M_g$};
  \draw [thick,teal] (0,2)--(4,4);
  \draw [thick,nicegreen] (4,4)--(0,7) node[black,midway,above,xshift=3pt] {$M_h$};
  \draw [thick,orange] (0,7)--(4,9);
  \draw [thick,purple] (4,4)--(4,9) node[black,midway,right] {$M_k$};
  \draw [thick,red] (4,9)--(4,11);
  \filldraw[fill=dandark, draw=black] (0,2) circle (8pt);
  \filldraw[fill=dandark, draw=black] (4,4) circle (8pt);
  \filldraw[fill=dandark, draw=black] (0,7) circle (8pt);
  \filldraw[fill=dandark, draw=black] (4,9) circle (8pt);
  \node[left] at (0,2) {$V_{g,hk}$};
  \node[right] at (4,4) {$V_{h,k}$};
  \node[left] at (0,7) {$\Lambda_{g,h}$};
  \node[right] at (4,9) {$\Lambda_{gh,k}$};
      \end{tikzpicture}
    }
    \end{array},\label{eqn:associator1} \\
\begin{array}{c}
\includeTikz{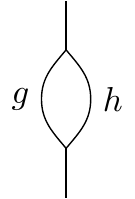}
{
  \begin{tikzpicture}[scale=0.5]
    \draw (0,-1) to[out=130,in=270] (-.5,0) to[out=90,in=230] (0,1);
\draw (0,1)--(0,2);
\draw (0,-1) to[out=90-40,in=270] (.5,0) to[out=90,in=270+40] (0,1);
\draw (0,-2)--(0,-1);
\node[left] at (-.5,0) {$g$};
\node[right] at (.5,0) {$h$};
\end{tikzpicture}
}
\end{array}
&\mapsto
\begin{array}{c}
\includeTikz{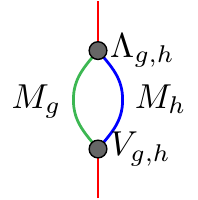}
{
  \begin{tikzpicture}[scale=0.5]
\draw[thick,nicegreen] (0,-1) to[out=130,in=270] (-.5,0) to[out=90,in=230] (0,1);
\draw[thick,red] (0,1)--(0,2);
\draw[thick,blue] (0,-1) to[out=90-40,in=270] (.5,0) to[out=90,in=270+40] (0,1);
\draw[thick,red] (0,-2)--(0,-1);
	\filldraw[fill=dandark] (0,-1) circle (.18);\filldraw[fill=dandark] (0,1) circle (.18);
	\node[inline text,left] at (-.5,0) {$M_g$};\node[inline text,right] at (.5,0) {$M_h$};
  \node[right] at (0,-1) {$V_{g,h}$};
  \node[right] at (0,1) {$\Lambda_{g,h}$};
\end{tikzpicture}
}
\end{array}.
\end{align}
This is closely related to the vanishing of the $O_3$ obstruction from \onlinecite{MR2677836}. For invertible bimodules $g,h,k$, the annular category representations Eqn.~\ref{eqn:annrep} define bimodule equivalences $M_g\otimes_\cat{C}M_g\cong M_{gh}$ (In Section 8 of \onlinecite{MR2677836}, these equivalences are called $M_{g,h}$). The domain wall structure Eqn.~\ref{eqn:associator1} then corresponds to $T_{g,h,k}$ in \onlinecite{MR2677836}.

 There are further obstructions called $O_4$ in \onlinecite{MR2677836}, which appear when scrutinizing the $3$-morphisms. It is not clear if these obstructions can be easily expressed in our framework. A good introduction to extension theory is \onlinecite{1711.00645} by Edie-Michell.

% This line sets the project root file.
% !TEX root = ../defects_domain_wall_structures.tex
% !TEX spellcheck = en_US
\clearpage
\section{Representation Tables} \label{sec:representation_tables}

This appendix records the irreducible representations for each annular category. In the following tables, we record the following data for each irreducible representation: 
\begin{itemize}
	\item A chosen basis for the representation.
	\item The action of a generating set of annular morphisms. 
\end{itemize}
For the bivalent vertices, we have tabulated the action by
\begin{align}
\begin{array}{c}
    \includeTikz{bivalent_action}{
    \begin{tikzpicture}[scale=.7,every node/.style={scale=.65}]
	\annparamss{}{}{}{};
	\annss{$g$}{$h$};
\end{tikzpicture}}
\end{array}.\label{eqn:biaction}
\end{align}
For the two-down-one-up trivalent vertices, we have tabulated the action by 
\begin{align}
\begin{array}{c}
  \includeTikz{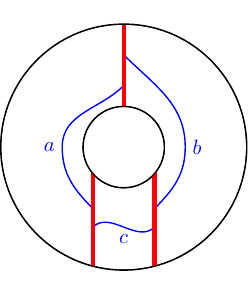}
{
	\begin{tikzpicture}[scale=1.25,every node/.style={scale=.6}]
\twodownoneuptube{$a$}{$b$}{$c$}{}{}{}
	\end{tikzpicture}	
}
\end{array}.\label{eqn:triactionA}
\end{align}
For the one-down-two-up trivalent vertices, we have tabulated the action by
\begin{align}
%
}
\caption{Bivalent representation tables. Upper table shows the chosen basis vectors for each bimodule pair. Rows correspond to lower bimod., columns label upper bimod. Lower table records the action of the annulus Eqn.~\ref{eqn:biaction} on the basis.}\label{tab:bivertextable}
\end{table}

\begin{table}
	\setlength{\tabheight}{.96\textheight}
	\resizebox*{!}{\tabheight}{
		\begin{tabular}{!{\vrule width 1pt}>{\columncolor[gray]{.9}[\tabcolsep]}c!{\vrule width 1pt}c!{\color[gray]{.8}\vrule}c!{\color[gray]{.8}\vrule}c!{\vrule width 1pt}}
			\toprule[1pt]
			&Decomposition&Basis vectors&Action\\
			\toprule[1pt]
			&$T\otimes_{\vvec{\ZZ{p}}}T$&
			$\begin{array}{c}\includeTikz{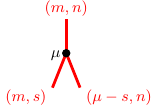}{
				\begin{tikzpicture}[scale=.7,every node/.style={scale=.5}]
				\ssstrivertex{(m,s)}{(\mu -s,n)}{(m,n)}{\mu};
				\end{tikzpicture}
			}\end{array}$&
			$\begin{array}{c}\includeTikz{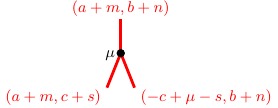}{
				\begin{tikzpicture}[scale=.7,every node/.style={scale=.5}]
				\ssstrivertex{(a+m,c+s)}{(-c+\mu -s,b+n)}{(a+m,b+n)}{\mu};
				\end{tikzpicture}
			}\end{array}$\\
			\greycline{2-4}
			&$T\otimes_{\vvec{\ZZ{p}}}L$&
			$\begin{array}{c}\includeTikz{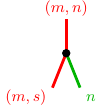}{
				\begin{tikzpicture}[scale=.7,every node/.style={scale=.5}]
				\ststrivertex{(m,s)}{n}{(m,n)}{};
				\end{tikzpicture}
			}\end{array}$&
			$\begin{array}{c}\includeTikz{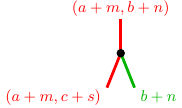}{
				\begin{tikzpicture}[scale=.7,every node/.style={scale=.5}]
				\ststrivertex{(a+m,c+s)}{b+n}{(a+m,b+n)}{};
				\end{tikzpicture}
			}\end{array}$\\
			\greycline{2-4}
			&$T\otimes_{\vvec{\ZZ{p}}}X_l$&
			$\begin{array}{c}\includeTikz{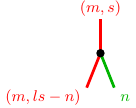}{
				\begin{tikzpicture}[scale=.7,every node/.style={scale=.5}]
				\ststrivertex{(m,l s-n)}{n}{(m,s)}{};
				\end{tikzpicture}
			}\end{array}$&
			$\begin{array}{c}\includeTikz{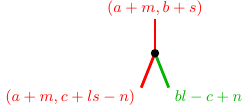}{
				\begin{tikzpicture}[scale=.7,every node/.style={scale=.5}]
				\ststrivertex{(a+m,c+l s-n)}{b l-c+n}{(a+m,b+s)}{};
				\end{tikzpicture}
			}\end{array}$\\
\end{tabular}
		\begin{tabular}{!{\vrule width 1pt}>{\columncolor[gray]{.9}[\tabcolsep]}c!{\vrule width 1pt}c!{\color[gray]{.8}\vrule}c!{\color[gray]{.8}\vrule}c!{\vrule width 1pt}}
			\toprule[1pt]
			&Decomposition&Basis vectors&Action\\
			\toprule[1pt]
			&$T\otimes_{\vvec{\ZZ{p}}}R$&
			$\begin{array}{c}\includeTikz{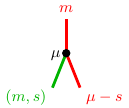}{
				\begin{tikzpicture}[scale=.7,every node/.style={scale=.5}]
				\tsstrivertex{(m,s)}{\mu -s}{m}{\mu};
				\end{tikzpicture}
			}\end{array}$&
			$\begin{array}{c}\includeTikz{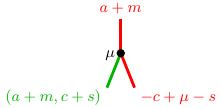}{
				\begin{tikzpicture}[scale=.7,every node/.style={scale=.5}]
				\tsstrivertex{(a+m,c+s)}{-c+\mu -s}{a+m}{\mu};
				\end{tikzpicture}
			}\end{array}$\\
			\greycline{2-4}
			&$T\otimes_{\vvec{\ZZ{p}}}F_0$&
			$\begin{array}{c}\includeTikz{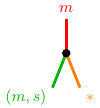}{
				\begin{tikzpicture}[scale=.7,every node/.style={scale=.5}]
				\rsttrivertex{(m,s)}{*}{m}{};
				\end{tikzpicture}
			}\end{array}$&
			$\begin{array}{c}\includeTikz{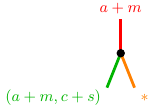}{
				\begin{tikzpicture}[scale=.7,every node/.style={scale=.5}]
				\rsttrivertex{(a+m,c+s)}{*}{a+m}{};
				\end{tikzpicture}
			}\end{array}$\\
			\greycline{2-4}
			&$T\otimes_{\vvec{\ZZ{p}}}F_r$&
			$\begin{array}{c}\includeTikz{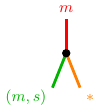}{
				\begin{tikzpicture}[scale=.7,every node/.style={scale=.5}]
				\rsttrivertex{(m,s)}{*}{m}{};
				\end{tikzpicture}
			}\end{array}$&
			$\omega^{b r s}\begin{array}{c}\includeTikz{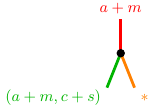}{
				\begin{tikzpicture}[scale=.7,every node/.style={scale=.5}]
				\rsttrivertex{(a+m,c+s)}{*}{a+m}{};
				\end{tikzpicture}
			}\end{array}$\\
\end{tabular}
	}
	\caption{2:1 Trivalent basis vectors and annular action. Grey column denotes the top bimodule, left column indicates the lower left$\otimes_{\vvec{\ZZ{p}}}$right bimodules, central column records the chosen basis and right column records the action of Eqn.~\ref{eqn:triactionA}.}\label{tab:trivertextable}
\end{table}

\begin{table}
	\setlength{\tabheight}{.96\textheight}
	\resizebox*{!}{\tabheight}{
		\begin{tabular}{!{\vrule width 1pt}>{\columncolor[gray]{.9}[\tabcolsep]}c!{\vrule width 1pt}c!{\color[gray]{.8}\vrule}c!{\color[gray]{.8}\vrule}c!{\vrule width 1pt}}
			\toprule[1pt]
			&Decomposition&Basis vectors&Action\\
			\toprule[1pt]
			&$T\otimes_{\vvec{\ZZ{p}}}T$&
			$\begin{array}{c}\includeTikz{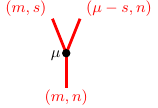}{
				\begin{tikzpicture}[scale=.7,every node/.style={scale=.5}]
				\ssstrivertexdual{(m,s)}{(\mu -s,n)}{(m,n)}{\mu};
				\end{tikzpicture}
			}\end{array}$&
			$\begin{array}{c}\includeTikz{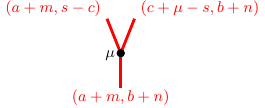}{
				\begin{tikzpicture}[scale=.7,every node/.style={scale=.5}]
				\ssstrivertexdual{(a+m,s-c)}{(c+\mu -s,b+n)}{(a+m,b+n)}{\mu};
				\end{tikzpicture}
			}\end{array}$\\
			\greycline{2-4}
			&$T\otimes_{\vvec{\ZZ{p}}}L$&
			$\begin{array}{c}\includeTikz{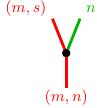}{
				\begin{tikzpicture}[scale=.7,every node/.style={scale=.5}]
				\ststrivertexdual{(m,s)}{n}{(m,n)}{};
				\end{tikzpicture}
			}\end{array}$&
			$\begin{array}{c}\includeTikz{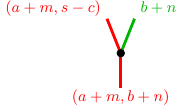}{
				\begin{tikzpicture}[scale=.7,every node/.style={scale=.5}]
				\ststrivertexdual{(a+m,s-c)}{b+n}{(a+m,b+n)}{};
				\end{tikzpicture}
			}\end{array}$\\
			\greycline{2-4}
			&$T\otimes_{\vvec{\ZZ{p}}}X_l$&
			$\begin{array}{c}\includeTikz{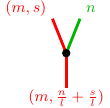}{
				\begin{tikzpicture}[scale=.7,every node/.style={scale=.5}]
				\ststrivertexdual{(m,s)}{n}{(m,\frac{n}{l}+\frac{s}{l})}{};
				\end{tikzpicture}
			}\end{array}$&
			$\begin{array}{c}\includeTikz{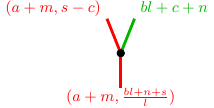}{
				\begin{tikzpicture}[scale=.7,every node/.style={scale=.5}]
				\ststrivertexdual{(a+m,s-c)}{b l+c+n}{(a+m,\frac{b l+n+s}{l})}{};
				\end{tikzpicture}
			}\end{array}$\\
\end{tabular}
		\begin{tabular}{!{\vrule width 1pt}>{\columncolor[gray]{.9}[\tabcolsep]}c!{\vrule width 1pt}c!{\color[gray]{.8}\vrule}c!{\color[gray]{.8}\vrule}c!{\vrule width 1pt}}
			\toprule[1pt]
			&Decomposition&Basis vectors&Action\\
			\toprule[1pt]
			&$T\otimes_{\vvec{\ZZ{p}}}R$&
			$\begin{array}{c}\includeTikz{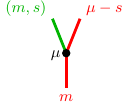}{
				\begin{tikzpicture}[scale=.7,every node/.style={scale=.5}]
				\tsstrivertexdual{(m,s)}{\mu -s}{m}{\mu};
				\end{tikzpicture}
			}\end{array}$&
			$\begin{array}{c}\includeTikz{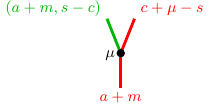}{
				\begin{tikzpicture}[scale=.7,every node/.style={scale=.5}]
				\tsstrivertexdual{(a+m,s-c)}{c+\mu -s}{a+m}{\mu};
				\end{tikzpicture}
			}\end{array}$\\
			\greycline{2-4}
			&$T\otimes_{\vvec{\ZZ{p}}}F_0$&
			$\begin{array}{c}\includeTikz{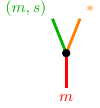}{
				\begin{tikzpicture}[scale=.7,every node/.style={scale=.5}]
				\rsttrivertexdual{(m,s)}{*}{m}{};
				\end{tikzpicture}
			}\end{array}$&
			$\begin{array}{c}\includeTikz{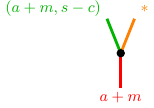}{
				\begin{tikzpicture}[scale=.7,every node/.style={scale=.5}]
				\rsttrivertexdual{(a+m,s-c)}{*}{a+m}{};
				\end{tikzpicture}
			}\end{array}$\\
			\greycline{2-4}
			&$T\otimes_{\vvec{\ZZ{p}}}F_r$&
			$\begin{array}{c}\includeTikz{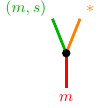}{
				\begin{tikzpicture}[scale=.7,every node/.style={scale=.5}]
				\rsttrivertexdual{(m,s)}{*}{m}{};
				\end{tikzpicture}
			}\end{array}$&
			$\omega^{-b r s}\begin{array}{c}\includeTikz{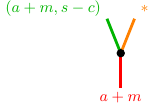}{
				\begin{tikzpicture}[scale=.7,every node/.style={scale=.5}]
				\rsttrivertexdual{(a+m,s-c)}{*}{a+m}{};
				\end{tikzpicture}
			}\end{array}$\\
\end{tabular}
	}
	\caption{1:2 Trivalent basis vectors and annular action. Grey column denotes the bottom bimodule, left column indicates the upper left$\otimes_{\vvec{\ZZ{p}}}$right bimodules, central column records the chosen basis and right column records the action of Eqn.~\ref{eqn:triactionB}.}\label{tab:dualtrivertextable}
\end{table}
% This line sets the project root file.
% !TEX root = ../defects_domain_wall_structures.tex
\clearpage
\section{Bimodule Associator Tables} \label{sec:bimodule_associator_tables}

\begin{align}
\begin{array}{c}
\includeTikz{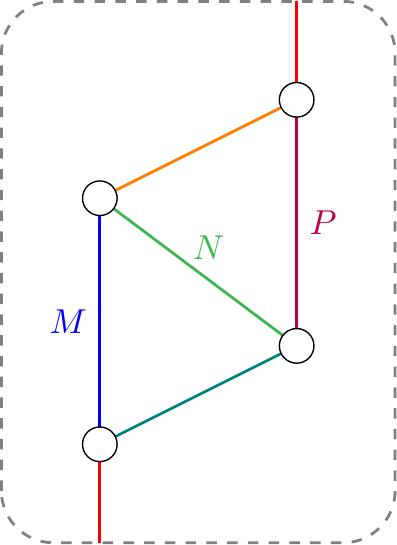}
{
	\begin{tikzpicture}[scale=0.5]
	\genericbimoduleassociator{$M$}{$N$}{$P$}
	\end{tikzpicture}
}
\end{array}\label{eqn:bimod_associator_def_appendix}
\end{align}
\setlength{\tabwidth}{.49\textwidth}
\setlength{\tabheight}{.2\textwidth}
\begin{table}[h!] 
	\renewcommand{\arraystretch}{2}
	\resizebox*{\tabwidth}{\tabheight}{
		\begin{tabular}{!{\vrule width 1pt}>{\columncolor[gray]{.9}[\tabcolsep]}c!{\vrule width 1pt}c !{\color[gray]{.8}\vrule} c !{\color[gray]{.8}\vrule} c !{\color[gray]{.8}\vrule} c!{\vrule width 1pt}c !{\color[gray]{.8}\vrule} c!{\vrule width 1pt}}
			\toprule[1pt]
			\rowcolor[gray]{.9}[\tabcolsep]\cellcolor{red!20} $\,\,T\,\,$ &$T$&$L$&$R$&$F_0$&$X_b$&$F_y$\\
			\toprule[1pt]
			$ T $&$ \delta_{\mu_0}^{\nu_0} \delta_{\mu_1}^{\nu_1} \defect{T}{T}{0}{0}{\mu_0,\nu_0,\mu_1,\nu_1} $&$ \delta_{\mu_1}^{\nu_1} \defect{T}{T}{0}{0}{\mu_1,\nu_1} $&$ \delta_{\mu_0}^{\nu_0} \defect{T}{T}{0}{0}{\mu_0,\nu_0} $&$ \defect{T}{T}{0}{0}{} $&$ \delta_{\mu_0}^{\nu_1 b} \defect{T}{T}{0}{0}{\mu_0,\nu_1} $&$ \defect{T}{T}{0}{0}{} $\\\greycline{2-7}
			$ L $&$ \delta_{\mu_0}^{\nu_0} \delta_{\mu_1}^{\nu_1} \defect{L}{L}{0}{0}{\mu_0,\nu_0,\mu_1,\nu_1} $&$ \delta_{\mu_1}^{\nu_1} \defect{L}{L}{0}{0}{\mu_1,\nu_1} $&$ \delta_{\mu_0}^{\nu_0} \defect{L}{L}{0}{0}{\mu_0,\nu_0} $&$ \defect{L}{L}{0}{0}{} $&$ \delta_{\mu_0}^{\nu_1 b} \defect{L}{L}{0}{0}{\mu_0,\nu_1} $&$ \defect{L}{L}{0}{0}{} $\\\greycline{2-7}
			$ R $&$ \delta_{\mu_1}^{\nu_1} \defect{T}{T}{0}{0}{\mu_1,\nu_1} $&$ \delta_{\mu_0}^{\nu_0} \delta_{\mu_1}^{\nu_1} \defect{T}{T}{0}{0}{\mu_0,\nu_0,\mu_1,\nu_1} $&$ \defect{T}{T}{0}{0}{} $&$ \delta_{\mu_0}^{\nu_0} \defect{T}{T}{0}{0}{\mu_0,\nu_0} $&$ \defect{T}{T}{0}{0}{} $&$ \delta_{\mu_0}^{\nu_1 y} \defect{T}{T}{0}{0}{\mu_0,\nu_1} $\\\greycline{2-7}
			$ F_{0} $&$ \delta_{\mu_1}^{\nu_1} \defect{L}{L}{0}{0}{\mu_1,\nu_1} $&$ \delta_{\mu_0}^{\nu_0} \delta_{\mu_1}^{\nu_1} \defect{L}{L}{0}{0}{\mu_0,\nu_0,\mu_1,\nu_1} $&$ \defect{L}{L}{0}{0}{} $&$ \delta_{\mu_0}^{\nu_0} \defect{L}{L}{0}{0}{\mu_0,\nu_0} $&$ \defect{L}{L}{0}{0}{} $&$ \delta_{\mu_0}^{\nu_1 y} \defect{L}{L}{0}{0}{\mu_0,\nu_1} $\\
\end{tabular}
	}
	\resizebox*{\tabwidth}{\tabheight}{
		\begin{tabular}{!{\vrule width 1pt}>{\columncolor[gray]{.9}[\tabcolsep]}c!{\vrule width 1pt}c !{\color[gray]{.8}\vrule} c !{\color[gray]{.8}\vrule} c !{\color[gray]{.8}\vrule} c!{\vrule width 1pt}c !{\color[gray]{.8}\vrule} c!{\vrule width 1pt}}
			\toprule[1pt]
			\rowcolor[gray]{.9}[\tabcolsep]\cellcolor{red!20} $\,\,L\,\,$ &$T$&$L$&$R$&$F_0$&$X_b$&$F_y$\\
			\toprule[1pt]
			$T $&$ \delta_{\mu_0}^{\nu_0} \defect{T}{T}{0}{0}{\mu_0,\nu_0} $&$ \defect{T}{T}{0}{0}{} $&$ \delta_{\mu_0}^{\nu_0} \delta_{\mu_1}^{\nu_1} \defect{T}{T}{0}{0}{\mu_0,\nu_0,\mu_1,\nu_1} $&$ \delta_{\mu_1}^{\nu_1} \defect{T}{T}{0}{0}{\mu_1,\nu_1} $&$ \defect{T}{T}{0}{0}{} $&$ \delta_{y \mu_0}^{\nu_1} \defect{T}{T}{0}{0}{\mu_0,\nu_1} $\\\greycline{2-7}
			$L $&$ \delta_{\mu_0}^{\nu_0} \defect{L}{L}{0}{0}{\mu_0,\nu_0} $&$ \defect{L}{L}{0}{0}{} $&$ \delta_{\mu_0}^{\nu_0} \delta_{\mu_1}^{\nu_1} \defect{L}{L}{0}{0}{\mu_0,\nu_0,\mu_1,\nu_1} $&$ \delta_{\mu_1}^{\nu_1} \defect{L}{L}{0}{0}{\mu_1,\nu_1} $&$ \defect{L}{L}{0}{0}{} $&$ \delta_{y \mu_0}^{\nu_1} \defect{L}{L}{0}{0}{\mu_0,\nu_1} $\\\greycline{2-7}
			$R $&$ \defect{T}{T}{0}{0}{} $&$ \delta_{\mu_0}^{\nu_0} \defect{T}{T}{0}{0}{\mu_0,\nu_0} $&$ \delta_{\mu_1}^{\nu_1} \defect{T}{T}{0}{0}{\mu_1,\nu_1} $&$ \delta_{\mu_0}^{\nu_0} \delta_{\mu_1}^{\nu_1} \defect{T}{T}{0}{0}{\mu_0,\nu_0,\mu_1,\nu_1} $&$ \delta_{b \mu_0}^{\nu_1} \defect{T}{T}{0}{0}{\mu_0,\nu_1} $&$ \defect{T}{T}{0}{0}{} $\\\greycline{2-7}
			$F_{0} $&$ \defect{L}{L}{0}{0}{} $&$ \delta_{\mu_0}^{\nu_0} \defect{L}{L}{0}{0}{\mu_0,\nu_0} $&$ \delta_{\mu_1}^{\nu_1} \defect{L}{L}{0}{0}{\mu_1,\nu_1} $&$ \delta_{\mu_0}^{\nu_0} \delta_{\mu_1}^{\nu_1} \defect{L}{L}{0}{0}{\mu_0,\nu_0,\mu_1,\nu_1} $&$ \delta_{b \mu_0}^{\nu_1} \defect{L}{L}{0}{0}{\mu_0,\nu_1} $&$ \defect{L}{L}{0}{0}{} $\\
\end{tabular}
	}
	\resizebox*{\tabwidth}{\tabheight}{
		\begin{tabular}{!{\vrule width 1pt}>{\columncolor[gray]{.9}[\tabcolsep]}c!{\vrule width 1pt}c !{\color[gray]{.8}\vrule} c !{\color[gray]{.8}\vrule} c !{\color[gray]{.8}\vrule} c!{\vrule width 1pt}c !{\color[gray]{.8}\vrule} c!{\vrule width 1pt}}
			\toprule[1pt]
			\rowcolor[gray]{.9}[\tabcolsep]\cellcolor{red!20} $\,\,R\,\,$ &$T$&$L$&$R$&$F_0$&$X_b$&$F_y$\\
			\toprule[1pt]
			$T $&$ \delta_{\mu_0}^{\nu_0} \delta_{\mu_1}^{\nu_1} \defect{R}{R}{0}{0}{\mu_0,\nu_0,\mu_1,\nu_1} $&$ \delta_{\mu_1}^{\nu_1} \defect{R}{R}{0}{0}{\mu_1,\nu_1} $&$ \delta_{\mu_0}^{\nu_0} \defect{R}{R}{0}{0}{\mu_0,\nu_0} $&$ \defect{R}{R}{0}{0}{} $&$ \delta_{\mu_0}^{\nu_1 b} \defect{R}{R}{0}{0}{\mu_0,\nu_1} $&$ \defect{R}{R}{0}{0}{} $\\\greycline{2-7}
			$	L $&$ \delta_{\mu_0}^{\nu_0} \delta_{\mu_1}^{\nu_1} \defect{F_{0}}{F_{0}}{0}{0}{\mu_0,\nu_0,\mu_1,\nu_1} $&$ \delta_{\mu_1}^{\nu_1} \defect{F_{0}}{F_{0}}{0}{0}{\mu_1,\nu_1} $&$ \delta_{\mu_0}^{\nu_0} \defect{F_{0}}{F_{0}}{0}{0}{\mu_0,\nu_0} $&$ \defect{F_{0}}{F_{0}}{0}{0}{} $&$ \delta_{\mu_0}^{\nu_1 b} \defect{F_{0}}{F_{0}}{0}{0}{\mu_0,\nu_1} $&$ \defect{F_{0}}{F_{0}}{0}{0}{} $\\\greycline{2-7}
			$	R $&$ \delta_{\mu_1}^{\nu_1} \defect{R}{R}{0}{0}{\mu_1,\nu_1} $&$ \delta_{\mu_0}^{\nu_0} \delta_{\mu_1}^{\nu_1} \defect{R}{R}{0}{0}{\mu_0,\nu_0,\mu_1,\nu_1} $&$ \defect{R}{R}{0}{0}{} $&$ \delta_{\mu_0}^{\nu_0} \defect{R}{R}{0}{0}{\mu_0,\nu_0} $&$ \defect{R}{R}{0}{0}{} $&$ \delta_{\mu_0}^{\nu_1 y} \defect{R}{R}{0}{0}{\mu_0,\nu_1} $\\\greycline{2-7}
			$	F_{0} $&$ \delta_{\mu_1}^{\nu_1} \defect{F_{0}}{F_{0}}{0}{0}{\mu_1,\nu_1} $&$ \delta_{\mu_0}^{\nu_0} \delta_{\mu_1}^{\nu_1} \defect{F_{0}}{F_{0}}{0}{0}{\mu_0,\nu_0,\mu_1,\nu_1} $&$ \defect{F_{0}}{F_{0}}{0}{0}{} $&$ \delta_{\mu_0}^{\nu_0} \defect{F_{0}}{F_{0}}{0}{0}{\mu_0,\nu_0} $&$ \defect{F_{0}}{F_{0}}{0}{0}{} $&$ \delta_{\mu_0}^{\nu_1 y} \defect{F_{0}}{F_{0}}{0}{0}{\mu_0,\nu_1} $\\
\end{tabular}
	}
	\resizebox*{\tabwidth}{\tabheight}{
		\begin{tabular}{!{\vrule width 1pt}>{\columncolor[gray]{.9}[\tabcolsep]}c!{\vrule width 1pt}c !{\color[gray]{.8}\vrule} c !{\color[gray]{.8}\vrule} c !{\color[gray]{.8}\vrule} c!{\vrule width 1pt}c !{\color[gray]{.8}\vrule} c!{\vrule width 1pt}}
			\toprule[1pt]
			\rowcolor[gray]{.9}[\tabcolsep]\cellcolor{red!20} $\,\,F_0\,\,$ &$T$&$L$&$R$&$F_0$&$X_b$&$F_y$\\
			\toprule[1pt]
			$	T $&$ \delta_{\mu_0}^{\nu_0} \defect{R}{R}{0}{0}{\mu_0,\nu_0} $&$ \defect{R}{R}{0}{0}{} $&$ \delta_{\mu_0}^{\nu_0} \delta_{\mu_1}^{\nu_1} \defect{R}{R}{0}{0}{\mu_0,\nu_0,\mu_1,\nu_1} $&$ \delta_{\mu_1}^{\nu_1} \defect{R}{R}{0}{0}{\mu_1,\nu_1} $&$ \defect{R}{R}{0}{0}{} $&$ \delta_{y \mu_0}^{\nu_1} \defect{R}{R}{0}{0}{\mu_0,\nu_1} $\\\greycline{2-7}
			$	L $&$ \delta_{\mu_0}^{\nu_0} \defect{F_{0}}{F_{0}}{0}{0}{\mu_0,\nu_0} $&$ \defect{F_{0}}{F_{0}}{0}{0}{} $&$ \delta_{\mu_0}^{\nu_0} \delta_{\mu_1}^{\nu_1} \defect{F_{0}}{F_{0}}{0}{0}{\mu_0,\nu_0,\mu_1,\nu_1} $&$ \delta_{\mu_1}^{\nu_1} \defect{F_{0}}{F_{0}}{0}{0}{\mu_1,\nu_1} $&$ \defect{F_{0}}{F_{0}}{0}{0}{} $&$ \delta_{y \mu_0}^{\nu_1} \defect{F_{0}}{F_{0}}{0}{0}{\mu_0,\nu_1} $\\\greycline{2-7}
			$	R $&$ \defect{R}{R}{0}{0}{} $&$ \delta_{\mu_0}^{\nu_0} \defect{R}{R}{0}{0}{\mu_0,\nu_0} $&$ \delta_{\mu_1}^{\nu_1} \defect{R}{R}{0}{0}{\mu_1,\nu_1} $&$ \delta_{\mu_0}^{\nu_0} \delta_{\mu_1}^{\nu_1} \defect{R}{R}{0}{0}{\mu_0,\nu_0,\mu_1,\nu_1} $&$ \delta_{b \mu_0}^{\nu_1} \defect{R}{R}{0}{0}{\mu_0,\nu_1} $&$ \defect{R}{R}{0}{0}{} $\\\greycline{2-7}
			$	F_{0} $&$ \defect{F_{0}}{F_{0}}{0}{0}{} $&$ \delta_{\mu_0}^{\nu_0} \defect{F_{0}}{F_{0}}{0}{0}{\mu_0,\nu_0} $&$ \delta_{\mu_1}^{\nu_1} \defect{F_{0}}{F_{0}}{0}{0}{\mu_1,\nu_1} $&$ \delta_{\mu_0}^{\nu_0} \delta_{\mu_1}^{\nu_1} \defect{F_{0}}{F_{0}}{0}{0}{\mu_0,\nu_0,\mu_1,\nu_1} $&$ \delta_{b \mu_0}^{\nu_1} \defect{F_{0}}{F_{0}}{0}{0}{\mu_0,\nu_1} $&$
			\defect{F_{0}}{F_{0}}{0}{0}{} $\\
\end{tabular}
	}
	\setlength{\tabheight}{.23\textwidth}
	\resizebox*{\tabwidth}{\tabheight}{
		\begin{tabular}{!{\vrule width 1pt}>{\columncolor[gray]{.9}[\tabcolsep]}c!{\vrule width 1pt}c !{\color[gray]{.8}\vrule} c !{\color[gray]{.8}\vrule} c !{\color[gray]{.8}\vrule} c!{\vrule width 1pt}c !{\color[gray]{.8}\vrule} c!{\vrule width 1pt}}
			\toprule[1pt]
			\rowcolor[gray]{.9}[\tabcolsep]\cellcolor{red!20} $\,\,F_z\,\,$ &$T$&$L$&$R$&$F_0$&$X_b$&$F_y$\\
			\toprule[1pt]
			$	T $&$ \delta_{\mu_0}^{\nu_0} \defect{R}{R}{0}{0}{\mu_0,\nu_0} $&$ \defect{R}{R}{0}{0}{} $&$ \delta_{\mu_0}^{\nu_0} \defect{T}{T}{0}{0}{\mu_0,\nu_0} $&$ \defect{T}{T}{0}{0}{} $&$ \defect{R}{R}{0}{0}{} $&$ \defect{T}{T}{0}{0}{} $\\\greycline{2-7}
			$	L $&$ \delta_{\mu_0}^{\nu_0} \defect{F_{0}}{F_{0}}{0}{0}{\mu_0,\nu_0} $&$ \defect{F_{0}}{F_{0}}{0}{0}{} $&$ \delta_{\mu_0}^{\nu_0} \defect{L}{L}{0}{0}{\mu_0,\nu_0} $&$ \defect{L}{L}{0}{0}{} $&$ \defect{F_{0}}{F_{0}}{0}{0}{} $&$ \defect{L}{L}{0}{0}{} $\\\greycline{2-7}
			$	R $&$ \defect{R}{R}{0}{0}{} $&$ \delta_{\mu_0}^{\nu_0} \defect{R}{R}{0}{0}{\mu_0,\nu_0} $&$ \defect{T}{T}{0}{0}{} $&$ \delta_{\mu_0}^{\nu_0} \defect{T}{T}{0}{0}{\mu_0,\nu_0} $&$ \defect{T}{T}{0}{0}{} $&$ \defect{R}{R}{0}{0}{} $\\\greycline{2-7}
			$	F_{0} $&$ \defect{F_{0}}{F_{0}}{0}{0}{} $&$ \delta_{\mu_0}^{\nu_0} \defect{F_{0}}{F_{0}}{0}{0}{\mu_0,\nu_0} $&$ \defect{L}{L}{0}{0}{} $&$ \delta_{\mu_0}^{\nu_0} \defect{L}{L}{0}{0}{\mu_0,\nu_0} $&$ \defect{L}{L}{0}{0}{} $&$ \defect{F_{0}}{F_{0}}{0}{0}{} $\\\toprule[1pt]
			$	X_{a} $&$ \defect{R}{R}{0}{0}{} $&$ \defect{F_{0}}{F_{0}}{0}{0}{} $&$ \defect{T}{T}{0}{0}{} $&$ \defect{L}{L}{0}{0}{} $&\cellcolor{blue!10}$ \defect{F_{(a b)^{-1}z}}{F_{(a b)^{-1}z}}{0}{0}{} $&\cellcolor{blue!10}$ \defect{X_{y^{-1}(a z)}}{X_{y^{-1}(a z)}}{0}{0}{} $\\\greycline{2-7}
			$	F_{x} $&$ \defect{F_{0}}{F_{0}}{0}{0}{} $&$ \defect{R}{R}{0}{0}{} $&$ \defect{L}{L}{0}{0}{} $&$ \defect{T}{T}{0}{0}{} $&\cellcolor{blue!10}$ \defect{X_{(b x)^{-1}z}}{X_{(b x)^{-1}z}}{0}{0}{} $&\cellcolor{blue!10}$ \defect{F_{y^{-1}(x z)}}{F_{y^{-1}(x z)}}{0}{0}{} $\\
\end{tabular}
	}
	\resizebox*{\tabwidth}{\tabheight}{
		\begin{tabular}{!{\vrule width 1pt}>{\columncolor[gray]{.9}[\tabcolsep]}c!{\vrule width 1pt}c !{\color[gray]{.8}\vrule} c !{\color[gray]{.8}\vrule} c !{\color[gray]{.8}\vrule} c!{\vrule width 1pt}c !{\color[gray]{.8}\vrule} c!{\vrule width 1pt}}
			\toprule[1pt]
			\rowcolor[gray]{.9}[\tabcolsep]\cellcolor{red!20} $\,\,X_c\,\,$ &
	$T$                                                         &
	$L$                                                         &
	$R$                                                                 &
	$F_0$                                                               &
	$X_b$                                                        &
	$F_y$                                                                  \\
			\toprule[1pt]
			$	T $                                                           &
	$ \delta_{\mu_0}^{\nu_0} \defect{T}{T}{0}{0}{\mu_0,\nu_0} $ &
	$ \defect{T}{T}{0}{0}{} $                                   &
	$ \delta_{\mu_0}^{\nu_0} \defect{R}{R}{0}{0}{\mu_0,\nu_0} $         &
	$ \defect{R}{R}{0}{0}{} $                                           &
	$ \defect{T}{T}{0}{0}{} $                                    &
	$ \defect{R}{R}{0}{0}{} $                                              \\
			\greycline{2-7}
			$	L $                                                           &
	$ \delta_{\mu_0}^{\nu_0} \defect{L}{L}{0}{0}{\mu_0,\nu_0} $ &
	$ \defect{L}{L}{0}{0}{} $                                   &
	$ \delta_{\mu_0}^{\nu_0} \defect{F_{0}}{F_{0}}{0}{0}{\mu_0,\nu_0} $ &
	$ \defect{F_{0}}{F_{0}}{0}{0}{} $                                   &
	$ \defect{L}{L}{0}{0}{} $                                    &
	$ \defect{F_{0}}{F_{0}}{0}{0}{} $                                      \\
\end{tabular}
	}
	\caption{Bimodule associator tables. All associators are trivial. Rows, columns, table label (top left) label $M,N,P$ respectively in Eqn.~\ref{eqn:bimod_associator_def_appendix}. Highlighted cells correspond to associators of invertible bimodules, the only previously known associators for this model.}\label{tab:bimod_associator_table}
	\renewcommand{\arraystretch}{1}
\end{table}
\end{document}